\documentclass[11pt]{article}
\usepackage{amsmath}
\usepackage{amssymb}
\usepackage{amsfonts}
\usepackage{mathrsfs}
\usepackage{shortvrb,psfrag}
\usepackage{enumerate}
\usepackage{color}
\usepackage[dvips]{graphicx}

\allowdisplaybreaks

\let\pa=\partial

\let\e=\varepsilon

\let\f=\frac

\def\R{\Bbb R}

\def\no{\noindent}

\def\endproof{\hphantom{MM}\hfill\llap{$\square$}\goodbreak}

\newcommand{\beq}{\begin{equation}}
\newcommand{\eeq}{\end{equation}}
\newcommand{\ben}{\begin{eqnarray}}
\newcommand{\een}{\end{eqnarray}}
\newcommand{\beno}{\begin{eqnarray*}}
\newcommand{\eeno}{\end{eqnarray*}}

\newtheorem{Theorem}{Theorem}[section]

\newtheorem{Proposition}[Theorem]{Proposition}
\newtheorem{Lemma}[Theorem]{Lemma}

\newtheorem{Remark}[Theorem]{Remark}

\begin{document}
\title{\bf Zero-viscosity limit of the Navier-Stokes equations with the Navier friction boundary condition}

\author{Tao Tao$^\dag$, Wendong Wang$^\ddag$ and Zhifei Zhang$^\sharp$
\\[2mm]
{\small $ ^\dag$ School of  Mathematical Sciences, Shandong  University, Jinan 250100 , China}\\
{\small E-mail:  taotao@sdu.edu.cn}\\[2mm]
{\small $ ^\ddag$ School of Mathematical Sciences, Dalian University of Technology, Dalian 116024,  China}\\
{\small E-mail: wendong@dlut.edu.cn}\\[2mm]
{\small $ ^\sharp$ School of  Mathematical Sciences, Peking University, Beijing 100871, China}\\
{\small E-mail: zfzhang@math.pku.edu.cn}}
\maketitle

\begin{abstract}
In this paper, we consider the zero-viscosity limit  of the Navier-Stokes equations in a half space with the Navier friction boundary condition
\beno
(\beta u^{\varepsilon}-\varepsilon^{\gamma}\partial_y u^{\varepsilon})|_{y=0}=0,
\eeno
where $\beta$ is a constant and $\gamma\in (0,1]$. In the case of $\gamma=1$, the convergence to the Euler equations and  the Prandtl equation with the Robin boundary condition is justified for the analytic data.  In the case of $\gamma\in (0,1)$, the convergence to the Euler equations and  the linearized Prandtl equation is justified for the data in the Gevrey  class $\f 1 {\gamma}$.
\end{abstract}

%
%

\setcounter{equation}{0}
\section{Introduction}
In this paper, we consider the Navier-Stokes equations in the half space $\R^2_+$:
\begin{align} \label{eq:NS}(\rm NS)\,\, \left\{
\begin{aligned}
&\partial_t u^{\varepsilon}+ u^{\varepsilon}\partial_x u^{\varepsilon}+ v^{\varepsilon}\partial_y u^{\varepsilon}+\partial_x p^{\varepsilon}=\varepsilon^2\triangle u^{\varepsilon},\\
&\partial_t v^{\varepsilon}+ u^{\varepsilon}\partial_x v^{\varepsilon}+ v^{\varepsilon}\partial_y v^{\varepsilon}+\partial_y p^{\varepsilon}=\varepsilon^2\triangle v^{\varepsilon},\\
&\partial_x u^{\varepsilon}+\partial_y v^{\varepsilon}=0,
\end{aligned}
\right.
\end{align}
with the Navier friction boundary condition
\ben\label{eq:NF-bc}
v^{\varepsilon}|_{y=0}=0,\quad \eta u^{\varepsilon}+\partial_y u^{\varepsilon}|_{y=0}=0,
\een
which was first proposed by Navier and derived for gases by Maxwell.
Here $(u^\e,v^\e)$ is the velocity field, $p^\e$ is a scalar pressure, $\e^2$ is the viscosity coefficient and $\eta$ is the slip length.

 As $v^\e|_{y=0}=0$, the boundary condition
$\eta u^{\varepsilon}+\partial_y u^{\varepsilon}|_{y=0}=0$ can be written as
\beno
\eta u^{\varepsilon}+(\partial_y u^{\varepsilon}+\pa_xv^\e)\big|_{y=0}=0,\eeno
which means that the rate of strain on the boundary is proportional to the tangential slip velocity. For $\eta=0$, the Navier friction boundary condition is just the Navier slip boundary condition; Letting  $\eta\rightarrow +\infty$,  we derive the non-slip boundary condition.
As mentioned in \cite{QWS}, the slip length should depend on the viscosity
coefficient. For simplicity, we  consider the slip length $\eta$ of the form
$\eta=-\beta\e^{-\gamma}$ for some $\gamma\ge 0$,  where $\beta$ is a constant independent of $\e$. Thus, the Navier friction boundary condition \eqref{eq:NF-bc} is reduced to
\ben\label{bc:Navier}
v^\e|_{y=0}=0,\quad \beta u^\e-\e^\gamma\pa_yv^\e|_{y=0}=0.
\een

We are concerned with the behaviour of the solution as $\e\to 0$, i.e., the zero viscosity limit. Formally, as $\varepsilon\rightarrow 0$, the solution of \eqref{eq:NS} will be approximated by the Euler equations:
\begin{align}\label{e:Euler equation}
(\rm E)
\left\{
\begin{aligned}
&\partial_t u^{e}+ u^{e}\partial_x u^{e}+ v^{e}\partial_y u^{e}+\partial_x p^{e}=0,\\
&\partial_t v^{e}+ u^{e}\partial_x v^{e}+ v^{e}\partial_y v^{e}+\partial_y p^{e}=0,\\
&\partial_x u^{e}+\partial_y v^{e}=0,\\
&v^{e}(t,x,0)=0.
\end{aligned}
\right.
\end{align}
In the absence of physical boundary, it has been proved that the Navier-Stokes equations indeed converge to the Euler equations in various functional settings \cite{K1, CW, AD, Mas}. However, in the presence of physical boundaries, this is a challenging problem due to the possible formation of boundary layer. As the boundary layer is weak for the Navier slip boundary condition, the limit from the Navier-Stokes equations to the Euler equations has been justified in \cite{IP, XX, IS, MR,WXZ}.
While, the boundary layer is strong for the non-slip boundary condition. Prandtl developed the boundary layer theory in \cite{Pr}, where the Prandtl equation was derived by the following asymptotic boundary layer expansion:
\ben\label{eq:Prandtl-expan}
\begin{split}
&u^\e(t,x,y)=u^e(t,x,y)+u^p\big(t,x,\f y \e\big)+O(\e),\\
&v^\e(t,x,y)=v^e(t,x,y)+\e v^p\big(t,x,\f y \e\big)+O(\e),
\end{split}
\een
where $(u^p,v^p)$ satisfies a Prandtl type equation.
Roughly speaking, it was expected that the Navier-Stokes equations when $\e$ is small can be approximated by the Euler equations away from the boundary, and by the Prandtl equation near the boundary.

To justify the Prandtl boundary layer expansion \eqref{eq:Prandtl-expan}, one of the key steps is to establish the well-posedness of the Prandtl equation. Up to now,  the well-posedness of Prandtl equation was only established in some special functional space. Under a monotonic assumption on the velocity of the outflow, Oleinik and Samokhin \cite{OS} established the local existence and uniqueness of classical solutions in 2-D. The global existence of weak solution was established for the favorable pressure by Xin and Zhang \cite{XZ}. Recently, Alexandre et. al. \cite{AWXY} and Masmoudi and Wong \cite{MW} independently proved the local well-posedness in Sobolev space by a direct energy method. Sammartino and Caflisch \cite{SC1} obtained the local existence and uniqueness of analytic solution for full analytic data, see \cite{LCS, ZZ} for tangential analytic data. On the other hand, Gerard-Varet and Dormy \cite{GD} proved the ill-posedness in Sobolev space for the linearized Prandtl equation around non-monotonic shear flows.

Although one has a good understanding for the Prandtl equation, there are few results on the rigorous verification of the Prandtl boundary layer expansion. In \cite{SC2}, Sammartino and Caflisch achieved this in the analytic setting, and Wang, Wang and Zhang \cite{WWZ} present a new proof based on a direct energy method. Recently, Maekawa \cite{Ma}  justified the Prandtl boundary layer expansion for the initial vorticity supported away from the boundary. Fei, Tao and Zhang  \cite{FTZ} generalized Maekawa's result to $\R^3_+$ by using a direct energy method. Very recently,  Gerard-Varet, Maekawa and Masmoudi \cite{GMM}  proved  its stability of a class of shear flows of Prandtl type in the Gevrey class. Let us also mention some conditional convergence results \cite{TW, TW2, Ke1, Ke2, CKV} initiated by Kato \cite{K2} and some convergence results for special flows \cite{LMT, MT1, MT2}. We refer to the review paper \cite{MM} for more results.

For the Navier friction boundary condition \eqref{bc:Navier}, Wang, Wang and Xin \cite{WWX} formally derived the boundary layer expansion by using the multi-scale analysis. The asymptotic behaviour of the solution  depends on the slip length. For $\gamma>1$, the behaviour is the same as the case of the non-slip boundary condition; For $\gamma=1$, the boundary layer equation is the Prandtl equation with the Robin boundary condition;  For $\gamma\in (0,1)$, the boundary layer equation is the linearized Prandtl-type equation.\smallskip

The goal of this paper is to justify the boundary layer expansion derived by Wang, Wang and Xin in the Gevrey class with the regularity exponent depending on $\gamma$.


\section{Error equations, functional spaces and main result}
\setcounter{equation}{0}

In this section, we will derive the error equations, introduce the Gevrey functional spaces and state our main result.

\subsection{The error equations}

Let us first assume that the approximate solution $(u^a,v^a,p^a)$ satisfies
\begin{align} \label{eq:approximate 1}\,\, \left\{
\begin{aligned}
&\partial_t u^a+ u^a\partial_x u^a+ (v^a-\varepsilon^2f(t,x)e^{-y})\partial_y u^a+\partial_x p^a-\varepsilon^2\triangle u^a=-R_1,\\
&\partial_t v^a+ u^a\partial_x v^a+ (v^a-\varepsilon^2f(t,x)e^{-y})\partial_y v^a+\partial_y p^a-\varepsilon^2\triangle v^a=-R_2,\\
&\partial_x u^a+\partial_y v^a=0,\\
&v^a(t,x,y)|_{y=0}=\varepsilon^2f(t,x),\\
&\partial_y u^a(t,x,y)|_{y=0}=\beta\varepsilon^{-\gamma} u^a(t,x,0)-\varepsilon g_0(t,x),\\
&(u^a,v^a)(t,x,y)|_{t=0}=(u_0(x,y),v_0(x,y)),
\end{aligned}
\right. \end{align}
where $g_0(t,x)$ satisfies $g_0(0,x)=0$ and $(R_1, R_2)$ are remainders which are small in some functional space.

We introduce the error between the solution and the approximate solution
\beno
u=u^{\varepsilon}-u^a,\quad v=v^{\varepsilon}-v^a,\quad p=p^{\varepsilon}-p^a.
\eeno
Thanks to (\ref{eq:NS}), (\ref{bc:Navier}) and $(\ref{eq:approximate 1})$, we deduce that $(u,v,p)$ satisfies
\begin{align} \label{eq:error}\,\, \left\{
\begin{aligned}
&\partial_t u+ u\partial_x u^\varepsilon+u^a\partial_x u+(v+\varepsilon^2f(t,x)e^{-y})\partial_y u^\varepsilon+(v^a-\varepsilon^2f(t,x)e^{-y})\partial_y u\\
&\qquad+\partial_x p-\varepsilon^2\triangle u=R_1,\\
&\partial_t v+ u\partial_x v^\varepsilon+u^a\partial_x v+(v+\varepsilon^2f(t,x)e^{-y})\partial_y v^\varepsilon+(v^a-\varepsilon^2f(t,x)e^{-y})\partial_y v\\
&\qquad+\partial_y p-\varepsilon^2\triangle v=R_2,\\
&\partial_x u+\partial_y v=0,
\end{aligned}
\right. \end{align}
with the boundary condition
\ben\label{eq:boundary condition of u}
v|_{y=0}=-\varepsilon^2f(t,x),\quad \partial_yu|_{y=0}=\beta\varepsilon^{-\gamma} u(t,x,0)+\varepsilon g_0(t,x),
\een
and zero initial data.

For simplicity, we write
\beno
&&U^a=(u^a,v^a),\quad U=(u,v),\quad \tilde{U}^a=(u^a,v^a-\varepsilon^2f(t,x)e^{-y})=(u^a,\tilde{v}^a),\\
&& \tilde{U}=(u,v+\varepsilon^2f(t,x)e^{-y})=(u,\tilde{v}),\quad \tilde{R}=(R_1,R_2).
\eeno
Then we have
\begin{align} \label{eq:error2}\,\, \left\{
\begin{aligned}
&\partial_t U+ \tilde{U}\cdot\nabla(U+U^a)+\tilde{U}^a\cdot \nabla U
+\nabla p-\varepsilon^2\triangle U=\tilde{R},\\
&{\rm div} U=0,\\
&\partial_yu|_{y=0}=\beta\varepsilon^{-\gamma} u(t,x,0)+\varepsilon g_0(t,x),\\
&v|_{y=0}=-\varepsilon^2 f ,\\
&U|_{t=0}=0.
\end{aligned}
\right. \end{align}

We introduce the vorticity $w, w^a$ defined by
\beno
\omega=\partial_yu-\partial_x v,\quad \omega^a=\partial_yu^a-\partial_x v^a.
\eeno
It is easy to see from (\ref{eq:error2})  that
\begin{align} \label{eq:vorticity}\,\, \left\{
\begin{aligned}
&\partial_t \omega-\varepsilon^2\triangle \omega+\tilde{U}\cdot\nabla(\omega+\omega^a)+\tilde{U}^a\cdot \nabla \omega\\
&=\partial_yR_1-\partial_xR_2+\varepsilon^2 f e^{-y}\partial_yu^a+\varepsilon^2 \partial_xf e^{-y}\partial_yv^a,\\
&w|_{y=0}=\beta\varepsilon^{-\gamma} u(t,x,0)+\varepsilon g_0(t,x)+\varepsilon^2 \partial_xf ,\\
&w|_{t=0}=0.
\end{aligned}
\right. \end{align}

Let $g(t,x)=g_0(t,x)+\varepsilon \partial_xf(t,x)$ and
\ben\label{eq:relationship of vorticity}
\eta=\omega-\beta \varepsilon^{-\gamma} u-\varepsilon g(t,x) e^{-y},\quad \eta^a=\omega^a-\beta\varepsilon^{-\gamma} u^a.
\een
It follows from (\ref{eq:error}), (\ref{eq:boundary condition of u}) and (\ref{eq:vorticity}) that
\begin{eqnarray}\label{eq:modified vorticity equation}
\left \{
\begin {array}{ll}
\partial_t\eta-\varepsilon^2\triangle \eta+u^a\partial_x\eta+\widetilde{v}^a\partial_y\eta+u\partial_x\eta^a+\widetilde{v}\partial_y\eta^a+u\partial_x\eta\\[5pt]
\qquad \qquad +\widetilde{v}\partial_y\eta-\beta\frac{\partial_xp}{\varepsilon^\gamma}=
\partial_yR_1-\partial_xR_2+\varepsilon^2fe^{-y}\partial_yu^a+\varepsilon^2\partial_xfe^{-y}\partial_yv^a-\frac{\beta R_1}{\varepsilon^\gamma}+h,\\[5pt]
\eta(0,x,z)=0,\\[5pt]
\eta|_{y=0}=0,
\end{array}
\right.
\end{eqnarray}
where
$$
h=-\varepsilon e^{-y}\Big[\partial_tg+u^a\partial_xg-\widetilde{v}^ag+u\partial_xg-\widetilde{v}g
-\varepsilon^2g-\varepsilon^2\partial_{xx}g\Big].
$$

\subsection{The Gevrey functional spaces}

As in \cite{MR,MW}, we introduce the conormal operator $Z=\psi(y)\partial_y$, where
$\psi(y)$ is a smooth function defined by
\begin{align*} \psi(y)=\left\{
\begin{aligned}
&\delta y,\quad {\rm for}\quad y\leq 1,\\
&\frac{\delta y}{1+y},\quad {\rm for}\quad y\geq 2,
\end{aligned}
\right. \end{align*}
where $\delta>0$ is to be decided later. We denote
\beno
Z^k=\psi(y)^k\partial_y^k,\quad \tilde{Z}^k=(\delta z)^k\partial_z^k.
\eeno
The conormal Sobolev space $\overline{H}^s(\R^2_+)$ for $s\in N$ is defined by
\beno
\overline{H}^s(\R^2_+)=\Big\{u; \|u\|_{\overline{H}^s}\doteq\sum_{k+\ell\leq s}\|Z^k\partial_x^{\ell} u\|_{L^2_{x,y}(\R^2_+)}<\infty\Big\}.
\eeno
We also denote $\|f\|_{L^p_{x,y}(\R^2_+)}$ by  $\|f\|_{p}$ for $1\leq p\leq \infty.$

Furthermore, for $\alpha\in N^2$ with  $\alpha=(\alpha_1,\alpha_2)$, let $\partial^\alpha=\partial^{\alpha_1}_xZ^{\alpha_2}$ and
$$\|\partial^ku\|_\infty=\Big(\sum_{|\alpha|=k}\|\partial^\alpha u\|^2_\infty\Big)^{\frac12},\quad k\in {\rm N}_+.$$
For $k\in {\rm N_+}$ and $\gamma\in (0,1]$, we define the conormal Gevrey norm as
\beno
\|u\|_{X^k}^2&=&\sum_{m=k}^{\infty}\frac{\rho(t)^{2(m-k)}}{((m-k)!)^{\frac{2}{\gamma}}}\sum_{|\alpha|=m}\big\|\partial^\alpha u\big\|^2_{2}+\|u\|^2_2,\\
\|u\|_{Y^k}^2&=&\sum_{m=k+1}^{\infty}\frac{\rho(t)^{2(m-k)}(m-k)}{((m-k)!)^{\frac{2}{\gamma}}}\sum_{|\alpha|=m}\big\|\partial^\alpha u\big\|^2_{2}+\|u\|^2_2.
\eeno
Here  $\rho(t)=\rho-\lambda t$ and  $\rho(t)\in [1,2].$ Especially, it is reduced to the analytic semi-norm as in \cite{WWZ}  when $\gamma=1.$ For simplicity, we set
$$|u|^2_{j,k}=\frac{\rho(t)^{2(j-k)}}{((j-k)!)^{\frac{2}{\gamma}}}\sum_{|\alpha|=j}\big\|\partial^\alpha u\big\|^2_{2}.$$
Then we have
\beno
\|u\|_{X^k}^2=\sum_{m=k}^{\infty}|u|^2_{m,k}+\|u\|_2^2,\quad
\|u\|_{Y^k}^2=\sum_{m=k+1}^{\infty}(m-k)|u|^2_{m,k}+\|u\|_2^2.
\eeno
We also introduce the Gevrey norm
\beno
\|u\|^2_{X_e^k}&=&\sum_{m=k}^\infty\frac{\rho(t)^{2(m-k)}}{((m-k)!)^{\frac{2}{\gamma}}}\sum_{|\alpha| =m}\big\|\partial_{x,y}^\alpha u\big\|^2_2+\|u\|^2_2,\\
\|u\|^2_{X_x^k}&=&\sum_{m=k}^\infty\frac{\rho(t)^{2(m-k)}}{((m-k)!)^{\frac{2}{\gamma}}}\big\|\partial_{x}^m u\big\|^2_2+ \|u\|^2_2\quad {\rm for}\quad u=u(x).
\eeno

Let us say $\beta\leq\alpha$ in  ${\rm N}^2$ if $\alpha=(\alpha_1,\alpha_2)$, $\beta=(\beta_1,\beta_2)$ satisfy $\beta_1\leq \alpha_1, \beta_2\leq\alpha_2.$ We denote $C_\alpha^\beta=C_{\alpha_1}^{\beta_1}C_{\alpha_2}^{\beta_2}$.
Let us state a useful lemma, which will be used frequently.
\begin{Lemma}\label{lem:identity}
\begin{itemize}

\item[(a)]
Let $\{x_\alpha\}_{ \alpha\in N^2}$ and $\{y_\beta\}_{ \beta\in N^2}$ be real numbers. Then we have
$$\sum_{|\alpha|=m}\sum_{|\beta|=j, \beta\leq\alpha}x_\beta y_{\alpha-\beta}=\Big(\sum_{|\alpha|=j}x_\alpha\Big)\Big(\sum_{|\beta|=m-j}y_\beta\Big).$$

\item[(b)] For $\alpha,\beta\in N^2$, $|\alpha|=m$ and $j\leq m$, there holds
$$\sum_{|\beta|=j, \beta\leq \alpha}C_\alpha^\beta=C_m^j.$$
\end{itemize}
\end{Lemma}

Part (a) can be found in \cite{KV}, and (b) can be obtained by computing the coefficient of $x^j$ in the binomial expansion of $(1+x)^{\alpha_1}(1+x)^{\alpha_2}$ and $(1+x)^m$, where $\alpha=(\alpha_1, \alpha_2)$ and $|\alpha|=m$. Here we omit the details.

\subsection{Main result}

Let us first introduce some assumptions (H) on  approximate solutions and remainders.\smallskip

\begin{itemize}

\item[(H1)] Formulation and uniform bounds of approximate solutions: let $z=\frac{y}{\varepsilon}$ and
\beno
&&u^{a}(t,x,y)= u^e(t,x,y)+\varepsilon^{1-\gamma} u^p(t,x,z),\\
&&v^{a}(t,x,y)= v^e(t,x,y)+\varepsilon^{2-\gamma} v^p(t,x,z),
\eeno
and
\textcolor[rgb]{0.00,0.00,0.00}{ there exist $a_0>0$ and }$T_a> 0$ such that for any $t\in [0,T_a]$,
\ben\label{e:uniform boundness for approximate solution}
&&\sum_{m=3}^{\infty}\frac{\rho(t)^{2(m-3)}}{((m-3)!)^{\frac{2}{\gamma}}}\sum_{m-3\leq |\alpha|\leq m+6}\big\|\partial_{x,y}^\alpha(u^{e},v^{e})(t,\cdot)\big\|^2_2\leq C_0,\nonumber\\
&&\sum_{m=3}^{\infty}\frac{\rho(t)^{2(m-3)}}{((m-3)!)^{\frac{2}{\gamma}}}\sum_{m-3\leq |\alpha|\leq m+6}\sum_{k=0}^2\|e^{a_0z^2}\widetilde{\partial}^\alpha\partial_{z}^k(u^{p},v^{p})(t,\cdot)\big\|_{2}\leq C_0,
\een
where
$
\widetilde{\partial}^\alpha=\partial_x^{\alpha_1}\widetilde{Z}^{\alpha_2},\widetilde{Z}^k=(\delta z)^k\partial_z^k.
$

\item[(H2)] Uniform bounds of the remainders:
 for $T_a> 0$ as in (H1), there holds
\ben\label{e:uniform boundness for error}
\big\|(R_1, R_2)(t,\cdot)\big\|^2_{X^3}\leq C_0\varepsilon^4,\quad \big\|\nabla(R_1, R_2)(t,\cdot)\big\|^2_{X^2}\leq C_0\varepsilon^2,
\een
 for any $t\in [0,T_a]$.

\item[(H3)]Formulation and uniform bounds of   ($f, g_0$): there exists $\overline{f}(t,x)$ such that
$f(t,x)=\partial_x\overline{f}$ with
\ben\label{eq:boundness of f}
\|\partial_tf(t,\cdot)\|_{X^3_x}+\|f(t,\cdot)\|_{X^5_x}+\|g_0(t,\cdot)\|_{X^5_x}+\|\partial_tg_0(t,\cdot)\|_{X^3_x}+\|\partial_t\bar{f}(t,\cdot)\|_{L^2_x}\leq C_0
\een
 for any $t\in [0,T_a]$.
\end{itemize}

Our main result is stated as follows.

\begin{Theorem}\label{thm:main}
Assume (H1)-(H3). There exist $T> 0$ and $C_0> 0$ independent of $\varepsilon$ such that for any sufficiently small $\varepsilon> 0$, the error equation (\ref{eq:error2}) admits a unique solution $(u, v)(t,\cdot)\in X^3$, which satisfies
\begin{align}
\sup_{0\leq t\leq T}\big\|(u,v)(t,\cdot)\big\|_{X^3}\leq C_0\varepsilon^2,\quad \sup_{0\leq t\leq T}\big\|w(t,\cdot)\big\|_{X^2}\leq C_0\varepsilon.\nonumber
\end{align}
In particular, we have
\beno
\sup_{0\leq t\leq T}\big\|(u,v)(t,\cdot)\big\|_{L^2\cap L^\infty(R^2_+)}\leq C_0\varepsilon.
\eeno
\end{Theorem}

\begin{Remark}
For every rational number $\gamma\in (0,1]$, we can construct the approximate solution of $(u^a, v^a)$ of (\ref{eq:NS}) by using the matched asymptotic expansion, which satisfies the assumption (H1)-(H3) if the initial data $(u_0,v_0)\in X_e^{k(\gamma)}$ for some $k(\gamma)\in {\rm N_+}$ depending on $\gamma$. Thus, for any rational number $\gamma\in (0,1]$ and $(u_0,v_0)\in X_e^{k(\gamma)}$, there exist $T> 0$ and $C_0> 0$ independent of $\varepsilon$ such that the Navier-Stokes equations (\ref{eq:NS}) admit a unique solution $(u^\varepsilon, v^\varepsilon)$ on $[0,T]$, which satisfies
\begin{align}
\sup_{0\leq t\leq T}\big\|(u^\varepsilon,v^\varepsilon)-(u^a,v^a)\big\|_{X^3}\leq C_0\varepsilon^2,\quad \sup_{0\leq t\leq T}\big\|w-w^a\big\|_{X^2}\leq C_0\varepsilon.\nonumber
\end{align}
In particular,
\beno
\sup_{0\leq t\leq T}\big\|(u^\varepsilon,v^\varepsilon)-(u^a,v^a)\big\|_{L^2\cap L^\infty(R^2_+)}\leq C_0\varepsilon,
\eeno
where $(u^a, v^a)$ has the form as (H1), and $(u^e, v^e)$ is the sum of solutions of the Euler equations and linearized Euler equations, $(u^p, v^p)$ is the sum of solutions of the \textcolor[rgb]{0.00,0.00,0.00}{  linearized Prandtl-type equations}(for $\gamma=1$, it is the sum of solutions of nonlinear Prandtl equation with Robin boundary condition and some linearized Prandtl equation).
\end{Remark}

\begin{Remark}
When $\gamma=0$, \eqref{eq:NS} is just the Navier-Stokes equations with the Navier-slip boundary condition. In this case, the boundary layer is weak so that the convergence from the Navier-Stokes equations to the Euler equations can be proved in Sobolev space. In fact, we can give an explanation from the viewpoint of asymptotic expansion. Observe that the approximate solutions $u^a(t,x,y)=u^e(t,x,y)+\varepsilon u^p(t,x,z), v^a(t,x,y)=v^e(t,x,y)+\varepsilon^2 v^p(t,x,z)$ for $\gamma=0$, and so $\partial_yu^a(t,x,y)=\partial_yu^e(t,x,y)+ \partial_zu^p(t,x,z), \partial_yv^a(t,x,y)=\partial_yv^e(t,x,y)+\varepsilon \partial_zu^p(t,x,z)$ are not singular terms, hence the error equation (\ref{eq:error2}) is not a singular equation and we can make the energy estimate in conormal Sobolev space.
\end{Remark}


\subsection{Proof of main result}

The key point is to prove that the solution of error equation (\ref{eq:error2}) is uniformly bounded in the suitable functional spaces.
However, from the assumption (H1) and the error equation (\ref{eq:error2}), we note that there are some singular terms such as
\beno
\tilde{v}\partial_y(\varepsilon^{1-\gamma}u^p(t,x,\frac{y}{\varepsilon}))=\frac{\tilde{v}}{\varepsilon^\gamma}\partial_zu^p(t,x,\frac {y}{\varepsilon}).
\eeno
This singular term is handled as follows: near boundary, using $\tilde{v}|_{y=0}=0$, it can be transformed to a term which loses $``\gamma"$-order derivative
\beno
\frac{\tilde{v}}{\varepsilon^\gamma}\partial_zu^p(t,x,\frac {y}{\varepsilon})=\frac{\tilde{v}}{y^\gamma}z^\gamma\partial_zu^p(t,x,\frac {y}{\varepsilon})\sim y^{1-\gamma}\partial_y\tilde{v}\partial_zu^p(t,x,\frac {y}{\varepsilon});
\eeno
Away from the boundary, this term is good and has a decay factor $y^{-\gamma}$. Hence, it is natural to work in the Gevrey  setting and we estimate this term near boundary and away from boundary in a different way respectively. On the other hand, if we directly take $\partial_y$ derivative to the error equation, the Prandtl part $u^p,v^p$ of the approximate solution will give rise to a bad factor $\varepsilon^{-1}$. To avoid this singularity, we use the conormal derivative motivated by \cite{MR, WWZ}. However, we can not obtain a control in the normal derivative near the boundary by using the conormal derivative. Motivated by \cite{Ma, WWZ}, we will use the vorticity equation to gain one order normal derivative. Finally, we complete our argument by combining the velocity estimate and vorticity estimate together.
\medskip

Let us complete the proof of main result by admitting Proposition \ref{prop:Velocity estimate in Gevrey class}, Proposition \ref{prop:Velocity estimate in Sobolev space}, Proposition \ref{prop:Vorticity estimate in Gevrey class} and Proposition \ref{prop:Vorticity estimate in Sobolev space}, which will be proved in the later sections.\smallskip

We introduce the energy functional
\begin{align}
E(t):=&\varepsilon^{-2}\Big(\big\|U\big\|^2_{X^3}+\varepsilon^4\Big)+\big\|\eta\big\|^2_{X^2}, \quad F(t):=\varepsilon^{-2}\big\|U\big\|^2_{Y^3}+\big\|\eta\big\|^2_{Y^2},\nonumber\\
G(t):=&\big\|\nabla U\big\|^2_{X^3}+\varepsilon^2\big\|\nabla \eta\big\|^2_{X^2}.\nonumber
\end{align}
Using that the fact $\gamma\leq 1$ and (\ref{eq:relationship of vorticity}),  we have
 \beno
\|\omega\big\|^2_{Y^2}\leq C_0\left(\|\eta\big\|^2_{Y^2}+\varepsilon^{-2}\|u\big\|^2_{Y^2}+\varepsilon^{2}\|ge^{-y}\big\|^2_{Y^2}\right).
\eeno
By integration by parts, we have
\beno
\|\psi(y)\partial_y u\|_2\leq C_0\left(\|u\|_2+\|\psi(y)^2\partial_y^2 u\|_2\right),
\eeno
and it follows from (\ref{eq:boundness of f}) that for $\delta<\frac14$
\beno
\|ge^{-y}\big\|^2_{Y^2}(t)&=&\sum_{m=3}^{\infty}\frac{\rho(t)^{2(m-2)}(m-2)}{((m-2)!)^{\frac{2}{\gamma}}}\sum_{|\alpha|=m}\big\|\partial^\alpha (ge^{-y})\big\|^2_{2}+\|ge^{-y}\|^2_2\\
&\leq& C_0\sum_{\alpha_1=0}^2\|\partial_x^{\alpha_1}g\|_2^2\delta^{-2\alpha_1}\sum_{m=3}^\infty\frac{\rho(t)^{2(m-2)}(m-2)}{((m-2)!)^{\frac{2}{\gamma}}}\delta^{2m}\\
&&+C_0\sum_{\alpha_1=3}^\infty\|\partial_x^{\alpha_1} g\|_2^2\delta^{-2\alpha_1}\sum_{m=\alpha_1}^\infty\frac{\rho(t)^{2(m-2)}(m-2)}{((m-2)!)^{\frac{2}{\gamma}}}\delta^{2m}+C_0\\
&\leq&C_0\|g\|_{X_x^3}(t)+C_0\leq C_0,
\eeno
hence
\beno
\|\omega\big\|^2_{Y^2}\leq
C_0(F(t)+\varepsilon^2),
\eeno
and similarly
\beno
 \|\omega\big\|^2_{X^2}\leq C_0E(t).
 \eeno
Thus, Proposition \ref{prop:Velocity estimate in Gevrey class} and Proposition \ref{prop:Velocity estimate in Sobolev space} ensure that
\ben\label{eq:estimate of prop 3.1}
&&\varepsilon^{-2}\frac{1}{2}\frac{d}{dt}\big\|U\big\|^2_{X^3}+\lambda\varepsilon^{-2}\big\|U\big\|^2_{Y^3}+\frac{1}{2}\big\|\nabla U\big\|^2_{X^3}
\nonumber\\[5pt]
&&\leq C_0\varepsilon^{\frac12}E(t)^{\frac12}(E(t)+F(t))+C_0\varepsilon^{-\frac12}E(t)^{\frac34}G(t)^{\frac14}(\varepsilon^2+\varepsilon E(t)^{\frac12}+G(t)^{\frac12} )\nonumber\\[5pt]
&&\quad+C_0(E(t)+F(t))(1+E(t)^{\frac12}+\varepsilon E(t) )+C_0[E(t)G(t)+\varepsilon^2],
\een
where we used $\partial_y u=w+\partial_x v$ and
\beno
X^3\hookrightarrow X^2\hookrightarrow \overline{H}^5,\quad Y^3\hookrightarrow Y^2.
\eeno
Similarly, Proposition \ref{prop:Vorticity estimate in Gevrey class} and Proposition \ref{prop:Vorticity estimate in Sobolev space} ensure that
\ben\label{eq:estimate of prop 3.2}
&&\frac{1}{2}\frac{d}{dt}\big\|\eta\big\|^2_{X^2}+\lambda\big\|\eta\big\|^2_{Y^2}+\frac{\varepsilon^2}{2}\big\|\nabla \eta\big\|^2_{X^2}
\nonumber\\[5pt]
&&\leq C_0(E(t)+F(t))(1+ E(t))+C_0 \delta\varepsilon^{4-4\gamma}(G(t)+\varepsilon^2)\nonumber\\[5pt]
&&\quad+C_0E(t)G(t)^{\frac12}+C_0\varepsilon^{\frac12}E(t)^{\frac54}G(t)^{\frac14}+\textcolor[rgb]{0.00,0.00,0.00}{ C_0 E(t)G(t)}+\frac{1}{10}G(t).
\een

Combining (\ref{eq:estimate of prop 3.1})-(\ref{eq:estimate of prop 3.2}) and taking $\delta$ small enough,  we arrive at
\begin{align}
\frac{d}{dt}E(t)+\lambda F(t)+\frac{1}{2}G(t)\leq C_0E(t)(1+E(t)+G(t))+C_0(1+E(t))F(t)+C_0\frac{E^3(t)}{\varepsilon^2}.\nonumber
\end{align}
Then the standard continuity argument yields that there exists $T$ independent of $\varepsilon$ and $\varepsilon_0>0$ such that
\begin{align}
\sup_{0\leq t\leq T}E(t)\leq C_0\varepsilon^2\quad {\rm for \,\, all}\quad \varepsilon\in (0,\varepsilon_0).\nonumber
\end{align}
Therefore, there holds
\beno
\sup_{0\leq t\leq T}\Big(\varepsilon^{-1}\|U\big\|_{X^3}+\big\|\eta\big\|_{X^2}\Big)\leq C_0\varepsilon\quad {\rm for \,\, all}\quad \varepsilon\in (0,\varepsilon_0).
\eeno
Recalling (\ref{eq:relationship of vorticity}), we arrive at
\beno
\sup_{0\leq t\leq T}\Big(\varepsilon^{-1}\|U\big\|_{X^3}+\big\|\omega\big\|_{X^2}\Big)\leq C_0\varepsilon,
\eeno
from which and Sobolev embedding, it follows that
$$\|U\|_{L^\infty(R^2_+)}\leq C_0\varepsilon.$$
The proof is completed.
\endproof

\section{Nonlinear estimates in Gevrey type spaces}
\setcounter{equation}{0}

Our main goal is to obtain the uniform estimates of the solution for the error equations (\ref{eq:error2}) and (\ref{eq:modified vorticity equation}) in the Gevrey norms. The key point is to estimate some linear or nonlinear terms, for example $u\partial_x u^\varepsilon+u^a\partial_x u+(v+\varepsilon^2f(t,x)e^{-y})\partial_y u^\varepsilon+(v^a-\varepsilon^2f(t,x)e^{-y})\partial_y u$. Generally speaking, there are four different types:
\\1) $u\partial_x u$,  or $u^a\partial_x u$;\\
2) $(v+\varepsilon^2f(t,x)e^{-y})\partial_y u$,  or $(v^a-\varepsilon^2f(t,x)e^{-y})\partial_y u$;\\
3) $u\partial_x u^a$;\quad 4) $v\partial_y u^a$.

In this section, we will deal with these terms. For simplicity, $\langle\cdot,\cdot\rangle$ means the inner product in $L^2_{x,y}(\R^2_+)$, and  we denote
\beno
\|f\|_{X\cap Y}=\|f\|_{X}+\|f\|_Y,\quad \|(f,g)\|_{X}=\|f\|_{X}+\|g\|_{X}.
\eeno

First of all, we deal with the terms in 1), for example, $u\partial_xv$.

\begin{Lemma}\label{lem:x product estimate}
Let
\beno
A=\sum_{m=3}^{\infty}\frac{\rho(t)^{2(m-3)}}{((m-3)!)^{\frac{2}{\gamma}}}\sum_{|\alpha|=m}\big|\langle \partial^\alpha(u\partial_xv),\partial^\alpha v\rangle\big|.
\eeno
Then there holds
\begin{itemize}
\item[(a)]
\beno
 A\leq C_0\|(u, v)\|_{X^3}^{\frac12}\| (\partial_yu,\partial_yv)\|_{X^2}^{\frac12}\big(\|v\|^2_{ Y^3}+\|(u,v)\|^2_{X^3}\big);
\eeno
\item[(b)]
\beno
 A&\leq&
C_0\|u\|_{X^4}^{\frac12}\|\partial_yu\|_{X^4}^{\frac12}\|v\|^2_{X^3\cap Y^3}.
\eeno
\end{itemize}
\end{Lemma}

\no{\bf Proof.} (a) Note that the Sobolev inequality implies that
\beno
\|\partial_x u\|_{\infty}\leq C_0\|u\|_{\overline{H}^2}^{\frac12}\|\partial_yu\|_{\overline{H}^2}^{\frac12},
\eeno
thus, by integration by parts we have
\beno
\sum_{m=3}^{\infty}\frac{\rho(t)^{2(m-3)}}{((m-3)!)^{\frac{2}{\gamma}}}\sum_{|\alpha|=m}\big|\langle u\partial_x\partial^\alpha v,\partial^\alpha v\rangle\big|
\leq C_0\|u\|_{\overline{H}^2}^{\frac12}\|\partial_yu\|_{\overline{H}^2}^{\frac12}\|v\|_{X^3}^2,
\eeno
then
\beno
\sum_{m=3}^{\infty}\frac{\rho(t)^{2(m-3)}}{((m-3)!)^{\frac{2}{\gamma}}}\sum_{|\alpha|=m}\big|\langle \partial^\alpha(u\partial_xv),\partial^\alpha v\rangle\big|
\leq I+C_0\|u\|_{\overline{H}^2}^{\frac12}\|\partial_yu\|_{\overline{H}^2}^{\frac12}\|v\|_{X^3}^2,
\eeno
where
\beno
I=\sum_{m=3}^{\infty}\frac{\rho(t)^{2(m-3)}}{((m-3)!)^{\frac{2}{\gamma}}}\sum_{|\alpha|=m}\sum_{0<\beta\leq\alpha}C_\alpha^\beta\big|\langle \partial^\beta u\partial_x\partial^{\alpha-\beta}v,\partial^\alpha v\rangle\big|,
\eeno
and $I$ could be decomposed as follows
\beno
I
&\leq&\sum_{m=3}^{\infty}\frac{\rho(t)^{2(m-3)}}{((m-3)!)^{\frac{2}{\gamma}}}\sum_{|\alpha|=m}\sum_{1\leq|\beta|\leq 2,\beta\leq\alpha}C_\alpha^\beta\|\partial^\beta u\partial_x\partial^{\alpha-\beta}v\|_{2}\|\partial^\alpha v\|_{2}\\
&&+\sum_{m=5}^{\infty}\frac{\rho(t)^{2(m-3)}}{((m-3)!)^{\frac{2}{\gamma}}}\sum_{|\alpha|=m}\sum_{j=3}^{m-2}\sum_{|\beta|=j,\beta\leq\alpha}C_\alpha^\beta\|\partial^\beta u\partial_x\partial^{\alpha-\beta}v\|_{2}\|\partial^\alpha v\|_{2}\\
&&+\sum_{m=3}^{\infty}\frac{\rho(t)^{2(m-3)}}{((m-3)!)^{\frac{2}{\gamma}}}\sum_{|\alpha|=m}\sum_{m-1\leq|\beta|\leq m,\beta\leq\alpha}C_\alpha^\beta\|\partial^\beta u\partial_x\partial^{\alpha-\beta}v\|_{2}\|\partial^\alpha v\|_{2}
=\sum_{i=1}^3I_{i}
\eeno

Next we handle them term by term.\smallskip

{\bf Step 1. Estimate of $I_1$.}
Let $I_1=I_{11}+I_{12}$ according to the value of $|\beta|=1,2$.

Firstly, applying H\"{o}lder inequality twice, we get
\ben\label{eq:beta=1}
&&\sum_{m=4}^{\infty}\frac{\rho(t)^{2(m-3)}}{((m-3)!)^{\frac{2}{\gamma}}}\sum_{|\alpha|=m}\sum_{|\beta|=1,\beta\leq\alpha}C_\alpha^\beta\|\partial^\beta u\partial_x\partial^{\alpha-\beta}v\|_{2}\|\partial^\alpha v\|_{2}\nonumber\\
&&\leq\sum_{m=4}^{\infty}\frac{\rho(t)^{(m-3)}}{((m-3)!)^{\frac{1}{\gamma}}}\Big(\sum_{|\alpha|=m}\Big(\sum_{|\beta|=1,\beta\leq\alpha}C_\alpha^\beta
\|\partial^\beta u\|_\infty\|\partial_x\partial^{\alpha-\beta}v\|_{2}\Big)^2\Big)^{\frac12}|v|_{m,3}\nonumber\\
&&\leq\sum_{m=4}^{\infty}\frac{\rho(t)^{(m-3)}}{((m-3)!)^{\frac{1}{\gamma}}}\Big[\sum_{|\alpha|=m}\Big(\sum_{|\beta|=1,\beta\leq\alpha}(C_\alpha^\beta)^2
\sum_{|\beta|=1,\beta\leq\alpha}\|\partial^\beta u\|_\infty^2\|\partial_x\partial^{\alpha-\beta}v\|_{2}^2\Big)\Big]^{\frac12}|v|_{m,3}\nonumber\\
&&\leq C_0\|\partial u\|_\infty\sum_{m=4}^{\infty}m|v|_{m,3}^2\leq C_0\|\partial u\|_\infty\|v\|^2_{Y^3},
\een
where we used Lemma \ref{lem:identity} in the last step, that is, for $j=1,\cdots,m$
\ben\label{eq:C alpha beta}
\sum_{|\beta|=j,\beta\leq\alpha}(C_\alpha^\beta)^2\leq \textcolor[rgb]{0.00,0.00,0.00}{ \Big(\sum_{|\beta|=j,\beta\leq\alpha}C_\alpha^\beta\Big)^2}\leq  (C_m^j)^2,
\een
and
\ben\label{eq:u alpha beta}
\sum_{|\alpha|=m}
\sum_{|\beta|=j,\beta\leq\alpha}\|\partial^\beta f\|_\infty^2\|\partial^{\alpha-\beta}\partial_xg\|_{2}^2\leq \sum_{|\beta|=j}\|\partial^\beta f\|_\infty^2\sum_{|\alpha|=m-j+1}\|\partial^{\alpha}g\|_{2}^2.
\een
Hence, we obtain
\beno
I_{11}&\leq&  C_0\|\partial u\|_\infty\|v\|^2_{Y^3}+ \sum_{|\alpha|=3}\sum_{|\beta|=1,\beta\leq\alpha}C_\alpha^\beta\|\partial^\beta u\partial_x\partial^{\alpha-\beta}v\|_{2}\|\partial^\alpha v\|_{2}\\
&\leq& C_0\|u\|_{\overline{H}^2}^{\frac12}\|\partial_yu\|_{\overline{H}^2}^{\frac12}\big(\|v\|^2_{Y^3}+\|v\|^2_{X^3}\big). \eeno

We remark here that the technique of (\ref{eq:beta=1}) includes  H\"{o}lder inequality(twice), (\ref{eq:C alpha beta}) and (\ref{eq:u alpha beta}), which will be used frequently, and we just mention (\ref{eq:beta=1}) later.

Secondly, for the term $I_{12}$,
 similar computations as (\ref{eq:beta=1}) yield that
\beno
&&\sum_{m=5}^{\infty}\frac{\rho(t)^{2(m-3)}}{((m-3)!)^{\frac{2}{\gamma}}}\sum_{|\alpha|=m}\sum_{|\beta|=2,\beta\leq\alpha}C_\alpha^\beta\|\partial^\beta u\partial_x\partial^{\alpha-\beta}v\|_{2}\|\partial^\alpha v\|_{2}\\
&&\leq C_0\|\partial^2 u\|_\infty\sum_{m=5}^{\infty}\frac{\rho(t)^{2(m-3)}C_m^2}{((m-3)!)^{\frac{2}{\gamma}}}
\Big(\sum_{|\alpha|=m}\|\partial^\alpha v\|^2_{2}\Big)^{\frac12}\Big(\sum_{|\alpha|=m-1}\|\partial^\alpha v\|^2_{2}\Big)^{\frac12}\\
&&\leq C_0\|\partial^2 u\|_\infty\|v\|^2_{Y^3},
\eeno
where we used $\gamma\leq 1$ and  \textcolor[rgb]{0.00,0.00,0.00}{ $\frac{m}{(m-3)^{\frac{1}{\gamma}}}\leq 3$ for $m\geq 5$}. Hence,
\beno
I_{12}\leq C_0\|\partial^2 u\|_\infty\|v\|^2_{Y^3}+\sum_{m=3}^{4}\frac{\rho(t)^{2(m-3)}}{((m-3)!)^{\frac{2}{\gamma}}}\sum_{|\alpha|=m}\sum_{|\beta|=2,\beta\leq\alpha}C_\alpha^\beta\|\partial^\beta u\partial_x\partial^{\alpha-\beta}v\|_{2}\|\partial^\alpha v\|_{2},
\eeno
which can be controlled by
\beno
C_0\|u\|_{\overline{H}^3}^{\frac12}\|\partial_yu\|_{\overline{H}^3}^{\frac12}\big(\|v\|^2_{Y^3}+\|v\|^2_{X^3}\big).
\eeno

Finally, collecting the estimates of $I_{11}$ and $I_{12}$, we obtain
\beno
I_1\leq C_0\|u\|_{\overline{H}^3}^{\frac12}\|\partial_yu\|_{\overline{H}^3}^{\frac12}\big(\|v\|^2_{Y^3}+\|v\|^2_{X^3}\big).
\eeno

{\bf Step 2. Estimate of $I_2$.} Firstly, according to the different values of $|\beta|$,  $I_2$ is divided into two terms
\beno
I_2&=&\sum_{m=6}^{\infty}\frac{\rho(t)^{2(m-3)}}{((m-3)!)^{\frac{2}{\gamma}}}\sum_{|\alpha|=m}\sum_{j=3}^{[\frac{m}{2}]}\sum_{|\beta|=j,\beta\leq\alpha}
C_\alpha^\beta\|\partial^\beta u\partial_x\partial^{\alpha-\beta}v\|_{2}\|\partial^\alpha v\|_{2}\\
&&+\sum_{m=5}^{\infty}\frac{\rho(t)^{2(m-3)}}{((m-3)!)^{\frac{2}{\gamma}}}\sum_{|\alpha|=m}\sum_{j=[\frac{m}{2}]+1}^{m-2}\sum_{|\beta|=j,\beta\leq\alpha}C_\alpha^\beta\|\partial^\beta u\partial_x\partial^{\alpha-\beta}v\|_{2}\|\partial^\alpha v\|_{2}=I_{21}+I_{22}.
\eeno
By Sobolev embedding, we have
\beno
\|\partial^\beta u\|_{\infty}\leq C_0\|\partial^\beta u\|^{\frac12}_2\|\partial_y \partial^\beta u\|^{\frac12}_2+C_0\|\partial_x\partial^\beta u\|^{\frac12}_2\|\partial_x\partial_y \partial^\beta u\|^{\frac12}_2,
\eeno
which gives
\beno
I_{21}
&\leq&\sum_{m=6}^{\infty}\frac{\rho(t)^{2(m-3)}}{((m-3)!)^{\frac{2}{\gamma}}}\sum_{j=3}^{[\frac{m}{2}]}\sum_{|\alpha|=m}
\sum_{|\beta|=j,\beta\leq\alpha}C_\alpha^\beta\|\partial^\beta u\|^{\frac12}_2\|\partial_y \partial^\beta u\|^{\frac12}_2\|\partial_x\partial^{\alpha-\beta}v\|_{2}\|\partial^\alpha v\|_{2}\\
&&+\sum_{m=6}^{\infty}\frac{\rho(t)^{2(m-3)}}{((m-3)!)^{\frac{2}{\gamma}}}\sum_{j=3}^{[\frac{m}{2}]}\sum_{|\alpha|=m}
\sum_{|\beta|=j,\beta\leq\alpha}C_\alpha^\beta\|\partial_x\partial^\beta u\|^{\frac12}_2\|\partial_x\partial_y \partial^\beta u\|^{\frac12}_2\|\partial_x\partial^{\alpha-\beta}v\|_{2}\|\partial^\alpha v\|_{2}.
\eeno
Obviously, it suffices to estimate the second term, since  \textcolor[rgb]{0.00,0.00,0.00}{ the order of derivative in the first term}  is more lower. Using the same argument as in (\ref{eq:beta=1}) and discrete Young inequality,
it can be controlled by
\ben\label{eq:middle term-Young inequ}
&& \sum_{m=6}^{\infty}\frac{\rho(t)^{(m-3)}}{((m-3)!)^{\frac{1}{\gamma}}}\sum_{j=3}^{[\frac{m}{2}]}
C_m^j\Big(\sum_{|\beta|=j}\|\partial^\beta\partial_x u\|_2\|\partial_y \partial^\beta\partial_x u\|_2\Big)^{\frac12}\Big(\sum_{|\alpha|=m+1-j}
\|\partial^\alpha v\|^2_{2}\Big)^{\frac12}|v|_{m,3}\nonumber\\
&&\leq C_0\sum_{m=6}^{\infty}\sum_{j=3}^{[\frac{m}{2}]}
\frac{C_m^j((m-j-2)!)^{\frac{1}{\gamma}}((j-2)!)^{\frac{1}{2\gamma}}((j-1)!)^{\frac{1}{2\gamma}}}{((m-3)!)^{\frac{1}{\gamma}}\sqrt{m-3}\sqrt{m-j-2}}\nonumber\\
&&\quad\times\frac{\rho(t)^{\frac{(j-1)}{2}}}{((j-1)!)^{\frac{1}{2\gamma}}}\Big(\sum_{|\alpha|=j+1}\|\partial_y \partial^\alpha u\|^2_2\Big)^{\frac14}|u|^{\frac12}_{j+1,3}\sqrt{m-j-2}|v|_{m+1-j,3}\sqrt{m-3}|v|_{m,3}\nonumber\\
&&\leq C_0\sum_{m=6}^{\infty}\sum_{j=3}^{[\frac{m}{2}]}\frac{1}{j}
\frac{\rho(t)^{\frac{(j-1)}{2}}}{((j-1)!)^{\frac{1}{2\gamma}}}\Big(\sum_{|\alpha|=j+1}\|\partial_y \partial^\alpha u\|^2_2\Big)^{\frac14}|u|^{\frac12}_{j+1,3}\sqrt{m-j-2}|v|_{m+1-j,3}\sqrt{m-3}|v|_{m,3}\nonumber\\
&&\leq C_0\Big(\sum_{m=6}^{\infty}\Big(\sum_{j=3}^{[\frac{m}{2}]}\frac{1}{j}
\frac{\rho(t)^{\frac{(j-1)}{2}}}{((j-1)!)^{\frac{1}{2\gamma}}}\Big(\sum_{|\alpha|=j+1}\|\partial_y \partial^\alpha u\|^2_2\Big)^{\frac14}|u|^{\frac12}_{j+1,3}\sqrt{m-j-2}|v|_{m+1-j,3}\Big)^2\Big)^{\frac12}\|v\|_{Y^3}\nonumber\\
&&\leq C_0\| u\|^{\frac12}_{X^3}\| \partial_yu\|_{X^2}^{\frac12}\|v\|^2_{Y^3},
\een
where we used i) the estimate:
$$\frac{C_m^j((m-j-2)!)^{\frac{1}{\gamma}}((j-1)!)^{\frac{1}{\gamma}}}{((m-3)!)^{\frac{1}{\gamma}}\sqrt{m-3}\sqrt{m-j-2}}\leq C_0\Big(C_{m-3}^{j-1}\Big)^{1-\frac{1}{\gamma}}j^{-1}\leq C_0j^{-1}$$
for $\gamma\leq 1$ and $j\in \{3,\cdots,[\frac{m}{2}]\}$;
ii)discrete Young inequality to estimate
\beno
&&\Big(\sum_{m=6}^{\infty}\Big(\sum_{j=3}^{[\frac{m}{2}]}\frac{1}{j}
\frac{\rho(t)^{\frac{(j-1)}{2}}}{((j-1)!)^{\frac{1}{2\gamma}}}\Big(\sum_{|\alpha|=j+1}\|\partial_y \partial^\alpha u\|^2_2\Big)^{\frac14}|u|^{\frac12}_{j+1,3}\sqrt{m-j-2}|v|_{m+1-j,3}\Big)^2\Big)^{\frac12}\\
&&\leq C\|v\|_{Y^3}\sum_{m=3}^\infty\frac{1}{m}
\frac{\rho(t)^{\frac{(m-1)}{2}}}{((m-1)!)^{\frac{1}{2\gamma}}}\Big(\sum_{|\alpha|=m+1}\|\partial_y \partial^\alpha u\|^2_2\Big)^{\frac14}|u|^{\frac12}_{m+1,3}\\
&&\leq C\|v\|_{Y^3}\|u\|_{X^3}^{\frac12}\Big(\sum_{m=3}^\infty
\frac{\rho(t)^{2(m-1)}}{((m-1)!)^{\frac{2}{\gamma}}}\sum_{|\alpha|=m+1}\|\partial_y \partial^\alpha u\|^2_2\Big)^{\frac14},
\eeno
iii) the commutator estimate: as
\beno
\|\partial_y\partial^\alpha u\|_{2}\leq \|\partial^\alpha \partial_yu\|_{2}+\|[\partial_y,\partial^\alpha] u\|_{2},
\eeno
thus,
\ben\label{e:y derivative estimate low order}
&&\sum_{m=4}^{\infty}\frac{\rho(t)^{2(m-2)}}{((m-2)!)^{\frac{2}{\gamma}}}\sum_{|\alpha|=m}\|\partial_y\partial^\alpha u\|^2_{2}\nonumber\\
&&\leq C_0\| \partial_yu\|^2_{X^2}+\sum_{m=4}^{\infty}\frac{\rho(t)^{2(m-2)}}{((m-2)!)^{\frac{2}{\gamma}}}\sum_{|\alpha|=m-1}m^2\|\partial^\alpha \partial_yu\|^2_{2}\leq C_0\| \partial_yu\|^2_{X^2}.
\een
Hence, we arrive at
$$I_{21}\leq C_0\| u\|^{\frac12}_{X^3}\| \partial_yu\|_{X^2}^{\frac12}\|v\|^2_{Y^3}.$$

Secondly, different from the estimate of $I_{21}$, for $I_{22}$, we estimate $L^{\infty}$ norm of $\partial_x\partial^{\alpha-\beta}v,$ and it can be bounded by
\beno
&&\sum_{m=5}^{\infty}\frac{\rho(t)^{2(m-3)}}{((m-3)!)^{\frac{2}{\gamma}}}\sum_{|\alpha|=m}\sum_{j=[\frac{m}{2}]+1}^{m-2}
\sum_{|\beta|=j,\beta\leq\alpha}C_\alpha^\beta\|\partial^\beta u\|_2\|\partial_x\partial^{\alpha-\beta}v\|^{\frac12}_{2}\|\partial_x\partial_y\partial^{\alpha-\beta}v\|^{\frac12}_{2}\|\partial^\alpha v\|_{2}\\
&&\quad+\sum_{m=5}^{\infty}\frac{\rho(t)^{2(m-3)}}{((m-3)!)^{\frac{2}{\gamma}}}\sum_{|\alpha|=m}\sum_{j=[\frac{m}{2}]+1}^{m-2}
\sum_{|\beta|=j,\beta\leq\alpha}C_\alpha^\beta\|\partial^\beta u\|_2\|\partial_{xx}\partial^{\alpha-\beta}v\|^{\frac12}_{2}\|\partial_{xx}\partial_y\partial^{\alpha-\beta}v\|^{\frac12}_{2}\|\partial^\alpha v\|_{2}.
\eeno
We only estimate the second term. Using similar arguments as in (\ref{eq:middle term-Young inequ}), it be controlled by
\beno
&& \sum_{m=5}^{\infty}\frac{\rho(t)^{(m-3)}}{((m-3)!)^{\frac{1}{\gamma}}}\sum_{j=[\frac{m}{2}]+1}^{m-2}
C_m^j\Big(\sum_{|\alpha|=m-j+2}\|\partial^\alpha v\|_2\|\partial_y \partial^\alpha v\|_2\Big)^{\frac12}\Big(\sum_{|\alpha|=j}
\|\partial^\alpha u\|^2_{2}\Big)^{\frac12}|v|_{m,3}\\
&&\leq C_0\|v\|^{\frac12}_{\overline{H}^4}\|\partial_yv\|^{\frac12}_{\overline{H}^4}\|u\|_{X^3} \|v\|_{X^3}\\ &&\quad+C_0\sum_{m=6}^{\infty}\sum_{j=[\frac{m}{2}]+1}^{m-2}
\frac{C_m^j((m-j-1)!)^{\frac{1}{2\gamma}}((m-j)!)^{\frac{1}{2\gamma}}((j-3)!)^{\frac{1}{\gamma}}}{((m-3)!)^{\frac{1}{\gamma}}\sqrt{m-3}}\\
&&\quad\times\frac{\rho(t)^{\frac{(m-j)}{2}}}{((m-j)!)^{\frac{1}{2\gamma}}}\Big(\sum_{|\alpha|=m-j+2}\|\partial_y \partial^\alpha v\|^2_2\Big)^{\frac14}|v|^{\frac12}_{m-j+2,3}|u|_{j,3}\sqrt{m-3}|v|_{m,3}\\
&&\leq C_0\sum_{m=6}^{\infty}\sum_{j=[\frac{m}{2}]+1}^{m-2}
\frac{(m-j)^{-1}\rho(t)^{\frac{(m-j)}{2}}}{((m-j)!)^{\frac{1}{2\gamma}}}\Big(\sum_{|\alpha|=m-j+2}\|\partial_y \partial^\alpha v\|^2_2\Big)^{\frac14}|v|^{\frac12}_{m-j+2,3}|u|_{j,3}\sqrt{m-3}|v|_{m,3}\\
&&\quad+C_0\|v\|^{\frac12}_{\overline{H}^4}\|\partial_yv\|^{\frac12}_{\overline{H}^4}\|u\|_{X^3} \|v\|_{X^3}\\
&&\leq C_0\| v\|^{\frac12}_{X^3}\| \partial_yv\|_{X^2}^{\frac12}\|u\|_{X^3}\|v\|_{X^3\cap Y^3},
\eeno
where we used
$$\frac{C_m^j((m-j-1)!)^{\frac{1}{2\gamma}}((m-j)!)^{\frac{1}{2\gamma}}((j-3)!)^{\frac{1}{\gamma}}}{((m-3)!)^{\frac{1}{\gamma}}\sqrt{m-3}}\leq C_0\Big(C_{m-3}^{j-3}\Big)^{1-\frac{1}{\gamma}}(m-j)^{-1}$$
for $j\in \{[\frac{m}{2}]+1,\cdots, m-2\}$
and (\ref{e:y derivative estimate low order}).

Therefore, we get
\beno
I_{22}\leq   C_0\| v\|^{\frac12}_{X^3}\| \partial_yv\|_{X^2}^{\frac12}\|u\|_{X^3}\|v\|_{X^3\cap Y^3}.
\eeno

Finally, collecting the estimates of $I_{21}$ and $I_{22}$, we obtain
\beno
I_2\leq  C_0\Big[\| u\|^{\frac12}_{X^3}\| \partial_yu\|_{X^2}^{\frac12}\|v\|^2_{Y^3}+
\| v\|^{\frac12}_{X^3}\| \partial_yv\|_{X^2}^{\frac12}\|u\|_{X^3}\|v\|_{X^3\cap Y^3}
\Big].
\eeno

{\bf Step 3. Estimate of $I_3$.} This is similar to $I_2$, but more easier. We rewrite $I_{31}$ as the term of $|\beta|=m-1$ in $I_3$, and $I_{32}$ for the term of $|\beta|=m.$

Again using the technique as in (\ref{eq:beta=1}), we get
\beno
&&\sum_{m=4}^{\infty}\frac{\rho(t)^{2(m-3)}}{((m-3)!)^{\frac{2}{\gamma}}}\sum_{|\alpha|=m}\sum_{|\beta|=m-1,\beta\leq\alpha}C_\alpha^\beta\|\partial^\beta u\|_2\|\partial_x\partial^{\alpha-\beta}v\|_\infty\|\partial^\alpha v\|_{2}\\
&&\leq C_0\|\partial^2v\|_\infty\sum_{m=4}^{\infty}\frac{\rho(t)^{2(m-3)}m}{((m-3)!)^{\frac{2}{\gamma}}}
\Big(\sum_{|\alpha|=m-1}\|\partial^\alpha u\|^2_2\Big)^{\frac12}\Big(\sum_{|\alpha|=m}\|\partial^\alpha v\|^2_{2}\Big)^{\frac12}\\
&&\leq C_0\|\partial^2v\|_\infty\|u\|_{X^3}\|v\|_{X^3}
\eeno
and
\beno
\sum_{|\alpha|=3}\sum_{|\beta|=2,\beta\leq\alpha}C_\alpha^\beta\|\partial^\beta u\|_2\|\partial_x\partial^{\alpha-\beta}v\|_\infty\|\partial^\alpha v\|_{2}\leq C_0\|\partial^2v\|_\infty\|\partial^2u\|_2\|v\|_{X^3},
\eeno
which give
\beno
I_{31}\leq C_0\|v\|_{\overline{H}^3}^{\frac12}\|\partial_yv\|_{\overline{H}^3}^{\frac12} \|u\|_{X^3}\|v\|_{X^3}.
\eeno

Similarly, there holds
\beno
I_{32}&\leq&\sum_{m=3}^{\infty}\frac{\rho(t)^{2(m-3)}}{((m-3)!)^{\frac{2}{\gamma}}}\sum_{|\alpha|=m}\sum_{|\beta|=m,\beta\leq\alpha}C_\alpha^\beta\|\partial^\beta u\|_2\|\partial_x\partial^{\alpha-\beta}v\|_\infty\|\partial^\alpha v\|_{2}\\
&\leq&C_0\|v\|_{\overline{H}^3}^{\frac12}\|\partial_yv\|_{\overline{H}^3}^{\frac12} \|u\|_{X^3}\|v\|_{X^3},
\eeno
which along with the estimate of $I_{31}$ implies that
\beno
I_3\leq C_0\|v\|_{\overline{H}^3}^{\frac12}\|\partial_yv\|_{\overline{H}^3}^{\frac12} \|u\|_{X^3}\|v\|_{X^3}.
\eeno

Collecting the estimates in Step 1-Step 3, the proof of the inequality (a) is completed.\smallskip

(b) This inequality is used to estimate the linear term like $u^a\partial_xv$. We use the same notations as in (a).
First of all, we know that
\beno
\sum_{m=3}^{\infty}\frac{\rho(t)^{2(m-3)}}{((m-3)!)^{\frac{2}{\gamma}}}\sum_{|\alpha|=m}\big|\langle \partial^\alpha(u\partial_xv),\partial^\alpha v\rangle\big|
\leq I+C_0\|u\|_{\overline{H}^3}^{\frac12}\|\partial_yu\|_{\overline{H}^3}^{\frac12}\|v\|_{X^3}^2
\eeno
and  $
I=\sum\limits_{i=1}^3I_{i}
$ as in (a).

The estimate of $I_1$ is as follows
\beno
I_1\leq C_0\|u\|_{\overline{H}^3}^{\frac12}\|\partial_yu\|_{\overline{H}^3}^{\frac12}\big(\|v\|^2_{Y^3}+\|v\|^2_{X^3}\big).
\eeno

{\bf  Estimate of $I_2$.} By Sobolev embedding, $I_2$ can be bounded by
\beno
&&\sum_{m=6}^{\infty}\frac{\rho(t)^{2(m-3)}}{((m-3)!)^{\frac{2}{\gamma}}}\sum_{|\alpha|=m}\sum_{j=3}^{m-3}
\sum_{|\beta|=j,\beta\leq\alpha}C_\alpha^\beta\|\partial^\beta u\|^{\frac12}_2\|\partial_y \partial^\beta u\|^{\frac12}_2\|\partial_x\partial^{\alpha-\beta}v\|_{2}\|\partial^\alpha v\|_{2}\\
&&+\sum_{m=6}^{\infty}\frac{\rho(t)^{2(m-3)}}{((m-3)!)^{\frac{2}{\gamma}}}\sum_{|\alpha|=m}\sum_{j=3}^{m-3}
\sum_{|\beta|=j,\beta\leq\alpha}C_\alpha^\beta\|\partial_x\partial^\beta u\|^{\frac12}_2\|\partial_x\partial_y \partial^\beta u\|^{\frac12}_2\|\partial_x\partial^{\alpha-\beta}v\|_{2}\|\partial^\alpha v\|_{2}\\
&&+\sum_{m=6}^{\infty}\frac{\rho(t)^{2(m-3)}}{((m-3)!)^{\frac{2}{\gamma}}}\sum_{|\alpha|=m}
\sum_{|\beta|=m-2,\beta\leq\alpha}C_\alpha^\beta\|\partial^\beta u\partial_x\partial^{\alpha-\beta}v\|_{2}\|\partial^\alpha v\|_{2}+C_0\sum_{|\beta|=3}\|\partial^\beta u\|_\infty\|v\|^2_{X^3}\\
&&=\sum_{i=1}^3I_{2i}+C_0\|u\|^{\frac12}_{\overline{H}^4}\|\partial_yu\|^{\frac12}_{\overline{H}^4}\|v\|^2_{X^3}.
\eeno
Now we estimate the term $I_{22}.$ It follows from Lemma \ref{lem:identity} and discrete convolution inequality as (\ref{eq:middle term-Young inequ}) that
\beno
I_{22}&\leq& \sum_{m=6}^{\infty}\frac{\rho(t)^{(m-3)}}{((m-3)!)^{\frac{1}{\gamma}}}\sum_{j=3}^{m-3}
C_m^j\Big(\sum_{|\beta|=j+1}\|\partial^\beta u\|_2\|\partial_y \partial^\beta u\|_2\Big)^{\frac12}\Big(\sum_{|\alpha|=m+1-j}
\|\partial^\alpha v\|^2_{2}\Big)^{\frac12}|v|_{m,3}\\
&\leq& C_0\sum_{m=6}^{\infty}\sum_{j=3}^{m-3}
\frac{C_m^j((m-j-2)!)^{\frac{1}{\gamma}}((j-3)!)^{\frac{1}{\gamma}}}{((m-3)!)^{\frac{1}{\gamma}}\sqrt{m-3}\sqrt{m-j-2}}
\frac{\rho(t)^{(j-3)}}{((j-3)!)^{\frac{1}{\gamma}}}\Big(\sum_{|\beta|=j+1}\|\partial^\beta u\|_2\|\partial_y \partial^\beta u\|_2\Big)^{\frac12}\\
&&\times\sqrt{m-j-2}
|v|_{m+1-j,3}\sqrt{m-3}|v|_{m,3}\\
&\leq& C_0\sum_{m=6}^{\infty}\sum_{j=3}^{m-3}
\frac{a_{m,j}\rho(t)^{(j-3)}}{((j-3)!)^{\frac{1}{\gamma}}}\Big(\sum_{|\beta|=j+1}\|\partial^\beta u\|_2\|\partial_y \partial^\beta u\|_2\Big)^{\frac12}\cdot\sqrt{m-j-2}
|v|_{m+1-j,3}\sqrt{m-3}|v|_{m,3}\\
&\leq& C_0\| u\|^{\frac12}_{X^4}\| \partial_yu\|_{X^4}^{\frac12}\|v\|^2_{Y^3},
\eeno
where we used
\begin{align}\label{eq:a m j}
a_{m,j}=\left\{
\begin{array}{lll}
\frac{1}{j^2},\quad j\in \Big\{3,...,\Big[\frac{m}{2}\Big]\Big\},\\
\frac{1}{(m-j)^2},\quad j\in \Big\{\Big[\frac{m}{2}\Big],...,m-3\Big\}
\end{array}
\right.
\end{align}
satisfying
\beno
\frac{C_m^j((m-j-2)!)^{\frac{1}{\gamma}}((j-3)!)^{\frac{1}{\gamma}}}{((m-3)!)^{\frac{1}{\gamma}}\sqrt{m-3}\sqrt{m-j-2}}\leq C_0a_{m,j}
\eeno
and the commutator estimate as in (\ref{e:y derivative estimate low order}).

The same argument gives
\beno
&&I_{21}\leq C_0\| u\|^{\frac12}_{X^3}\| \partial_yu\|_{X^3}^{\frac12}\|v\|^2_{Y^3}.
\eeno
Thus, we arrive at
\beno
I_2\leq C_0\| u\|_{X^4}^{\frac12}\| \partial_yu\|_{X^4}^{\frac12}\|v\|^2_{Y^3\cap X^3}
+I_{23},
\eeno
and $I_{23}$ will be estimated in the next step.

{\bf Estimate of $I_3$ and $I_{23}$.} We first decompose $I_3$ as in (a).
By similar arguments as in (\ref{eq:beta=1}), we get
\beno
&&\sum_{m=4}^{\infty}\frac{\rho(t)^{2(m-3)}}{((m-3)!)^{\frac{2}{\gamma}}}\sum_{|\alpha|=m}\sum_{|\beta|=m-1,\beta\leq\alpha}C_\alpha^\beta\|\partial^\beta u\|_{L^\infty_yL^2_x}\|\partial_x\partial^{\alpha-\beta}v\|_{L^2_yL^\infty_x}\|\partial^\alpha v\|_{2}\\
&&\leq C_0\|v\|_{\overline{H}^3}\sum_{m=4}^{\infty}\frac{\rho(t)^{2(m-3)}m}{((m-3)!)^{\frac{2}{\gamma}}}
\Big(\sum_{|\alpha|=m}\|\partial^\alpha v\|^2_2\Big)^{\frac12}\Big(\sum_{|\alpha|=m-1}\|\partial^\alpha u\|_{2}\|\partial_y\partial^\alpha u\|_{2}\Big)^{\frac12}\\
&&\leq C_0\| u\|_{X^3}^{\frac12}\| \partial_yu\|_{X^3}^{\frac12}\|v\|_{X^3}\|v\|_{\overline{H}^3}.
\eeno
Hence, we have
\beno
I_{31}&\leq& C_0\| u\|_{X^3}^{\frac12}\| \partial_yu\|_{X^3}^{\frac12}\|v\|_{X^3}\|v\|_{\overline{H}^3} +C_0\|\partial^2u\|_\infty\|\partial^2v\|_2\|v\|_{X^3}\\
&\leq&  C_0\| u\|_{X^3}^{\frac12}\| \partial_yu\|_{X^3}^{\frac12}\|v\|_{X^3}\|v\|_{\overline{H}^3}.
\eeno
The same argument gives
\beno
I_{32}&\leq&C_0\| u\|_{X^3}^{\frac12}\| \partial_yu\|_{X^3}^{\frac12}\|v\|_{X^3}\|v\|_{\overline{H}^2},\\
I_{23}&\leq&C_0\| u\|_{X^3}^{\frac12}\| \partial_yu\|_{X^3}^{\frac12}\|v\|_{X^3}\|v\|_{\overline{H}^4}.
\eeno
Finally, there holds
\beno
I_3+I_{23}\leq C_0\| u\|_{X^3}^{\frac12}\| \partial_yu\|_{X^3}^{\frac12}\|v\|_{X^3}\|v\|_{\overline{H}^4}.
\eeno
Recalling the estimates of $I_1$ and $I_2$, the proof is completed.\endproof\smallskip

\smallskip
To estimate the term in 2) like $v\partial_yu$, we need the following lemma.

\begin{Lemma}\label{e:y product estimate}
For $k=2,3$ and the suitable functions $u, v$ with $u|_{y=0}=0$, let
\beno
B_k&=&\sum_{m=k}^{\infty}\frac{\rho(t)^{2(m-k)}}{((m-k)!)^{\frac{2}{\gamma}}}\sum_{|\alpha|=m}\big|\langle \partial^\alpha(u\partial_yv),\partial^\alpha v\rangle\big|.
\eeno
There holds
\begin{itemize}
\item[(a)]
\beno
 B&\leq &C_0(\delta)\|(u,\partial_yu)\|^{\frac12}_{\overline{H}^3}\|(\partial_yu,\partial_{yy}u)\|^{\frac12}_{\overline{H}^3}
\| v\|^2_{X^k\cap Y^k}\\
&&+C_0\| u\|^{\frac12}_{X^k}\|\partial_yu\|_{X^k}^{\frac12}
\|\partial_yv\|_{X^k}\|v\|_{X^k}+C_0(\delta)\| v\|^{\frac32}_{X^k}\| \partial_yv\|_{X^k}^{\frac12}\|(u,\partial_yu)\|_{X^k}.
\eeno
\item[(b)]
\beno
B\leq C_0(\delta)\|(u,\partial_yu)\|_{X^{k+1}}^{\frac12}\|(\partial_yu,\partial_{yy}u)\|_{X^{k+1}}^{\frac12}
\| v\|^2_{X^k\cap Y^k}.
\eeno
\end{itemize}
\end{Lemma}

\no{\bf Proof.} We only give a proof for $k=3$. The estimate for $k=2$ can be obtained by the same argument.

(a) Let
\beno
&&\sum_{m=3}^{\infty}\frac{\rho(t)^{2(m-3)}}{((m-3)!)^{\frac{2}{\gamma}}}\sum_{|\alpha|=m}\big|\langle \partial^\alpha(u\partial_yv),\partial^\alpha v\rangle\big|\leq\tilde{I}+\widetilde{II},
\eeno
where
\beno
\tilde{I}&=&\sum_{m=3}^{\infty}\frac{\rho(t)^{2(m-3)}}{((m-3)!)^{\frac{2}{\gamma}}}\sum_{|\alpha|=m}\sum_{0<\beta\leq\alpha}C_\alpha^\beta\big|\langle \partial^\beta u\partial^{\alpha-\beta}\partial_yv,\partial^\alpha v\rangle\big|,\\
\widetilde{II}&=&\sum_{m=3}^{\infty}\frac{\rho(t)^{2(m-3)}}{((m-3)!)^{\frac{2}{\gamma}}}\sum_{|\alpha|=m}\big|\langle u\partial^\alpha\partial_y v,\partial^\alpha v\rangle\big|.
\eeno

We first estimate the term $\widetilde{II}$. Thanks to the definition of $\psi$, we have
\beno
\Big\|\frac{u}{\psi}\Big\|_\infty\leq C_0(\delta)(\|\partial_yu\|_\infty+\|u\|_\infty)\leq C_0(\delta)\|(u,\partial_yu)\|_{\overline{H}^1}^{\frac12}\|(\partial_yu,\partial_{yy}u)\|_{\overline{H}^1}^{\frac12}.
\eeno
Using $u|_{y=0}=0$ and integration by parts, we get
\beno
\widetilde{II}&=&\sum_{m=3}^{\infty}\frac{\rho(t)^{2(m-3)}}{((m-3)!)^{\frac{2}{\gamma}}}\sum_{|\alpha|=m}\Big[\langle u\partial_y \partial^\alpha v,\partial^\alpha v\rangle
-\langle u[\partial_y, \partial^\alpha ]v,\partial^\alpha v\rangle\Big]\\
&\leq &\|\partial_yu\|_{L^\infty}\|v\|_{X^3}^2+\Big\|\frac{u}{\psi}\Big\|_\infty\Big(\sum_{m=4}^{\infty}\frac{\rho(t)^{2(m-3)}m}{((m-3)!)^{\frac{2}{\gamma}}}
\sum_{|\alpha|=m} \|\partial^\alpha v\|^2_{2}+\sum_{|\alpha|=3}\|\partial^\alpha v\|^2_{2}\Big)\\
&\leq& C_0(\delta)\|(u,\partial_yu)\|_{\overline{H}^1}^{\frac12}\|(\partial_yu,\partial_{yy}u)\|_{\overline{H}^1}^{\frac12} \|v\|^2_{X^3\cap Y^3}.
\eeno
Then, similar to the estimate of $I$ in Lemma \ref{lem:x product estimate}, we decompose $\tilde{I}$ into three terms according to the value of $|\beta|$:
\beno
\tilde{I}
&\leq&\sum_{m=3}^{\infty}\frac{\rho(t)^{2(m-3)}}{((m-3)!)^{\frac{2}{\gamma}}}\sum_{|\alpha|=m}\sum_{1\leq |\beta|\leq 2,\beta\leq\alpha}
C_\alpha^\beta\|\partial^\beta u\partial^{\alpha-\beta}\partial_yv\|_{2}\|\partial^\alpha v\|_{2}\\
&&+\sum_{m=5}^{\infty}\frac{\rho(t)^{2(m-3)}}{((m-3)!)^{\frac{2}{\gamma}}}\sum_{|\alpha|=m}\sum_{j=3}^{m-2}\sum_{|\beta|=j,\beta\leq\alpha}
C_\alpha^\beta\|\partial^\beta u\partial^{\alpha-\beta}\partial_yv\|_{2}\|\partial^\alpha v\|_{2}\\
&&+\sum_{m=3}^{\infty}\frac{\rho(t)^{2(m-3)}}{((m-3)!)^{\frac{2}{\gamma}}}\sum_{|\alpha|=m}\sum_{m-1\leq |\beta|\leq m,\beta\leq\alpha}
C_\alpha^\beta\|\partial^\beta u\partial^{\alpha-\beta}\partial_yv\|_{2}\|\partial^\alpha v\|_{2}
=\sum_{i=1}^3\tilde{I}_{i}
\eeno
We handle them term by term.\smallskip

{\bf Step 1. Estimate of $\tilde{I}_1$.}
We denote $\tilde{I}_1=\tilde{I}_{11}+\tilde{I}_{12}$ according to $|\beta|=1$ or $|\beta|=2$, where
\beno
\tilde{I}_{11}&=&\sum_{m=4}^{\infty}\frac{\rho(t)^{2(m-3)}}{((m-3)!)^{\frac{2}{\gamma}}}\sum_{|\alpha|=m}\sum_{|\beta|=1,\beta\leq\alpha}
C_\alpha^\beta\|\partial^\beta u\partial^{\alpha-\beta}\partial_yv\|_{2}\|\partial^\alpha v\|_{2}\\
&&+\sum_{|\alpha|=3}\sum_{|\beta|=1,\beta\leq\alpha}C_\alpha^\beta\|\partial^\beta u\partial^{\alpha-\beta}\partial_yv\|_{2}\|\partial^\alpha v\|_{2}
\eeno
Using $\partial^\beta u|_{y=0}=0$ and Sobolev embedding, the second term of $\tilde{I}_{11}$ can be controlled by
\beno
&&C_0\sum_{|\alpha|=3}\sum_{|\beta|=1,\beta\leq\alpha}\Big\|\frac{\partial^\beta u}{\psi}\psi\partial^{\alpha-\beta}\partial_yv\Big\|_{2}\|\partial^\alpha v\|_{2}\\
 &\leq&C_0\Big\|\frac{\partial^\beta u}{\psi}\Big\|_\infty\|v\|^2_{X^3}
\leq C_0(\delta)\|(u,\partial_yu)\|^{\frac12}_{\overline{H}^2}\|(\partial_yu,\partial_{yy}u)\|^{\frac12}_{\overline{H}^2}\|v\|^2_{X^3}.
\eeno
As in (\ref{eq:beta=1}), the first term of $\tilde{I}_{11}$ can be controlled by
\beno
&&\sum_{m=4}^{\infty}\frac{\rho(t)^{2(m-3)}}{((m-3)!)^{\frac{2}{\gamma}}}\sum_{|\alpha|=m}\sum_{|\beta|=1,\beta\leq\alpha}C_\alpha^\beta
\Big\|\frac{\partial^\beta u}{\psi}\Big\|_\infty\|\psi\partial^{\alpha-\beta}\partial_yv\|_{2}\|\partial^\alpha v\|_{2}\\
&&\leq C_0(\delta)\|(u,\partial_yu)\|^{\frac12}_{\overline{H}^2}\|(\partial_yu,\partial_{yy}u)\|^{\frac12}_{\overline{H}^2}\sum_{m=4}^{\infty}\frac{\rho(t)^{2(m-3)}m}{((m-3)!)^{\frac{2}{\gamma}}}\sum_{|\alpha|=m}
\|\partial^\alpha v\|^2_{2}\\
&&\leq C_0(\delta)\|(u,\partial_yu)\|^{\frac12}_{\overline{H}^2}\|(\partial_yu,\partial_{yy}u)\|^{\frac12}_{\overline{H}^2}\|v\|^2_{Y^3}.
\eeno
This shows that
$$\tilde{I}_{11}\leq C_0(\delta)\|(u,\partial_yu)\|^{\frac12}_{\overline{H}^2}\|(\partial_yu,\partial_{yy}u)\|^{\frac12}_{\overline{H}^2}\|v\|^2_{X^3\cap Y^3}. $$

For the case of $|\beta|=2$, similar arguments imply
\beno
\tilde{I}_{12}&=&\sum_{m=5}^{\infty}\frac{\rho(t)^{2(m-3)}}{((m-3)!)^{\frac{2}{\gamma}}}\sum_{|\alpha|=m}\sum_{|\beta|=2,\beta\leq\alpha}
C_\alpha^\beta\|\partial^\beta u\partial^{\alpha-\beta}\partial_yv\|_{2}\|\partial^\alpha v\|_{2}\\
&&+\sum_{m=3}^{4}\frac{\rho(t)^{2(m-3)}}{((m-3)!)^{\frac{2}{\gamma}}}\sum_{|\alpha|=m}\sum_{|\beta|=2,\beta\leq\alpha}C_\alpha^\beta\|\partial^\beta u\partial^{\alpha-\beta}\partial_yv\|_{2}\|\partial^\alpha v\|_{2},
\eeno
and the second term of the right hand is obviuously bounded by
$$C_0(\delta)\|(u,\partial_yu)\|^{\frac12}_{\overline{H}^3}\|(\partial_yu,\partial_{yy}u)\|^{\frac12}_{\overline{H}^3}
\|v\|^2_{X^3}.$$
The first term is bounded by
\beno
&&C_0\sum_{m=5}^{\infty}\frac{\rho(t)^{2(m-3)}}{((m-3)!)^{\frac{2}{\gamma}}}\sum_{|\alpha|=m}\sum_{|\beta|=2,\beta\leq\alpha}C_\alpha^\beta
\Big\|\frac{\partial^\beta u}{\psi}\Big\|_\infty\|\partial^{\alpha-\beta}Zv\|_{2}\|\partial^\alpha v\|_{2}\\
&\leq&C_0\Big\|\frac{\partial^\beta u}{\psi}\Big\|_\infty\sum_{m=5}^{\infty}\frac{\rho(t)^{2(m-3)}m(m-1)}{((m-3)!)^{\frac{2}{\gamma}}}
\Big(\sum_{|\alpha|=m}\|\partial^\alpha v\|^2_{2}\Big)^{\frac12}\Big(\sum_{|\alpha|=m-1}\|\partial^\alpha v\|^2_{2}\Big)^{\frac12}\\
&\leq& C_0(\delta)\|(u,\partial_yu)\|^{\frac12}_{\overline{H}^3}\|(\partial_yu,\partial_{yy}u)\|^{\frac12}_{\overline{H}^3}\|v\|^2_{Y^3}.
\eeno
This gives
\beno
\tilde{I}_{12}
\leq C_0(\delta)\|(u,\partial_yu)\|^{\frac12}_{\overline{H}^3}\|(\partial_yu,\partial_{yy}u)\|^{\frac12}_{\overline{H}^3}
\|v\|^2_{X^3\cap Y^3}.
\eeno

Collecting the results of $\tilde{I}_{11}$ and $\tilde{I}_{12}$ together, we obtain
\ben
\tilde{I}_1
\leq C_0(\delta)\|(u,\partial_yu)\|^{\frac12}_{\overline{H}^3}\|(\partial_yu,\partial_{yy}u)\|^{\frac12}_{\overline{H}^3}\|v\|^2_{X^3\cap Y^3}.
\een

{\bf Step 2. Estimate of $\tilde{I}_2$.} As in Lemma \ref{lem:x product estimate}, we decompose $\tilde{I}_2$ into two terms $\tilde{I}_{21}$ and $\tilde{I}_{22}$, where
\beno
\tilde{I}_{21}&=&\sum_{m=6}^{\infty}\frac{\rho(t)^{2(m-3)}}{((m-3)!)^{\frac{2}{\gamma}}}\sum_{|\alpha|=m}\sum_{j=3}^{[\frac{m}{2}]}\sum_{|\beta|=j,\beta\leq\alpha}
C_\alpha^\beta\|\partial^\beta u\partial^{\alpha-\beta}\partial_yv\|_{2}\|\partial^\alpha v\|_{2}
\eeno
and
\beno
\tilde{I}_{22}=\sum_{m=5}^{\infty}\frac{\rho(t)^{2(m-3)}}{((m-3)!)^{\frac{2}{\gamma}}}\sum_{|\alpha|=m}\sum_{j=[\frac{m}{2}]+1}^{m-2}\sum_{|\beta|=j,\beta\leq\alpha}
C_\alpha^\beta\|\partial^\beta u\partial^{\alpha-\beta}\partial_yv\|_{2}\|\partial^\alpha v\|_{2}.
\eeno
By Sobolev embedding, we have
\beno
\tilde{I}_{21}
&\leq&\sum_{m=6}^{\infty}\frac{\rho(t)^{2(m-3)}}{((m-3)!)^{\frac{2}{\gamma}}}\sum_{|\alpha|=m}\sum_{j=3}^{[\frac{m}{2}]}
\sum_{|\beta|=j,\beta\leq\alpha}C_\alpha^\beta\|\partial^\beta u\|^{\frac12}_2\|\partial_y \partial^\beta u\|^{\frac12}_2\|\partial^{\alpha-\beta}\partial_yv\|_{2}\|\partial^\alpha v\|_{2}\\
&&+\sum_{m=6}^{\infty}\frac{\rho(t)^{2(m-3)}}{((m-3)!)^{\frac{2}{\gamma}}}\sum_{|\alpha|=m}\sum_{j=3}^{[\frac{m}{2}]}
\sum_{|\beta|=j,\beta\leq\alpha}C_\alpha^\beta\|\partial_x\partial^\beta u\|^{\frac12}_2\|\partial_y \partial_x\partial^\beta u\|^{\frac12}_2\|\partial^{\alpha-\beta}\partial_yv\|_{2}\|\partial^\alpha v\|_{2}.
\eeno

Similar to (\ref{e:y derivative estimate low order}), the first term of the right hand side can be bounded by
\beno
&&\sum_{m=6}^{\infty}\frac{\rho(t)^{(m-3)}}{((m-3)!)^{\frac{1}{\gamma}}}\sum_{j=3}^{[\frac{m}{2}]}C_m^j
\Big(\sum_{|\beta|=j}\|\partial^\beta u\|_2\|\partial_{y} \partial^\beta u\|_2\Big)^{\frac12}\Big(\sum_{|\alpha|=m-j}\|\partial^{\alpha}\partial_yv\|^2_{2}\Big)^{\frac12}|v|_{m,3}\\
&&\leq C_0\sum_{m=6}^{\infty}\sum_{j=3}^{[\frac{m}{2}]}
\frac{C_m^j((j-3)!)^{\frac{1}{\gamma}}((m-j-3)!)^{\frac{1}{\gamma}}}
{((m-3)!)^{\frac{1}{\gamma}}}\\
&&\quad\times\frac{\rho(t)^{\frac{(j-3)}{2}}}{((j-3)!)^{\frac{1}{2\gamma}}}\Big(\sum\limits_{|\beta|=j}\|\partial_{y} \partial^\beta u\|^2_2\Big)^{\frac14}
|u|_{j,3}^{\frac12}|\partial_yv|_{m-j,3}
|v|_{m,3}\\
&&\leq C_0\| u\|^{\frac12}_{X^3}\|\partial_yu\|_{X^3}^{\frac12}\|\partial_yv\|_{X^3}\|v\|_{X^3},
\eeno
where we used the fact that
$$\frac{C_m^j((j-3)!)^{\frac{1}{\gamma}}((m-j-3)!)^{\frac{1}{\gamma}}}
{((m-3)!)^{\frac{1}{\gamma}}}\leq C_0\Big(C_{m-3}^{j-1}\Big)^{1-\frac{1}{\gamma}}j^{-1-\frac{2}{\gamma}}\leq C_0j^{-1-\frac{2}{\gamma}}$$
for $j\in \{3,\cdots,[\frac{m}{2}]\}$
and
\ben\label{eq:partial y commute}
\sum_{m=3}^{\infty}\frac{\rho(t)^{2(m-3)}}{((m-3)!)^{\frac{2}{\gamma}}}\sum_{|\alpha|=m}\|\partial_{y}\partial^\alpha u\|^2_{L^2}\leq
 C_0\| \partial_yu\|^2_{X^3}.
 \een
The second term is similar, and thus we get
\beno
\tilde{I}_{21}\leq C_0\| u\|^{\frac12}_{X^3}\|\partial_yu\|_{X^3}^{\frac12}
\|\partial_yv\|_{X^3}\|v\|_{X^3}.
\eeno
By Sobolev embedding and Hardy inequality, we have
\beno
&&\tilde{I}_{22}
\leq C_0(\delta)\sum_{m=5}^{\infty}\frac{\rho(t)^{2(m-3)}}{((m-3)!)^{\frac{2}{\gamma}}}\sum_{|\alpha|=m}\sum_{j=[\frac{m}{2}]+1}^{m-2}
\sum_{|\beta|=j,\beta\leq\alpha}C_\alpha^\beta(\|\partial^\beta u\|_2+\|\partial_y\partial^\beta u\|_2)\|\partial^\alpha v\|_{2}\\
&&\qquad\times\Big(\|\psi\partial^{\alpha-\beta}\partial_yv\|^{\frac12}_{2}\|\partial_y(\psi\partial^{\alpha-\beta}\partial_yv)\|^{\frac12}_{2}
+\|\partial_{x}(\psi\partial^{\alpha-\beta}\partial_yv)\|^{\frac12}_{2}
\|\partial_{x}\partial_y(\psi\partial^{\alpha-\beta}\partial_yv)\|^{\frac12}_{2}\Big).
\eeno
We estimate the term $\|\psi\partial^{\alpha-\beta}\partial_yv\|^{\frac12}_{2}\|\partial_y(\psi\partial^{\alpha-\beta}\partial_yv)\|^{\frac12}_{2}$. By similar arguments as in (\ref{eq:middle term-Young inequ}), it can be bounded by
\beno
&& \sum_{m=5}^{\infty}\frac{\rho(t)^{(m-3)}}{((m-3)!)^{\frac{1}{\gamma}}}\sum_{j=[\frac{m}{2}]+1}^{m-2}
C_m^j\Big(\sum_{|\alpha|=m-j+1}\|\partial^\alpha v\|_2\|\partial_y \partial^\alpha v\|_2\Big)^{\frac12}\Big(\sum_{|\alpha|=j}
\|\partial^\alpha u\|^2_{2}+\|\partial_y\partial^\alpha u\|^2_{2}\Big)^{\frac12}|v|_{m,3}\\
&&\leq C_0(\delta)\sum_{m=5}^{\infty}\sum_{j=[\frac{m}{2}]+1}^{m-2}
\frac{C_m^j((m-j-2)!)^{\frac{1}{\gamma}}((j-3)!)^{\frac{1}{\gamma}}}{((m-3)!)^{\frac{1}{\gamma}}}|v|_{m-j+1,3}^{\frac12}|v|_{m,3}\\
&&\quad\times\frac{\rho(t)^{\frac{m-j-2}{2}}}{((m-j-2)!)^{\frac{1}{2\gamma}}}\Big(\sum_{|\alpha|=m-j+1}\|\partial_y\partial^\alpha v\|^2_2\Big)^{\frac14}\frac{\rho(t)^{j-3}}{((j-3)!)^{\frac{1}{\gamma}}}\Big(\sum_{|\alpha|=j}\|\partial^\alpha u\|^2_{2}+\|\partial_y\partial^\alpha u\|^2_2\Big)^{\frac12}\\
&&\leq C_0(\delta)\sum_{m=6}^{\infty}\sum_{j=[\frac{m}{2}]+1}^{m-2}\frac{1}{(m-j)^2}
|v|_{m-j+1,3}^{\frac12}|v|_{m,3}
\frac{\rho(t)^{\frac{m-j-2}{2}}}{((m-j-2)!)^{\frac{1}{2\gamma}}}\Big(\sum_{|\alpha|=m-j+1}\|\partial_y\partial^\alpha v\|^2_2\Big)^{\frac14}\\
&&\quad\times \frac{\rho(t)^{j-3}}{((j-3)!)^{\frac{1}{\gamma}}}\Big(\sum_{|\alpha|=j}\|\partial^\alpha u\|^2_{2}+\|\partial_y\partial^\alpha u\|^2_2\Big)^{\frac12}\\
&&\leq C_0(\delta)\| v\|^{\frac32}_{X^3}\| \partial_yv\|_{X^3}^{\frac12}\|(u,\partial_yu)\|_{X^3},
\eeno
where we used (\ref{eq:partial y commute}) and
$$\frac{C_m^j((m-j-2)!)^{\frac{1}{\gamma}}((j-3)!)^{\frac{1}{\gamma}}}{((m-3)!)^{\frac{1}{\gamma}}}\leq C_0\Big(C_{m-3}^{j-1}\Big)^{1-\frac{1}{\gamma}}\frac{1}{(m-j)^2}$$
for $j\in \{[\frac{m}{2}]+1,\cdots,m-2\}$.

The other terms  can also be bounded by
\beno
C_0(\delta)\| v\|^{\frac32}_{X^3}\| \partial_yv\|_{X^3}^{\frac12}\|(u,\partial_yu)\|_{X^3}.
\eeno
This shows that
\beno
\tilde{I}_{22}\leq  C_0(\delta)\| v\|^{\frac32}_{X^3}\| \partial_yv\|_{X^3}^{\frac12}\|(u,\partial_yu)\|_{X^3}.
\eeno

Collecting $\tilde{I}_{21}$ and $\tilde{I}_{22}$ together, we arrive at
\beno
\tilde{I}_2\leq C_0\| u\|^{\frac12}_{X^3}\|\partial_yu\|_{X^3}^{\frac12}
\|\partial_yv\|_{X^3}\|v\|_{X^3}+C_0(\delta)\| v\|^{\frac32}_{X^3}\| \partial_yv\|_{X^3}^{\frac12}\|(u,\partial_yu)\|_{X^3}.
\eeno

{\bf Step 3. Estimate of $\tilde{I}_3$.} We first decompose $\tilde{I}_3=\tilde{I}_{31}+\tilde{I}_{32}$ according to  the value of $|\beta|$ as in Lemma \ref{lem:x product estimate}. By similar computations as in (\ref{eq:beta=1}) and Hardy inequality, we get
\beno
\tilde{I}_{31}&\leq& C_0(\delta)\sum_{m=4}^{\infty}\frac{\rho(t)^{2(m-3)}}{((m-3)!)^{\frac{2}{\gamma}}}\sum_{|\alpha|=m}\sum_{|\beta|=m-1,\beta\leq\alpha}
C_\alpha^\beta(\|\partial^\beta u\|_2+\|\partial_y\partial^\beta u\|_2)\|\psi\partial^{\alpha-\beta}\partial_yv\|_\infty\|\partial^\alpha v\|_{2}\\
&&+C_0\sum_{|\alpha|=3,|\beta|=2,\beta\leq \alpha}\|\partial^\beta u \partial^{\alpha-\beta}\partial_yv\|_2\|\partial^\alpha v\|_2\\
&\leq& C_0(\delta)\|\partial^2 v\|_\infty\sum_{m=4}^{\infty}\frac{\rho(t)^{2(m-3)}m}{((m-3)!)^{\frac{2}{\gamma}}}
\Big(\sum_{|\alpha|=m}
\|\partial^\alpha v\|_{2}\Big)^{\frac12}\Big(\sum_{|\beta|=m-1}(\|\partial^\beta u\|_2^2+\|\partial_y\partial^\beta u\|^2_2)\Big)^{\frac12}\\
&&+C_0(\delta)\|\partial^2 v\|_\infty\|v\|_{X^3}\sum_{|\beta|=2}(\|\partial^\beta u\|_2+\|\partial_y\partial^\beta u \|_2)\\
&\leq& C_0(\delta)\|v\|^{\frac12}_{\overline{H}^3}\|\partial_yv\|^{\frac12}_{\overline{H}^3}\|(u,\partial_yu)\|_{X^3}\|v\|_{X^3}.
\eeno
Similarly, for the term $\tilde{I}_{32}$, there holds
\beno
\tilde{I}_{32}&\leq&C_0(\delta)\sum_{m=4}^{\infty}\frac{\rho(t)^{2(m-3)}}{((m-3)!)^{\frac{2}{\gamma}}}\sum_{|\alpha|=m}\sum_{|\beta|=m,\beta\leq\alpha}
C_\alpha^\beta(\|\partial^\beta u\|_2+\|\partial_y\partial^\beta u\|_2)\|\psi\partial^{\alpha-\beta}\partial_yv\|_\infty\|\partial^\alpha v\|_{2}\\
&&+C_0\sum_{|\alpha|=3}\|\partial^\alpha u \partial_yv\|_2\|\partial^\alpha v\|_2\\
&\leq& C_0(\delta)\|Z v\|_\infty\sum_{m=4}^{\infty}\frac{\rho(t)^{2(m-3)}}{((m-3)!)^{\frac{2}{\gamma}}}\Big(\sum_{|\alpha|=m}(\|\partial^\alpha u\|_2^2+\|\partial_y\partial^\alpha u\|^2_2)\Big)^{\frac12}\Big(\sum_{|\alpha|=m}\|\partial^\alpha v\|^2_{2}\Big)^{\frac12}\\
&&+C_0\sum_{|\alpha|=3}\Big\|\frac{\partial^\alpha u}{\psi}\Big \|_2\|Zv\|_\infty\|\partial^\alpha v\|_2\\
&\leq&C_0(\delta)\|v\|^{\frac12}_{\overline{H}^3}\|\partial_yv\|^{\frac12}_{\overline{H}^3}\|(u,\partial_yu)\|_{X^3}\|v\|_{X^3}.
\eeno
Finally, collecting the $\tilde{I}_{31}$ and $\tilde{I}_{32}$ together, we obtain
\beno
\tilde{I}_3\leq C_0(\delta)\|v\|^{\frac12}_{\overline{H}^3}\|\partial_yv\|^{\frac12}_{\overline{H}^3}\|(u,\partial_yu)\|_{X^3}\|v\|_{X^3}.
\eeno

Collecting the estimates in Step 1-Step 3, we complete the proof of the first inequality.\smallskip

(b) We estimate the second inequality which is used to deal with the linear term like $v_a\partial_yu$. As in (a), we first have
\beno
&&\sum_{m=3}^{\infty}\frac{\rho(t)^{2(m-3)}}{((m-3)!)^{\frac{2}{\gamma}}}\sum_{|\alpha|=m}\big|\langle \partial^\alpha(u\partial_yv),\partial^\alpha v\rangle\big|\\
&&\leq \sum_{m=3}^{\infty}\frac{\rho(t)^{2(m-3)}}{((m-3)!)^{\frac{2}{\gamma}}}\sum_{|\alpha|=m}\sum_{0<\beta\leq\alpha}C_\alpha^\beta\big|\langle \partial^\beta u\partial^{\alpha-\beta}\partial_yv,\partial^\alpha v\rangle\big|\\
&&\quad+C_0(\delta)\|(u,\partial_yu)\|_{\overline{H}^1}^{\frac12}\|(\partial_yu,\partial_{yy}u)\|_{\overline{H}^1}^{\frac12} \big(\|v\|^2_{Y^3}+\|v\|^2_{X^3}\big),
\eeno
where we used
\beno
|\partial^{\alpha}\partial_yv-\partial_y\partial^{\alpha}v|\leq C_0|\alpha||\partial^{\beta}\partial_yv|,\quad |\beta|=m-1,\,\,\beta\leq \alpha.
\eeno
Then, recall the definition of $\tilde{I}$ and
$
\tilde{I}\leq\sum\limits_{i=1}^3\tilde{I}_{i}.
$
We will handle them term by term.\smallskip

{\bf Step 1: Estimate of $\tilde{I}_1$.}
As in (a), we know that
\beno
\tilde{I}_1
\leq C_0(\delta)\|(u,\partial_yu)\|_{\overline{H}^3}^{\frac12}\|(\partial_yu,\partial_{yy}u)\|_{\overline{H}^3}^{\frac12}\|v\|^2_{X^3\cap Y^3}.
\eeno

{\bf Step 2. Estimate of $\tilde{I}_2$.} By Sobolev embedding, we have
\beno
\tilde{I}_2
&\leq&\sum_{m=6}^{\infty}\frac{\rho(t)^{2(m-3)}}{((m-3)!)^{\frac{2}{\gamma}}}\sum_{|\alpha|=m}\sum_{j=3}^{m-3}
\sum_{|\beta|=j,\beta\leq\alpha}C_\alpha^\beta\Big\|\frac{\partial^\beta u}{\psi}\Big\|_\infty\|\psi\partial^{\alpha-\beta}\partial_yv\|_{2}\|\partial^\alpha v\|_{2}\\
&&+\sum_{m=6}^{\infty}\frac{\rho(t)^{2(m-3)}}{((m-3)!)^{\frac{2}{\gamma}}}\sum_{|\alpha|=m}
\sum_{|\beta|=m-2,\beta\leq\alpha}C_\alpha^\beta\|\partial^\beta u \partial^{\alpha-\beta}\partial_yv\|_{2}\|\partial^\alpha v\|_{2}\\
&&+C_0(\delta)\|(u,\partial_yu)\|_{\overline{H}^4}^{\frac12}\|(\partial_yu,\partial_{yy}u)\|_{\overline{H}^4}^{\frac12}\|v\|^2_{X^3}\\
&=&\sum_{i=1}^2\tilde{I}_{2i}+C_0(\delta)\|(u,\partial_yu)\|_{\overline{H}^4}^{\frac12}\|(\partial_yu,\partial_{yy}u)\|_{\overline{H}^4}^{\frac12}\|v\|^2_{X^3}.
\eeno

Applying the same arguments in (\ref{eq:beta=1}) and discrete Young convolution inequality, we have
\beno
\tilde{I}_{21}&\leq&\sum_{m=6}^{\infty}\frac{\rho(t)^{(m-3)}}{((m-3)!)^{\frac{1}{\gamma}}}\sum_{j=3}^{m-3}C_m^j
\Big(\sum_{|\beta|=j}\Big\|\frac{\partial^\beta u}{\psi}\Big\|^2_\infty\Big)^{\frac12}\Big(\sum_{|\alpha|=m+1-j}\|\partial^{\alpha}v\|^2_{2}\Big)^{\frac12}|v|_m\\
&\leq&C_0\sum_{m=6}^{\infty}\sum_{j=3}^{m-3}
\frac{C_m^j((j-3)!)^{\frac{1}{\gamma}}((m-j-2)!)^{\frac{1}{\gamma}}}
{((m-3)!)^{\frac{1}{\gamma}}\sqrt{m-j-2}\sqrt{m-3}}\\
&&\times\frac{\rho(t)^{(j-3)}}{((j-3)!)^{\frac{1}{\gamma}}}\Big(\sum_{|\beta|=j}\Big\|\frac{\partial^\beta u}{\psi}\Big\|^2_\infty\Big)^{\frac12}
\sqrt{m-j-2}|v|_{m-j+1}\sqrt{m-3}|v|_m\\
&\leq&C_0(\delta)\|(u,\partial_y u)\|_{X^4}^{\frac12}\|(\partial_{y} u,\partial_{yy} u)\|_{X^4}^{\frac12}
\|v\|^2_{Y^3},
\eeno
where we used the natation (\ref{eq:a m j})
\beno
&&\frac{C_m^j((m-j-2)!)^{\frac{1}{\gamma}}((j-3)!)^{\frac{1}{\gamma}}}{((m-3)!)^{\frac{1}{\gamma}}\sqrt{m-3}\sqrt{m-j-2}}\leq C_0 a_{m,j},
\eeno
and
\beno
\sum_{m=3}^{\infty}\frac{\rho(t)^{2(m-3)}}{((m-3)!)^{\frac{2}{\gamma}}}\sum_{|\alpha|=m}\Big\|\frac{\partial^\beta u}{\psi}\Big\|^2_\infty
&\leq&
 C_0(\delta)\|(u,\partial_y u)\|_{X^4}\|(\partial_{y} u,\partial_{yy} u)\|_{X^4}.
\eeno

Summing up these estimate, we obtain
\beno
\tilde{I}_{2}\leq C_0(\delta)\|(u,\partial_y u)\|_{X^4}^{\frac12}\|(\partial_{y} u,\partial_{yy} u)\|_{X^4}^{\frac12}
\| v\|^2_{Y^3\cap X^3}+\tilde{I}_{22}.
\eeno

{\bf Step 3. Estimate of $\tilde{I}_3$ and $\tilde{I}_{22}$.} We first write $\tilde{I}_3=\tilde{I}_{31}+\tilde{I}_{32}$ according to the value of $|\beta|$ as above. The argument for  $\tilde{I}_{31}$ is similar to $I_{31}$:
\beno
&&\sum_{m=4}^{\infty}\frac{\rho(t)^{2(m-3)}}{((m-3)!)^{\frac{2}{\gamma}}}\sum_{|\alpha|=m}\sum_{|\beta|=m-1,\beta\leq\alpha}
C_\alpha^\beta\Big\|\frac{\partial^\beta u}{\psi}\Big\|_{L^\infty_yL^2_x}\|Z\partial^{\alpha-\beta}v\|_{L^2_yL^\infty_x}\|\partial^\alpha v\|_{2}\\
&&\leq C_0(\delta)\| v\|_{\overline{H}^3}\sum_{m=4}^{\infty}\frac{\rho(t)^{2(m-3)}m}{((m-3)!)^{\frac{2}{\gamma}}}
\Big(\sum_{|\alpha|=m}
\|\partial^\alpha v\|^2_{2}\Big)^{\frac12}\Big(\sum_{|\beta|=m-1}\Big\|\frac{\partial^\beta u}{\psi} \Big\|^2_{L^\infty_yL^2_x}\Big)^{\frac12}\\
&&\leq C_0(\delta)\| v\|_{\overline{H}^3}\| v\|_{X^3}\Big(\sum_{m=4}^{\infty}\frac{\rho(t)^{2(m-4)}}{((m-4)!)^{\frac{2}{\gamma}}}
\sum_{|\beta|=m-1}\Big\|\frac{\partial^\beta u}{\psi} \Big\|^2_{L^\infty_yL^2_x}\Big)^{\frac12}\\
&&\leq C_0(\delta)\|(u,\partial_y u)\|_{X^3}^{\frac12}\|(\partial_{y} u,\partial_{yy} u)\|_{X^3}^{\frac12}
\| v\|_{\overline{H}^3}\| v\|_{X^3},
\eeno
which gives
\beno
\tilde{I}_{31}\leq C_0(\delta)\|(u,\partial_y u)\|_{X^3}^{\frac12}\|(\partial_{y} u,\partial_{yy} u)\|_{X^3}^{\frac12}
\| v\|_{\overline{H}^3}\| v\|_{X^3}.
\eeno
Then a proof similar to that used to treat $\tilde{I}_{32}$ and $\tilde{I}_{22}$ yields that
\beno
\tilde{I}_{32}&\leq&C_0(\delta)\|(u,\partial_y u)\|_{X^3}^{\frac12}\|(\partial_{y} u,\partial_{yy} u)\|_{X^3}^{\frac12}
\| v\|_{\overline{H}^2}\| v\|_{X^3},\\
\tilde{I}_{22}&\leq&C_0(\delta)\|(u,\partial_y u)\|_{X^3}^{\frac12}\|(\partial_{y} u,\partial_{yy} u)\|_{X^3}^{\frac12}
\| v\|_{\overline{H}^4}\| v\|_{X^3}.
\eeno

Collecting $\tilde{I}_{31}, \tilde{I}_{32}, \tilde{I}_{22}$ together, we obtain
\beno
\tilde{I}_{31}+\tilde{I}_{32}+\tilde{I}_{22}
\leq C_0(\delta)\|(u,\partial_y u)\|_{X^3}^{\frac12}\|(\partial_{y} u,\partial_{yy} u)\|_{X^3}^{\frac12}
\| v\|_{\overline{H}^4}\| v\|_{ X^3}.
\eeno

Collecting the estimates in Step 1-Step 3, the proof is completed.
\endproof

\smallskip

The goal of the following lemma is to deal with the terms in 3) like $u\partial_xu^a$.
\begin{Lemma}\label{e:x sublinear product estimate}
For $k=2,3$, there holds
\beno
I_k'&=&\sum_{m=k}^{\infty}\frac{\rho(t)^{2(m-k)}}{((m-k)!)^{\frac{2}{\gamma}}}\sum_{|\alpha|=m}\big|\langle \partial^\alpha(uv),\partial^\alpha w\rangle_{L^2}\big|\\
&\leq&  C_0\|v\|_{X^{k+1}}^{\frac12}\|\partial_yv\|_{{X^{k+1}}}^{\frac12}
\big(\|u\|^2_{X^k}+\|w\|^2_{X^k\cap Y^k}\big).
\eeno
\end{Lemma}

\no{\bf Proof}.
Similarly as above, the argument for $k=2$ and $k=3$ is the same, and we only give a proof for $k=3.$ We divide $I'$ into three terms:
\beno
\sum_{i=1}^3I_{i}'&=&\sum_{m=3}^{\infty}\frac{\rho(t)^{2(m-3)}}{((m-3)!)^{\frac{2}{\gamma}}}\sum_{|\alpha|=m}\sum_{|\beta|\leq 3,\beta\leq\alpha}C_\alpha^\beta\|\partial^\beta u\partial^{\alpha-\beta}v\|_{2}\|\partial^\alpha w\|_{2}\\
&&+\sum_{m=7}^{\infty}\frac{\rho(t)^{2(m-3)}}{((m-3)!)^{\frac{2}{\gamma}}}\sum_{|\alpha|=m}\sum_{j=4}^{m-3}
\sum_{|\beta|=j,\beta\leq\alpha}C_\alpha^\beta\|\partial^\beta u\partial^{\alpha-\beta}v\|_{2}\|\partial^\alpha w\|_{2}\\
&&+\sum_{m=3}^{\infty}\frac{\rho(t)^{2(m-3)}}{((m-3)!)^{\frac{2}{\gamma}}}\sum_{|\alpha|=m}\sum_{m-2\leq|\beta|\leq m,\beta\leq\alpha}C_\alpha^\beta\|\partial^\beta u\partial^{\alpha-\beta}v\|_{2}\|\partial^\alpha w\|_{2}.
\eeno

{\bf Step 1. Estimate of $I_1'$.}
Again, we denote $I_1'=\sum\limits_{i=0}^3I_{1i}'$ according to the value of $|\beta|$. It suffices to estimate $I_{13}'$, since the others are similar.
By Sobolev embedding, $I_{13}'$ is clearly bounded by
\beno
\sum_{m=6}^{\infty}\frac{\rho(t)^{2(m-3)}}{((m-3)!)^{\frac{2}{\gamma}}}\sum_{|\alpha|=m}\sum_{|\beta|=3,\beta\leq\alpha}C_\alpha^\beta\|\partial^\beta u\partial^{\alpha-\beta}v\|_{2}\|\partial^\alpha w\|_{2}
+C_0\|v\|_{\overline{H}^4}^{\frac12}\|\partial_yv\|_{\overline{H}^4}^{\frac12}\|u\|_{\overline{H}^3}\|w\|_{X^3}.
\eeno
Applying the same technique as (\ref{eq:beta=1}), we get
\beno
&&\sum_{m=6}^{\infty}\frac{\rho(t)^{2(m-3)}}{((m-3)!)^{\frac{2}{\gamma}}}\sum_{|\alpha|=m}\sum_{|\beta|=3,\beta\leq\alpha}C_\alpha^\beta\|\partial^\beta u\partial^{\alpha-\beta}v\|_{2}\|\partial^\alpha w\|_{2}\\
&&\leq C_0\| u\|_{\overline{H}^4}\sum_{m=6}^{\infty}\frac{\rho(t)^{(m-3)}C_m^3}{((m-3)!)^{\frac{1}{\gamma}}}\Big(\sum_{|\alpha|=m-3}\|\partial^\alpha v\|^2_{L^\infty_yL^2_x}\Big)^{\frac12}|w|_{m,3}\\
&&\leq C_0\|v\|_{X^3}^{\frac12}\|\partial_yv\|_{X^3}^{\frac12}\| u\|_{\overline{H}^4}\|w\|_{X^3}.
\eeno
Thus, we have
$$I_{13}'\leq  C_0\|v\|_{X^3}^{\frac12}\|\partial_yv\|_{X^3}^{\frac12}\| u\|_{\overline{H}^4}\|w\|_{X^3}. $$
The same argument implies that
\beno
&&I_{11}'+I_{12}'\leq  C_0\|v\|_{X^3}^{\frac12}\|\partial_yv\|_{X^3}^{\frac12}\| u\|_{\overline{H}^3}\|w\|_{X^3},\\
&&I_{10}'\leq  C_0\|v\|_{X^3}^{\frac12}\|\partial_yv\|_{X^3}^{\frac12}\| u\|_{\overline{H}^1}\|w\|_{X^3}.
\eeno
Therefore, we obtain
\beno
I_1'\leq  C_0\|v\|_{X^3}^{\frac12}\|\partial_yv\|_{X^3}^{\frac12}\| u\|_{\overline{H}^4}\|w\|_{X^3}.
\eeno

{\bf Step 2. Estimate of $I_2'$.} By lemma \ref{lem:identity} and  \textcolor[rgb]{0.00,0.00,0.00}{ discrete Young inequality}, we have
\beno
I_{2}'&\leq& \sum_{m=7}^{\infty}\frac{\rho(t)^{(m-3)}}{((m-3)!)^{\frac{1}{\gamma}}}\sum_{j=4}^{m-3}
C_m^j\Big(\sum_{|\beta|=j}\|\partial^\beta u\|^2_2\Big)^{\frac12}\Big(\sum_{|\alpha|=m-j}
\|\partial^\alpha v\|^2_\infty\Big)^{\frac12}|w|_{m,3}\\
&\leq& C_0\sum_{m=7}^{\infty}\sum_{j=4}^{m-3}
\frac{C_m^j((m-j-2)!)^{\frac{1}{\gamma}}((j-3)!)^{\frac{1}{\gamma}}}{((m-3)!)^{\frac{1}{\gamma}}}\frac{\rho(t)^{(m-j-2)}}{((m-j-2)!)^{\frac{1}{\gamma}}}
\Big(\sum_{|\beta|=m-j}\|\partial^\beta v\|^2_\infty\Big)^{\frac12}|u|_{j,3}|w|_{m,3}\\
&\leq& C_0\sum_{m=7}^{\infty}\sum_{j=4}^{m-3}
a_{m,j}\frac{\rho(t)^{(m-j-3)}}{((m-j-2)!)^{\frac{1}{\gamma}}}
\sqrt{m-j-2}\Big(\sum_{|\beta|=m-j}\|\partial^\beta v\|^2_\infty\Big)^{\frac12}|u|_{j,3}\sqrt{m-3}|w|_{m,3}\\
&\leq& C_0\|v\|_{X^4}^{\frac12}\|\partial_yv\|_{X^4}^{\frac12}\| u\|_{X^3}\|w\|_{Y^3},
\eeno
where we used (\ref{eq:a m j})
\beno
&&\frac{C_m^j((m-j-2)!)^{\frac{1}{\gamma}}((j-3)!)^{\frac{1}{\gamma}}}{((m-3)!)^{\frac{1}{\gamma}}\sqrt{m-j-2}\sqrt{m-3}}\leq C_0a_{m,j}.
\eeno

{\bf Step 3. Estimate of $I_3'$.} Similarly as above, we write $I_3'$ into three terms:
\beno
I_{31}'+I_{32}'+I_{33}'&=&\sum_{m=3}^{\infty}\frac{\rho(t)^{2(m-3)}}{((m-3)!)^{\frac{2}{\gamma}}}\sum_{|\alpha|=m}\sum_{|\beta|=m-2 ,\beta\leq\alpha}C_\alpha^\beta\|\partial^\beta u\partial^{\alpha-\beta}v\|_{2}\|\partial^\alpha w\|_{2}\\
&&+\sum_{m=3}^{\infty}\frac{\rho(t)^{2(m-3)}}{((m-3)!)^{\frac{2}{\gamma}}}\sum_{|\alpha|=m}\sum_{|\beta|=m-1 ,\beta\leq\alpha}C_\alpha^\beta\|\partial^\beta u\partial^{\alpha-\beta}v\|_{2}\|\partial^\alpha w\|_{2}\\
&&+\sum_{m=3}^{\infty}\frac{\rho(t)^{2(m-3)}}{((m-3)!)^{\frac{2}{\gamma}}}\sum_{|\alpha|=m}\sum_{|\beta|=m,\beta\leq\alpha}C_\alpha^\beta\|\partial^\beta u\partial^{\alpha-\beta}v\|_{2}\|\partial^\alpha w\|_{2}
\eeno
Similar to the estimate of (a) in Lemma \ref{lem:x product estimate}, we know that
\beno
I_{32}'\leq C_0\|v\|_{\overline{H}^3}^{\frac12}\|\partial_yv\|_{\overline{H}^3}^{\frac12} \|u\|_{X^3}\|w\|_{X^3}.
\eeno
Similarly, there holds
\beno
I_{31}'+I_{33}'\leq C_0\|v\|_{\overline{H}^3}^{\frac12}\|\partial_yv\|_{\overline{H}^3}^{\frac12} \|u\|_{X^3}\|w\|_{X^3}.
\eeno
Finally, we obtain
\beno
I_3'\leq C_0\|v\|_{\overline{H}^3}^{\frac12}\|\partial_yv\|_{\overline{H}^3}^{\frac12} \|(u,w)\|^2_{X^3}.
\eeno
Collecting the estimates in Step 1-Step 3, the proof is completed.
\endproof
\smallskip

 Finally, we deal with the terms in 4). Recall the assumption $(H1)$ and notice that
  $\varepsilon^{1-\gamma}\partial_y(u^{p}(t,x,z))=\varepsilon^{-\gamma}\partial_zu^{p}(t,x,z)$ is of $\varepsilon^{-\gamma}$ order. Then we have the following estimate.

\begin{Lemma}\label{e:singular term estimate}
For $\gamma\leq 1$, there hold
\beno
J=\sum_{m=3}^{\infty}\frac{\rho(t)^{2(m-3)}\varepsilon^{-\gamma}}{((m-3)!)^{\frac{2}{\gamma}}}\sum_{|\alpha|=m}\big|\langle \partial^\alpha(\tilde{v}\partial_zu^{p}),\partial^\alpha u\rangle\big|
\leq C_0(\|(u,v)\|^2_{X^3\cap Y^3}+\varepsilon^4),
\eeno
and
\beno
J'=\sum_{m=2}^{\infty}\frac{\rho(t)^{2(m-2)}\varepsilon^{-\gamma}}{((m-2)!)^{\frac{2}{\gamma}}}\sum_{|\alpha|=m}\big|\langle \partial^\alpha(\tilde{v}\partial_zu^{p}),\partial^\alpha \eta\rangle\big|\leq C_0(\|(\eta,U)\|^2_{X^2\cap Y^2}+\varepsilon^4).
\eeno
\end{Lemma}
\no{\bf Proof.}
 $J$ can be controlled by the sum of $J_1$ and $J_2$, where
\beno
J_1&=&\sum_{m=3}^{\infty}\frac{\rho(t)^{2(m-3)}\varepsilon^{-\gamma}}{((m-3)!)^{\frac{2}{\gamma}}}  \sum_{|\alpha|=m}\sum_{|\beta|\leq 3,\beta\leq\alpha}C_\alpha^\beta\big|\langle \partial^\beta \tilde{v}\partial^{\alpha-\beta}\partial_zu^{p},\partial^\alpha u\rangle\big|,\\
J_2&=&\sum_{m=4}^{\infty}\frac{\rho(t)^{2(m-3)}\varepsilon^{-\gamma}}{((m-3)!)^{\frac{2}{\gamma}}} \textcolor[rgb]{0.00,0.00,0.00}{ \sum_{|\alpha|=m}\sum_{\beta\leq \alpha, |\beta|\geq 4}}C_\alpha^\beta\big|\langle \partial^\beta \tilde{v}\partial^{\alpha-\beta}\partial_zu^{p},\partial^\alpha u\rangle\big|.
\eeno

For $J_1$, we have
\beno
J_1
&=&\sum_{m=3}^{6}\frac{\rho(t)^{2(m-3)}}{((m-3)!)^{\frac{2}{\gamma}}} \textcolor[rgb]{0.00,0.00,0.00}{ \sum_{|\alpha|=m}\sum_{|\beta|\leq 3, \beta\leq \alpha}}C_\alpha^\beta\Big|\Big<
\frac{\partial^\beta \tilde{v}}{y^\gamma}z^\gamma\partial^{\alpha-\beta}\partial_zu^{p},\partial^\alpha u\Big>\Big|\\
&&+\sum_{m=7}^{\infty}\frac{\rho(t)^{2(m-3)}}{((m-3)!)^{\frac{2}{\gamma}}}\textcolor[rgb]{0.00,0.00,0.00}{ \sum_{|\alpha|=m}\sum_{|\beta|\leq 3, \beta\leq \alpha}}C_\alpha^\beta\Big|\Big<
\frac{\partial^\beta \tilde{v}}{y^\gamma}z^\gamma\partial^{\alpha-\beta}\partial_zu^{p},\partial^\alpha u\Big>\Big|.
\eeno
Note that (\ref{e:uniform boundness for approximate solution}) implies that
\beno
\|z^\gamma\partial^{\alpha-\beta}\partial_zu^{p}\|_{\infty}\leq C_0,
\eeno
thus the first term of the right hand side can be controlled by
\beno
 C_0(\|(u,v)\|^2_{X^3}+\varepsilon^4).
\eeno
Applying the technique from Lemma \ref{lem:identity} and (\ref{e:uniform boundness for approximate solution}), the second term can be bounded by
\beno
&&C_0\big(\|(u,v)\|_{\overline{H}^4}+\varepsilon^2\big)\sum_{m=7}^{\infty}\sum_{j=0}^3\frac{C_m^j((m-j-3)!)^{\frac{1}{\gamma}}}{(m-3)!)^{\frac{1}{\gamma}}}
|v|_{m,3}\\
&&\quad\times\frac{\rho(t)^{(m-j-3)}}{((m-j-3)!)^{\frac{1}{\gamma}}}\Big(\sum_{|\alpha|=m-j}\|z^\gamma\partial^\alpha \partial_zu^{p}\|^2_{2}\Big)^{\frac12}\leq C_0\big(\|(u,v)\|_{\overline{H}^4}+\varepsilon^2\big)\|v\|_{X^3}.
\eeno
This gives
$$J_1\leq C_0(\|(u,v)\|^2_{X^3}+\varepsilon^4).$$

Now we turn to deal with the term $J_2$.
\beno
J_2
&=&\sum_{m=7}^{\infty}\frac{\rho(t)^{2(m-3)}\varepsilon^{-\gamma}}{((m-3)!)^{\frac{2}{\gamma}}}\sum_{|\alpha|=m}\sum_{j=4}^{m-3} \textcolor[rgb]{0.00,0.00,0.00}{\sum_{\beta\leq \alpha, |\beta|= j}}C_\alpha^\beta\big|\langle \partial^\beta \tilde{v}\partial^{\alpha-\beta}\partial_zu^{p},\partial^\alpha u\rangle\big|\\
&&+\sum_{m=4}^{\infty}\frac{\rho(t)^{2(m-3)}\varepsilon^{-\gamma}}{((m-3)!)^{\frac{2}{\gamma}}}\sum_{|\alpha|=m}\sum_{j=m-2}^m\sum_{\beta\leq \alpha, |\beta|= j}C_\alpha^\beta\big|\langle \partial^\beta \tilde{v}\partial^{\alpha-\beta}\partial_zu^{p},\partial^\alpha u\rangle\big|\\
&=&J_{21}+J_{22}.
\eeno

For $J_{21}$, it can be controlled by
\beno
&& \sum_{m=7}^{\infty}\sum_{j=4}^{m-3}\frac{C_m^j|u|_{m,3}\rho(t)^{(m-3)}\varepsilon^{-\gamma}}{((m-3)!)^{\frac{1}{\gamma}}}\cdot\Big(\sum_{|\alpha|=m}\sum_{\beta\leq \alpha, |\beta|= j}\|\partial^\beta \tilde{v}\partial^{\alpha-\beta}\partial_zu^{p}\|_{2}^2\Big)^{\frac12}\\
&&\leq\sum_{m=7}^{\infty}\sum_{j=4}^{m-3}\Bigg\{\frac{C_m^j((j-2)!)^{\frac{1}{\gamma}}((m-j-3)!)^{\frac{1}{\gamma}}
\|\chi_j(y)y^{1-\gamma}\|_\infty}{((m-3)!)^{\frac{1}{\gamma}}\sqrt{m-3}\sqrt{j-2}}
\frac{\rho(t)^{(j-2)}\sqrt{j-2}}{((j-2)!)^{\frac{1}{\gamma}}}\cdot\Big(\sum\limits_{|\alpha|=j}\|\partial_y\partial^\alpha \tilde{v}\|_2^2\Big)^{\frac12}\\
&&\quad+\frac{C_m^j\|(1-\chi_j(y))y^{-\gamma}\|_\infty((j-3)!)^{\frac{1}{\gamma}}((m-j-3)!)^{\frac{1}{\gamma}}}{((m-3)!)^{\frac{1}{\gamma}}\sqrt{m-3}\sqrt{j-3}}
\frac{\rho(t)^{(j-3)}\sqrt{j-3}}{((j-3)!)^{\frac{1}{\gamma}}}\cdot\Big(\sum\limits_{|\alpha|=j}\|\partial^\alpha \tilde{v}\|_2^2\Big)^{\frac12}\Bigg\}\\
&&\quad\times
\sqrt{m-3}|u|_{m,3}\frac{\rho(t)^{(m-j-3)}}{((m-j-3)!)^{\frac{1}{\gamma}}}
\Big(\sum_{|\alpha|=m-j}\|z^\gamma\partial^\alpha\partial_zu^{p}\|^2_{\infty}\Big)^{\frac12},
\eeno
where we used Hardy inequality and
\begin{align}
\|\partial^\beta \tilde{v}\partial^{\alpha-\beta}\partial_zu^{p}\|^2_{2}=&\int_{R^2_+}\Big|\frac{\partial^\beta \tilde{v}}{y^\gamma}z^\gamma\partial^{\alpha-\beta}\partial_zu^{p}\Big|^2dxdy\nonumber\\[5pt]
\leq&\|\chi_j(y)y^{1-\gamma}\|^2_\infty\|\partial_y\partial^\beta \tilde{v}\|_2^2\|z^\gamma\partial^{\alpha-\beta}\partial_zu^{p}\|^2_{\infty}\nonumber\\[5pt]
&+\|(1-\chi_j(y))y^{-\gamma}\|^2_\infty\|\partial^\beta \tilde{v}\|_2^2\|z^\gamma\partial^{\alpha-\beta}\partial_zu^{p}\|^2_{\infty}.\nonumber
\end{align}
Here $\chi_j(y)$ is a cut-off function.
Case i) $4\leq j\leq [\frac{m}{2}]$:  $\chi_j(y)$ is chosen as
\begin{align}
\chi_j(y)=\left\{
\begin{array}{ll}
1,\quad 0\leq y\leq 1,\\
0,\quad  y\geq 2,\nonumber
\end{array}
\right.
\end{align}
which implies that
\beno
\frac{C_m^j((j-3)!)^{\frac{1}{\gamma}}((m-j-3)!)^{\frac{1}{\gamma}}\big((j-2)^{\frac1{\gamma}}\|\chi_j(y)y^{1-\gamma}\|_\infty+\|(1-\chi_j(y))y^{-\gamma}\|_\infty \big)}{((m-3)!)^{\frac{1}{\gamma}}\sqrt{m-3}\sqrt{j-2}}\leq C_0j^{-2}.
\eeno
Case ii) $[\frac{m}{2}]\leq j\leq m-3$: take $\chi_j$ such that $\|\chi_j(y)y^{1-\gamma}\|_\infty\approx j^{1-\frac{1}{\gamma}}$, thus
\beno
\frac{C_m^j((j-3)!)^{\frac{1}{\gamma}}((m-j-3)!)^{\frac{1}{\gamma}}\big((j-2)^{\frac1{\gamma}}\|\chi_j(y)y^{1-\gamma}\|_\infty+\|(1-\chi_j(y))y^{-\gamma}\|_\infty \big)}{((m-3)!)^{\frac{1}{\gamma}}\sqrt{m-3}\sqrt{j-2}}\leq C_0.
\eeno
Thus, by discrete convolution inequality and Lemma \ref{e:uniform boundness for approximate solution}, we always have
\beno
J_{21}\leq C_0\|v\|_{Y^3}\|u\|_{Y^3}\sum_{m=3}\frac{\rho(t)^{(m-3)}}{(m-3)!}\Big(\sum_{|\alpha|=m}\|z^\gamma\partial^\alpha\partial_zu^{p}\|^2_{L^\infty}\Big)^{\frac12}
\leq C_0(\|(u,v)\|^2_{Y^3}+\varepsilon^4).
\eeno

Next, we deal with the term $J_{22}.$ Recall that
\beno
J_{22}&=&\sum_{m=4}^{\infty}\frac{\rho(t)^{2(m-3)}\varepsilon^{-\gamma}}{((m-3)!)^{\frac{2}{\gamma}}}\sum_{|\alpha|=m}\sum_{j=m-2}^m\sum_{\beta\leq \alpha,
|\beta|= j}C_\alpha^\beta\big|\langle \partial^\beta \tilde{v}\partial^{\alpha-\beta}\partial_zu^{p},\partial^\alpha u\rangle\big|.
\eeno
Firstly, we deal with the term with $|\beta|=m-2$. Due to $\partial_y v=-\partial_x u$, we have
\beno
&&\sum_{m=4}^{\infty}\frac{\rho(t)^{2(m-3)}\varepsilon^{-\gamma}}{((m-3)!)^{\frac{2}{\gamma}}}\sum_{|\alpha|=m}\sum_{\beta\leq \alpha, |\beta|= m-2}C_\alpha^\beta\big|\langle \partial^\beta \tilde{v}\partial^{\alpha-\beta}\partial_zu^{p},\partial^\alpha u\rangle\big|\\
&&\leq\sum_{m=4}^{\infty}\frac{\rho(t)^{2(m-3)}}{((m-3)!)^{\frac{2}{\gamma}}}\sum_{|\alpha|=m}\sum_{\beta\leq \alpha, |\beta|= m-2}C_\alpha^\beta \Big\{\|\chi_{m-2}(y)y^{1-\gamma}\|_\infty\|\partial_y\partial^\beta \tilde{v}\|_2\\
&&
\quad+\|(1-\chi_{m-2}(y))y^{-\gamma}\|_\infty\|\partial^\beta \tilde{v}\|_2\Big\}\|z^\gamma\partial^{\alpha-\beta}\partial_zu^{p}\|_\infty\|\partial^\alpha u\|_2\\
&&\leq \sum_{m=6}^{\infty} \Big\{\frac{C_m^{m-2}\|\chi_{m-2}(y)y^{1-\gamma}\|_\infty}{(m-3)^{\frac{1}{\gamma}}\sqrt{m-4}\sqrt{m-3}}\frac{\rho(t)^{(m-4)}\sqrt{m-4}}{((m-4)!)^{\frac{1}{\gamma}}}
\Big(\sum_{|\alpha|=m-2}\|\partial_y\partial^\alpha \tilde{v}\|^2_2\Big)^{\frac12}\\
&&\quad+ \frac{C_m^{m-2}\|(1-\chi_{m-2}(y)))y^{-\gamma}\|_\infty}{(m-3)^{\frac{1}{\gamma}}(m-4)^{\frac{1}{\gamma}}\sqrt{m-5}\sqrt{m-3}}\frac{\rho(t)^{(m-5)}
\sqrt{m-5}}{((m-5)!)^{\frac{1}{\gamma}}}\Big(\sum_{|\alpha|=m-2}\|\partial^\alpha \tilde{v}\|^2_2\Big)^{\frac12}\Big\}\\
&&\qquad\times \sqrt{m-3}|u|_{m,3}\Big(\sum\limits_{|\alpha|=2}\|z\partial^\alpha\partial_zu^{p})\|^2_\infty\Big)^{\frac12}
+C_0(\|(u,v)\|^2_{X^3}+\varepsilon^4)\\
&&\leq C_0(\|(u,v)\|^2_{X^3\cap Y^3}+\varepsilon^4),
\eeno
where we take $\chi_{m-2}(y)$ to be a smooth function satisfying $0\leq \chi_{m-2}(y)\leq 1$,
\begin{align}
\chi_{m-2}(y)=\left\{
\begin{array}{ll}
1,\quad 0\leq y\leq m^{-\frac{1}{\gamma}}, \\
0,\quad y\geq 2m^{-\frac{1}{\gamma}},
\end{array}
\right.\nonumber
\end{align}
and use the equality
$$\frac{C_m^{m-2}\|\chi_{m-2}(y)y^{1-\gamma}\|_\infty}{(m-3)^{\frac{1}{\gamma}}\sqrt{m-4}\sqrt{m-3}}\approx \frac{C_m^{m-2}\|(1-\chi_{m-2}(y)))y^{-\gamma}\|_\infty}{(m-3)^{\frac{1}{\gamma}}(m-4)^{\frac{1}{\gamma}}\sqrt{m-5}\sqrt{m-3}}\approx m^{2(1-\frac{1}{\gamma})}.$$

Secondly, for the case of $|\beta|=m$, there holds
\beno
&&\sum_{m=4}^{\infty}\frac{\rho(t)^{2(m-3)}\varepsilon^{-\gamma}}{((m-3)!)^{\frac{2}{\gamma}}}\sum_{|\alpha|=m}\sum_{\beta\leq \alpha, |\beta|= m}C_\alpha^\beta\langle \partial^\beta \tilde{v}\partial^{\alpha-\beta}\partial_zu^{p},\partial^\alpha u\rangle\\
&&=\sum_{m=4}^{\infty}\frac{\rho(t)^{2(m-3)}}{((m-3)!)^{\frac{2}{\gamma}}}\sum_{|\alpha|=m}\Big<\frac{\partial^\alpha \tilde{v}}{y^\gamma}z^\gamma\partial_zu^{p},\partial^\alpha u\Big>\\
&&\leq C_0\sum_{m=4}^{\infty}\frac{(m-2)^{\frac{1}{\gamma}}\|\chi_m(y)y^{1-\gamma}\|\infty}{\sqrt{m-3}\sqrt{m-2}}
\frac{\rho(t)^{(m-2)}\sqrt{m-2}}{((m-2)!)^{\frac{1}{\gamma}}}\Big(\sum_{|\alpha|=m}\|\partial_y\partial^\alpha \tilde{v}\|_2^2\Big)^{\frac12}\sqrt{m-3}|u|_{m,3}\\
&&\quad+C_0\sum_{m=4}^{\infty}\frac{\|(1-\chi_m(y))y^{-\gamma}\|\infty}{m-3}
\frac{\rho(t)^{(m-3)}\sqrt{m-3}}{((m-3)!)^{\frac{1}{\gamma}}}\Big(\sum_{|\alpha|=m}\|\partial^\alpha \tilde{v}\|_2^2\Big)^{\frac12}
\sqrt{m-3}|u|_{m,3}\\
&&\leq  C_0\|v\|_{Y^3}\|u\|_{Y^3},
\eeno
where we take $\chi_m(y)$ as above
and use the inequality
$$\frac{(m-2)^{\frac{1}{\gamma}}\|\chi_m(y)y^{1-\gamma}\|\infty}{\sqrt{m-3}\sqrt{m-2}}\approx \frac{\|(1-\chi_m(y))y^{-\gamma}\|\infty}{m-3}\approx 1.$$

The same argument holds for the case of $|\beta|=m-1$.
Collecting these estimates, we obtain
$$J_{22}\leq C_0(\|(u,v)\|^2_{X^3\cap Y^3}+\varepsilon^4).$$

This together with the estimate of $J_{21}$ gives
\beno
J_2\leq C_0(\|(u,v)\|^2_{X^3\cap Y^3}+\varepsilon^4).
\eeno

Finally, collecting the estimates of $J_1$ and $J_2$, we obtain the estimate of $J$. By the same argument, we can obtain the estimate of $J'$.  The proof is completed.
\endproof

\section{Velocity estimates in Gevrey class and Sobolev space}

\subsection{Velocity estimate in Gevrey class}

We introduce the following energy quantities:
\beno
E_v(t)=\|U\|^2_{X^3}+\varepsilon^4,\quad E_\omega(t)=\|\omega\|^2_{X^2}.
\eeno

\begin{Proposition}[Velocity estimate in Gevrey class]\label{prop:Velocity estimate in Gevrey class}
There exist $\delta>0$ and  $\varepsilon_0>0$ such that for any $\varepsilon\in (0, \varepsilon_0)$, there holds
\beno
&&\frac{1}{2}\frac{d}{dt}\big(\big\|U\big\|^2_{X^3}-\big\|U\big\|^2_2\big)+\lambda\Big(\big\|U\big\|^2_{Y^3}-\big\|U\big\|^2_2\big)
+\frac{\varepsilon^2}{2}\big(\big\|\nabla U\big\|^2_{X^3}-\big\|\nabla U\big\|^2_2\big)\\
&&\leq C_0E_v(t)^{\frac14}\big(E_w(t)+E_v(t)\big)^{\frac14}\big(\|U\|^2_{Y^3}+E_v(t)\big)\\
&&
\quad+C_0(\delta)E_v(t)^{\frac34}\|\partial_yU\|_{X^3}^{\frac12}\big(\varepsilon^2+\|(U,\partial_yU)\|_{X^3}\big)\\
&&
\quad+C_0\big(\|U\|^2_{X^3\cap Y^3}+\varepsilon^4\big)\Big[1+\|(U,\partial_yU)\|_{\overline{H}^4}+\|U\|_{X^3}\|\partial_yU\|_{X^2}\Big]\\
&&\quad+C_0\varepsilon^{2(2-\gamma)}\big(\|u\|^2_{Y^3}+\|\partial_yu\|^2_{Y^2}+\big\|\nabla U\big\|^2_{\overline{H}^3}\big)  +C_0\big(\|U\|^2_{\overline{H}^4}+\|\nabla U\|^2_{\overline{H}^4}+1\big)\Big(\big\|U\big\|^2_{\overline{H}^4}+\varepsilon^4\Big).
\eeno
\end{Proposition}

Proposition \ref{prop:Velocity estimate in Gevrey class} is a direct result of Lemma \ref{e:Gevery derivate estimate on velocity}, Lemma \ref{lem:pressure estimate} and Lemma \ref{lem:low order pressure part}.

\begin{Lemma}\label{e:Gevery derivate estimate on velocity}
There exist $\delta_0,\varepsilon_0>0$ such that for any $\delta\in (0,\delta_0)$ and $\varepsilon\in (0, \varepsilon_0)$, there holds
\beno
&&\frac{1}{2}\frac{d}{dt}\big(\big\|U\big\|^2_{X^3}-\big\|U\big\|^2_2\big)+\lambda\Big(\big\|U\big\|^2_{Y^3}-\big\|U\big\|^2_2\big)
+\frac{\varepsilon^2}{2}\big(\big\|\nabla U\big\|^2_{X^3}-\big\|\nabla U\big\|^2_2\big)\\
&&\leq C_0E_v(t)^{\frac14}\big(E_w(t)+E_v(t)\big)^{\frac14}\big(\|U\|^2_{Y^3}+E_v(t)\big)
+C_0(\delta)E_v(t)^{\frac34}\|\partial_yU\|_{X^3}^{\frac12}\big(\varepsilon^2+\|(U,\partial_yU)\|_{X^3}\big)\\[3pt]
&&\quad
+C_0\big(\|U\|^2_{X^3\cap Y^3}+\varepsilon^4\big)\Big[1+\|(U,\partial_yU)\|_{\overline{H}^4}\Big]\\
&&\quad
+C_0\delta \sum_{m=3}^{\infty}\frac{\rho(t)^{2(m-3)}m}{((m-3)!)^{\frac{2}{\gamma}}}\sum_{|\alpha|=m-1}\big\|\partial^\alpha(\nabla p)\big\|_2^2.
\eeno
\end{Lemma}

\no{\bf Proof.}
Acting $\partial^\alpha$ on both sides of (\ref{eq:error2}), taking $L^2$ inner product with \textcolor[rgb]{0.00,0.00,0.00}{ $\frac{\rho(t)^{2(|\alpha|-3)}}{((|\alpha|-3)!)^{\frac{2}{\gamma}}}\partial^\alpha U$}, and summing up about $|\alpha|=m$ for $m=3,4,\cdots$, we have
\beno
&&\frac{1}{2}\frac{d}{dt}\big(\big\|U\big\|^2_{X^3}-\big\|U\big\|^2_2\big)+\lambda\Big(\big\|U\big\|^2_{Y^3}-\big\|U\big\|^2_2\big)-\varepsilon^2
\sum_{m=3}^{\infty}\frac{\rho(t)^{2(m-3)}}{((m-3)!)^{\frac{2}{\gamma}}}\sum_{|\alpha|=m}\big<\partial^\alpha \triangle U, \partial^\alpha  U\big>\\
&&\leq\Big|\sum_{m=3}^{\infty}\frac{\rho(t)^{2(m-3)}}{((m-3)!)^{\frac{2}{\gamma}}}\sum_{|\alpha|=m}\big<\partial^\alpha(
(\tilde{U}+\tilde{U}^a)\cdot\nabla U), \partial^\alpha U\big>\Big|\\
&&\quad+\Big|\sum_{m=3}^{\infty}\frac{\rho(t)^{2(m-3)}}{((m-3)!)^{\frac{2}{\gamma}}}\sum_{|\alpha|=m}\big<\partial^\alpha(
\tilde{U}\cdot\nabla U^a), \partial^\alpha  U\big>\Big|
+\Big|\sum_{m=3}^{\infty}\frac{\rho(t)^{2(m-3)}}{((m-3)!)^{\frac{2}{\gamma}}}\sum_{|\alpha|=m}\big<\partial^\alpha(
\nabla p), \partial^\alpha  U\big>\Big|\\
&&\quad+\Big|\sum_{m=3}^{\infty}\frac{\rho(t)^{2(m-3)}}{((m-3)!)^{\frac{2}{\gamma}}}\sum_{|\alpha|=m}\big<\partial^\alpha \widetilde{R}, \partial^\alpha U\big>\Big|=\sum_{i=1}^4K_i.
\eeno

{\bf Step 1. Estimate of $K_1.$}
\beno
K_1&\leq&\Big|\sum_{m=3}^{\infty}\frac{\rho(t)^{2(m-3)}}{((m-3)!)^{\frac{2}{\gamma}}}\sum_{|\alpha|=m}\big<\partial^\alpha(
u\partial_x u+\tilde{v}\partial_yu), \partial^\alpha u\big>\Big|\nonumber\\
&&+\Big|\sum_{m=3}^{\infty}\frac{\rho(t)^{2(m-3)}}{((m-3)!)^{\frac{2}{\gamma}}}\sum_{|\alpha|=m}\big<\partial^\alpha(
u\partial_x v+\tilde{v}\partial_yv), \partial^\alpha v\big>\Big|\\
&&+\Big|\sum_{m=3}^{\infty}\frac{\rho(t)^{2(m-3)}}{((m-3)!)^{\frac{2}{\gamma}}}\sum_{|\alpha|=m}\big<\partial^\alpha(
u^a\partial_x u+\tilde{v}^a\partial_yu), \partial^\alpha u\big>\Big|\nonumber\\
&&+\Big|\sum_{m=3}^{\infty}\frac{\rho(t)^{2(m-3)}}{((m-3)!)^{\frac{2}{\gamma}}}\sum_{|\alpha|=m}\big<\partial^\alpha(
u^a\partial_x v+\tilde{v}^a\partial_yv), \partial^\alpha v\big>\Big|=\sum_{i=1}^4K_{1i}.
\eeno

Firstly, by (a) of Lemma \ref{lem:x product estimate} and (a) of Lemma \ref{e:y product estimate}, we obtain
\beno
K_{11}&\leq& C_0(\delta)E_v(t)^{\frac14}(E_w(t)+E_v(t))^{\frac14}\big(\|u\|^2_{Y^3}+E_v(t)\big)
+C_0(\|(U,\partial_yU)\|_{\overline{H}^4}+\varepsilon^2)\|u\|^2_{X^3\cap Y^3}\\[3pt]
&&+C_0(\delta)E_v(t)^{\frac34}\Big[\|\partial_yu\|_{X^3}^{\frac12}(\varepsilon^2+\|(v,\partial_yv)\|_{X^3})
+(\|\partial_yv\|_{X^3}^{\frac12}+\varepsilon)\|\partial_yu\|_{X^3}\Big].\\
K_{12}&\leq&C_0E_v(t)^{\frac14}(E_w(t)+E_v(t))^{\frac14}\big(\|v\|^2_{Y^3}+E_v(t)\big)+C_0(\|(U,\partial_yU)\|_{\overline{H}^4}+\varepsilon^2)\|v\|^2_{X^3\cap Y^3}\\
&&+C_0E_v(t)^{\frac34}
\|\partial_yv\|_{X^3}^{\frac12}(\varepsilon^2+\|(v,\partial_yv)\|_{X^3}),
\eeno
where we used $\partial_yu=\omega+\partial_xv$ and $\partial_yv=-\partial_xu$.

Secondly,  by (b) of Lemma \ref{lem:x product estimate}, (b) of Lemma \ref{e:y product estimate} and (\ref{e:uniform boundness for approximate solution}), we get
\beno
K_{13}+
K_{14}\leq C_0\|U\|^2_{X^3\cap Y^3}.
\eeno

It follows from the estimates on $K_{11}-K_{14}$ that
\beno
K_1&\leq& C_0(\delta)E_v(t)^{\frac14}(E_w(t)+E_v(t))^{\frac14}\big(\|U\|^2_{Y^3}+E_v(t)\big)\\
&&+C_0(\delta)E_v(t)^{\frac34}\|\partial_yU\|_{X^3}^{\frac12}(\varepsilon^2+\|(U,\partial_yU)\|_{X^3})
+C_0\|U\|^2_{X^3\cap Y^3}\Big[1+\|(U,\partial_yU)\|_{\overline{H}^4}\Big].
\eeno

{\bf Step 2. Estimate of $K_2.$}
By Lemma \ref{e:x sublinear product estimate} and (\ref{e:uniform boundness for approximate solution}), we have
\beno
&&\sum_{m=3}^{\infty}\frac{\rho(t)^{2(m-3)}}{((m-3)!)^{\frac{2}{\gamma}}}\sum_{|\alpha|=m}\big|\langle \partial^\alpha(u\partial_xu^a),\partial^\alpha u\rangle\big|\leq C_0\|U\|^2_{X^3\cap Y^3},\\
&&\sum_{m=3}^{\infty}\frac{\rho(t)^{2(m-3)}}{((m-3)!)^{\frac{2}{\gamma}}}\sum_{|\alpha|=m}\big|\langle \partial^\alpha(u\partial_xv^a),\partial^\alpha v\rangle\big|\leq C_0\|U\|^2_{X^3\cap Y^3},\\
&&\sum_{m=3}^{\infty}\frac{\rho(t)^{2(m-3)}}{((m-3)!)^{\frac{2}{\gamma}}}\sum_{|\alpha|=m}\big|\langle \partial^\alpha(\tilde{v}\partial_yv^a),\partial^\alpha v\rangle\big|\leq C_0\Big(\varepsilon^4+\|U\|^2_{X^3\cap Y^3}\Big),\\
&&\sum_{m=3}^{\infty}\frac{\rho(t)^{2(m-3)}}{((m-3)!)^{\frac{2}{\gamma}}}\sum_{|\alpha|=m}\big|\langle \partial^\alpha(\tilde{v}\partial_yu^{e}),\partial^\alpha u\rangle\big|\leq C_0\Big(\varepsilon^4+\|U\|^2_{X^3\cap Y^3}\Big).
\eeno
Thanks to $u^a=u^e+\varepsilon^{1-\gamma}u^p$, it suffices to estimate the remaining term
\beno
&&\sum_{m=3}^{\infty}\frac{\rho(t)^{2(m-3)}\varepsilon^{1-\gamma}}{((m-3)!)^{\frac{2}{\gamma}}}\sum_{|\alpha|=m}\big|\langle \partial^\alpha((v+\varepsilon^2f(t,x)e^{-y})\partial_yu^{p}),\partial^\alpha u\rangle\big|,
\eeno
which is bounded from Lemma \ref{e:singular term estimate} by
\beno
C_0\big(\varepsilon^4+\|U\|^2_{X^3\cap Y^3}\big).
\eeno

Consequently, we obtain
\beno
K_2\leq C_0\big(\varepsilon^4+\|U\|^2_{X^3\cap Y^3}\big).
\eeno

{\bf Step 3. Estimate of $K_3$.}
Note that the following fact holds
\ben\label{eq:y alpha commutate}
\sum_{|\alpha|=m}|[\partial^\alpha,\partial_y]f|\leq C_0\delta m\sum_{|\alpha|=m-1}|\partial^\alpha \partial_yf|,
\een
and we have
\beno
K_3&=&\Big|\sum_{m=3}^{\infty}\frac{\rho(t)^{2(m-3)}}{((m-3)!)^{\frac{2}{\gamma}}}\sum_{|\alpha|=m}\big<\partial^\alpha(\nabla p), \partial^\alpha  U\big>\Big|\\
&\leq&\Big|\sum_{m=3}^{\infty}\frac{\rho(t)^{2(m-3)}}{((m-3)!)^{\frac{2}{\gamma}}}\sum_{|\alpha|=m}\big<\partial^\alpha p, \partial^\alpha  ({\rm div}U)\big>\Big|
+\Big|\sum_{m=3}^{\infty}\frac{\rho(t)^{2(m-3)}}{((m-3)!)^{\frac{2}{\gamma}}}\sum_{|\alpha|=m}\big<[\partial^\alpha,\partial_y] p, \partial^\alpha  v\big>\Big|\\
&&+\Big|\sum_{m=3}^{\infty}\frac{\rho(t)^{2(m-3)}}{((m-3)!)^{\frac{2}{\gamma}}}\sum_{|\alpha|=m}\big<\partial^\alpha p, [\partial^\alpha,\partial_y]  v\big>\Big|+\Big|\sum_{m=3}^{\infty}\frac{\rho(t)^{2(m-3)}}{((m-3)!)^{\frac{2}{\gamma}}}\int_{\partial  R^2_+} \partial_x^{m-1}p \partial_x^{m+1}vdx\Big|\\
&\leq& C_0\delta \sum_{m=3}^{\infty}\frac{\rho(t)^{2(m-3)}m}{((m-3)!)^{\frac{2}{\gamma}}}\sum_{|\alpha|=m-1}\big\|\partial^\alpha(\nabla p)\big\|_2^2
+C_0(\delta)\|v\|^2_{Y^3}+C_0(\delta)\varepsilon^4,
\eeno
where we estimate the boundary term as follows: due to $v|_{y=0}=-\varepsilon^2f(t,x)$ and the assumption (\ref{eq:boundness of f}), we get
\beno
&&\varepsilon^2\Big|\sum_{m=3}^{\infty}\frac{\rho(t)^{2(m-3)}}{((m-3)!)^{\frac{2}{\gamma}}}\int_{\partial  R^2_+} \partial_x^{m-1}p \partial_x^{m+1}fdx\Big|\\
&&=\varepsilon^2\Big|\sum_{m=3}^{\infty}\frac{\rho(t)^{2(m-3)}}{((m-3)!)^{\frac{2}{\gamma}}}\int_{  R^2_+} \partial_y(\partial_x^{m-1}p \partial_x^{m+1}f e^{-y})dxdy\Big|\\
&&\leq C_0\delta \sum_{m=3}^{\infty}\frac{\rho(t)^{2(m-3)}m}{((m-3)!)^{\frac{2}{\gamma}}}\sum_{|\alpha|=m-1}\big\|\partial^\alpha(\nabla p)\big\|_2^2
+C_0(\delta)\varepsilon^4.
\eeno


{\bf Step 4. Estimate of $K_4$.}
By (\ref{e:uniform boundness for error}), we obtain
\beno
K_4=\Big|\sum_{m=3}^{\infty}\frac{\rho(t)^{2(m-3)}}{((m-3)!)^{\frac{2}{\gamma}}}\sum_{|\alpha|=m}\big<\partial^\alpha \tilde{R}, \partial^\alpha  U\big>_2\Big|\leq C_0\varepsilon^4+C_0\|U\|^2_{X^3}.
\eeno

{\bf Step 5. Estimate of dissipative term.}
Integration by parts gives
\beno
E&=&-\sum_{m=3}^{\infty}\frac{\rho(t)^{2(m-3)}}{((m-3)!)^{\frac{2}{\gamma}}}\sum_{|\alpha|=m}\big<\partial^\alpha \triangle U, \partial^\alpha  U\big>\\
&\geq&\Big(\big\|\nabla U\big\|^2_{X^3}-\big\|\nabla U\big\|^2_2\Big)-
\Big|\sum_{m=3}^{\infty}\frac{\rho(t)^{2(m-3)}}{((m-3)!)^{\frac{2}{\gamma}}}\sum_{|\alpha|=m}\big<[\partial^\alpha ,\partial_y] \partial_yU, \partial^\alpha U\big>\Big|\\
&&-
\Big|\sum_{m=3}^{\infty}\frac{\rho(t)^{2(m-3)}}{((m-3)!)^{\frac{2}{\gamma}}}\sum_{|\alpha|=m}\big<\partial^\alpha\partial_yU, [\partial^\alpha,\partial_y]  U\big>\Big|\\
&&-\Big|\sum_{m=3}^{\infty}\frac{\rho(t)^{2(m-3)}}{((m-3)!)^{\frac{2}{\gamma}}}\int_{\partial R^2_+}\partial_y\partial_x^m u \partial_x^mu\Big|-\Big|\sum_{m=3}^{\infty}\frac{\rho(t)^{2(m-3)}}{((m-3)!)^{\frac{2}{\gamma}}}\int_{\partial R^2_+}\partial_y\partial_x^m v \partial_x^mv\Big|.
\eeno
Recalling (\ref{eq:y alpha commutate}) and the boundary conditions of (\ref{eq:error2}), we have
\beno
E&\geq& (1-C_0\delta)\Big(\big\|\nabla U\big\|^2_{X^3}-\big\|\nabla U\big\|^2_2\Big)-\Big|\sum_{m=3}^{\infty}\frac{\rho(t)^{2(m-3)}}{((m-3)!)^{\frac{2}{\gamma}}}\int_{\partial R^2_+}\partial_x^mu(\beta\varepsilon^{-\gamma}\partial_x^m u+\varepsilon\partial^m_x g_0)dx\Big|\\
&&-\Big|\sum_{m=3}^{\infty}\frac{\rho(t)^{2(m-3)}}{((m-3)!)^{\frac{2}{\gamma}}}\int_{\partial R^2_+}\partial_x\partial_x^m u \partial_x^m(\varepsilon^2f)\Big|
\eeno
It follows from Sobolev inequality and Young's inequality that
\beno
E&\geq& \Big(\frac23-C_0\delta\Big)\Big(\big\|\nabla U\big\|^2_{X^3}-\big\|\nabla U\big\|^2_2\Big)-C_0\varepsilon^{-2\gamma}\|
u\big\|^2_{X^3}\\
&&-C_0\varepsilon^2\sum_{m=3}^{\infty}\frac{\rho(t)^{2(m-3)}}{((m-3)!)^{\frac{2}{\gamma}}}
(\big\|\partial_x^{m+1}f\big\|_2^2+\big\|\partial_x^{m}\partial_yf\big\|_2^2).
\eeno
Thus, by taking $\delta$ small and using (\ref{eq:boundness of f}),  we obtain
\beno
\frac{\varepsilon^2}2\Big(\big\|\nabla U\big\|^2_{X^3}-\big\|\nabla U\big\|^2_2\Big)\leq -\varepsilon^2\sum_{m=3}^{\infty}\frac{\rho(t)^{2(m-3)}}{((m-3)!)^{\frac{2}{\gamma}}}\sum_{|\alpha|=m}\big<\partial^\alpha \triangle U, \partial^\alpha  U\big>+C_0\varepsilon^4+C_0\|u\big\|^2_{X^3}.
\eeno

Finally, the proof is completed by collecting these estimates in Step 1-Step 5.
\endproof

\subsection{Pressure estimate in Gevrey class and Sobolev space}
In order to close the estimates of $(u,v)$ of the last subsection, we need to estimate the pressure in Gevrey class. The proof will show why our assumptions are made.
We write
\ben\label{q:some new quantity}
F&=&u^a\partial_xu+\widetilde{v}^a\partial_yu+u\partial_xu^a+\widetilde{v}\partial_yu^a+u\partial_xu+\widetilde{v}\partial_yu,\nonumber\\
G&=&u^a\partial_xv+\widetilde{v}^a\partial_yv+u\partial_xv^a+\widetilde{v}\partial_yv^a
+u\partial_xv+\widetilde{v}\partial_yv.
\een
Thus, the system (\ref{eq:error2}) can be expressed as follows
\begin{eqnarray}
\left \{
\begin {array}{ll}
\partial_tu-\varepsilon^2
\Delta  u+F+\partial_{x}p=R_{1},\\[3pt]
 \partial_tv-\varepsilon^2
\Delta v+G+\partial_yp=R_{2},\\ [3pt]
\partial_xu+\partial_yv=0,\\[3pt]
(\partial_yu, v)(t,x,0)=(\beta\varepsilon^{-\gamma}u+\varepsilon g_0(t,x), -\varepsilon^2f(t,x)),\\[3pt]
(u,v)(0,x,y)=0.
\end{array}
\right.\nonumber
\end{eqnarray}

Taking the divergence operator on both sides, we obtain the equation of the pressure $p$:
\begin{eqnarray}\label{e:pressure equation}
\left \{
\begin {array}{ll}
-\triangle p=(\partial_xF+\partial_yG)-(\partial_xR_{1}+\partial_yR_{2}),\\[3pt]
\partial_yp=-\beta\varepsilon^{2-\gamma}\partial_xu-\varepsilon^3\partial_xg_0+R_{2}+\varepsilon^2\partial_tf-\varepsilon^4\partial_{xx}f
-u^a\partial_xv\quad {\rm on \quad  y=0},
\end{array}
\right.
\end{eqnarray}
where we used $(v+v_a)|_{y=0}=0$ and
$$G(t,x,0)=u^a\partial_xv(t,x,0),\quad v(t,x,0)=-\varepsilon^2f(t,x). $$

\begin{Lemma}[Pressure estimate in Gevrey class]\label{lem:pressure estimate}
\textcolor[rgb]{0.00,0.00,0.00}{Under the assumptions (H)},  we have
\beno
E'&=&\sum_{m=3}^{\infty}\frac{\rho(t)^{2(m-3)}m}{((m-3)!)^{\frac{2}{\gamma}}}\sum_{|\alpha|=m-1}\big\|\partial^\alpha(\nabla p)\big\|_2^2\\
&\leq&C_0\varepsilon^{2(2-\gamma)}\big(\|u\|^2_{Y^3}+\|\partial_yu\|^2_{Y^2}\big)+C_0\|\nabla p\|_2
+C_0(\|U\|_{X^3}\|\partial_yU\|_{X^2}+1)\big(\varepsilon^4+\|U\|^2_{X^3\cap Y^3}\big).
\eeno
\end{Lemma}
{\bf Proof.}
Acting $\partial^\alpha$ on both sides of (\ref{eq:error2}), taking $L^2$ inner product with $\frac{\rho(t)^{2(|\alpha|-3)}}{((|\alpha|-3)!)^{\frac{2}{\gamma}}}\partial^\alpha p$, and summing over $|\alpha|\geq 3$, we have
\beno
&&\sum_{m=3}^{\infty}\frac{\rho(t)^{2(m-3)}m}{((m-3)!)^{\frac{2}{\gamma}}}\sum_{|\alpha|=m-1}\langle -\partial^\alpha\triangle p,\partial^\alpha p\rangle\\
&&=\sum_{m=3}^{\infty}\frac{\rho(t)^{2(m-3)}m}{((m-3)!)^{\frac{2}{\gamma}}}\sum_{|\alpha|=m-1}\langle \partial^\alpha(\partial_xF+\partial_yG)-\partial^\alpha(\partial_xR_1+\partial_yR_2),\partial^\alpha p\rangle
\eeno
Especially, similar to the commutator estimate of dissipative term in Lemma \ref{e:Gevery derivate estimate on velocity} and by  (\ref{e:pressure equation}), the left side is bigger than
\beno
&&\Big(\frac23-C_0\delta\Big)E'-C_0\|\nabla p\|_{\overline{H}^1}^2
+\sum_{m=3}^{\infty}\frac{\rho(t)^{2(m-3)}m}{((m-3)!)^{\frac{2}{\gamma}}}\int_{\partial R^2_+} \partial_x^{m-1}R_2\partial_x^{m-1} pdx\\
&&-\sum_{m=3}^{\infty}\frac{\rho(t)^{2(m-3)}m}{((m-3)!)^{\frac{2}{\gamma}}}\int_{\partial R^2_+} \partial_x^{m-1}(u^a\partial_xv)\partial_x^{m-1} pdx\\
&&-\sum_{m=3}^{\infty}\frac{\rho(t)^{2(m-3)}m}{((m-3)!)^{\frac{2}{\gamma}}}\int_{\partial R^2_+} \partial_x^{m-1}(\beta\varepsilon^{2-\gamma}\partial_xu)\partial_x^{m-1} pdx-C_0\varepsilon^4,
\eeno
where we used integration by parts and
(\ref{eq:boundness of f}).

Moreover, integration by parts, (\ref{eq:y alpha commutate}) and (\ref{e:uniform boundness for error}) yield that
\beno
&&\sum_{m=3}^{\infty}\frac{\rho(t)^{2(m-3)}m}{((m-3)!)^{\frac{2}{\gamma}}}\sum_{|\alpha|=m-1}\langle \partial^\alpha(\partial_xF+\partial_yG)-\partial^\alpha(\partial_xR_1+\partial_yR_2),\partial^\alpha p\rangle\\
&&\leq \frac{1}{10}E'+C_0\varepsilon^4+C_0\sum_{m=3}^{\infty}\frac{\rho(t)^{2(m-3)}m}{((m-3)!)^{\frac{2}{\gamma}}}\sum_{|\alpha|=m-1}\big\| \partial^\alpha (F,G)\big\|_2^2+C_0\|\nabla p\|_{\overline{H}^1}^2\\
&&\quad-\sum_{m=3}^{\infty}\frac{\rho(t)^{2(m-3)}m}{((m-3)!)^{\frac{2}{\gamma}}}\int_{\partial R^2_+} \partial_x^{m-1}(u^a\partial_xv)\partial_x^{m-1} p dx
+\sum_{m=3}^{\infty}\frac{\rho(t)^{2(m-3)}m}{((m-3)!)^{\frac{2}{\gamma}}}\int_{\partial R^2_+} \partial_x^{m-1}R_2\partial_x^{m-1} p dx
\eeno

Canceling the boundary term and by (\ref{e:uniform boundness for error}), we arrive at
\ben\label{eq:pressure middle}
&&\Big(\frac12-C_0\delta\Big)\sum_{m=3}^{\infty}\frac{\rho(t)^{2(m-3)}m}{((m-3)!)^{\frac{2}{\gamma}}}\sum_{|\alpha|=m-1}\big\| \partial^\alpha (\nabla p)\big\|_2^2\nonumber\\
&&\leq C_0\varepsilon^4+C_0\sum_{m=3}^{\infty}\frac{\rho(t)^{2(m-3)}m}{((m-3)!)^{\frac{2}{\gamma}}}\sum_{|\alpha|=m-1}\big\| \partial^\alpha (F,G)\big\|_2^2+C_0\|\nabla p\|_{\overline{H}^1}^2\nonumber\\
&&\quad+\sum_{m=3}^{\infty}\frac{\rho(t)^{2(m-3)}m}{((m-3)!)^{\frac{2}{\gamma}}}\Big|\int_{\partial R^2_+} \partial_x^{m-1}(\beta\varepsilon^{2-\gamma}\partial_xu)\partial_x^{m-1} pdx\Big|.
\een

Recalling (\ref{q:some new quantity}), we have
\ben\label{eq:F G ESTIMATE}
\sum_{m=3}^{\infty}\frac{\rho(t)^{2(m-3)}m}{((m-3)!)^{\frac{2}{\gamma}}}\sum_{|\alpha|=m-1}\big\| \partial^\alpha (F,G)\big\|_2^2
\leq C_0(\|U\|_{X^3}\|\partial_yU\|_{X^2}+1)\big(\varepsilon^4+\|U\|^2_{X^3\cap Y^3}\big).
\een
In fact, it is sufficient to handle the terms $\tilde{v}^a\partial_yu$, $\tilde{v}\partial_yu^a$ and $\tilde{v}\partial_yu$, and other terms are similar.

Firstly, by the assumption (H1), (\ref{e:uniform boundness for approximate solution}) and $\tilde{v}_a|_{y=0}=0$, and following the proof of Lemma \ref{e:x sublinear product estimate}, we have
\beno
&&\sum_{m=3}^{\infty}\frac{\rho(t)^{2(m-3)}m}{((m-3)!)^{\frac{2}{\gamma}}}\sum_{|\alpha|=m-1}\big\| \partial^\alpha (\tilde{v}^a\partial_yu)\big\|_2^2\\
&&\leq \sum_{m=3}^{\infty}\frac{\rho(t)^{2(m-3)}m}{((m-3)!)^{\frac{2}{\gamma}}}\sum_{|\alpha|=m-1}\sum_{\beta\leq \alpha}(C_\alpha^\beta)^2\Big\| \frac{\partial^\beta \tilde{v}^a}{\psi}({\psi}\partial^{\alpha-\beta}\partial_yu)\Big\|_2^2
\leq C_0\|u\|^2_{X^3\cap Y^3}.
\eeno
Next, by the assumption (H1), (\ref{e:uniform boundness for approximate solution}), Hardy inequality and $\tilde{v}|_{y=0}=0$, similar argument as in Lemma \ref{e:x sublinear product estimate} gives
\beno
&&\sum_{m=3}^{\infty}\frac{\rho(t)^{2(m-3)}m}{((m-3)!)^{\frac{2}{\gamma}}}\sum_{|\alpha|=m-1}\big\| \partial^\alpha (\tilde{v}\partial_yu^a)\big\|_2^2\\
&&\leq \sum_{m=3}^{\infty}\frac{\rho(t)^{2(m-3)}m}{((m-3)!)^{\frac{2}{\gamma}}}\sum_{|\alpha|=m-1}\sum_{\beta\leq \alpha}(C_\alpha^\beta)^2\Big\| \frac{\partial^\beta \tilde{v}}{\psi}({\psi}\partial^{\alpha-\beta}\partial_y)u^a\Big\|_2^2
\leq C_0\Big(\varepsilon^4+\|U\|^2_{X^3\cap Y^3}\Big).
\eeno
Thirdly, by Sobolev inequality, Hardy inequality and $\tilde{v}|_{y=0}=0$, similar argument as in Lemma \ref{e:y product estimate} gives
\beno
&&\sum_{m=3}^{\infty}\frac{\rho(t)^{2(m-3)}m}{((m-3)!)^{\frac{2}{\gamma}}}\sum_{|\alpha|=m-1}\big\| \partial^\alpha (\tilde{v}\partial_yu)\big\|_2^2\\
&&\leq \sum_{m=3}^{\infty}\frac{\rho(t)^{2(m-3)}m}{((m-3)!)^{\frac{2}{\gamma}}}\sum_{|\alpha|=m-1}\sum_{\beta\leq \alpha}(C_\alpha^\beta)^2\Big\| \frac{\partial^\beta \tilde{v}}{\psi}({\psi}\partial^{\alpha-\beta}\partial_yu)\Big\|_2^2
\leq C_0\|U\|_{X^3}\|\partial_yU\|_{X^2}\|U\|^2_{X^3\cap Y^3}.
\eeno

Finally, there holds
\ben\label{eq:boundary of pressure}
&&\sum_{m=3}^{\infty}\frac{\rho(t)^{2(m-3)}m}{((m-3)!)^{\frac{2}{\gamma}}}\Big|\int_{\partial R^2_+} \partial_x^{m-1}(\beta\varepsilon^{2-\gamma}\partial_xu)\partial_x^{m-1} pdx\Big|\nonumber\\
&&\leq\sum_{m=3}^{\infty}\frac{\rho(t)^{2(m-3)}m}{((m-3)!)^{\frac{2}{\gamma}}}\Big|\int_{ R^2_+} \big(\beta\varepsilon^{2-\gamma}\partial_x^{m}u\partial_x^{m-1}\partial_y p+\beta\varepsilon^{2-\gamma}\partial_x^{m}\partial_{y}u\partial_x^{m-1} p\big)dxdy\Big|\nonumber\\
&&\leq \frac{1}{10}\sum_{m=3}^{\infty}\frac{\rho(t)^{2(m-3)}m}{((m-3)!)^{\frac{2}{\gamma}}}\sum_{|\alpha|=m-1}\big\| \partial^\alpha (\nabla p)\big\|_2^2
+C_0\varepsilon^{2(2-\gamma)}\big(\|u\|^2_{Y^3}+\|\partial_yu\|^2_{Y^2}\big).
\een

Putting (\ref{eq:F G ESTIMATE}) and (\ref{eq:boundary of pressure}) into (\ref{eq:pressure middle}) and
taking $\delta$ small, the proof is completed.
\endproof
\smallskip

Next we deal with the pressure estimate in $L^2$, which is similar to Lemma \ref{lem:pressure estimate} but more delicate. Especially, it is necessary to use the condition $f(t,x)=\partial_x\overline{f}$ of the Assumption (H3).

\begin{Lemma}[Pressure estimate in $L^2$]\label{lem:low order pressure part}
\textcolor[rgb]{0.00,0.00,0.00}{Under the  assumptions (H)},  we have
\begin{align}
\big\|\nabla p\big\|^2_2\leq C_0\big(\|U\|^2_{\overline{H}^4}+\|\nabla U\|^2_{\overline{H}^4}+1\big)\Big(\big\|U\big\|^2_{\overline{H}^4}+\varepsilon^4\Big)+C_0\varepsilon^{2(2-\gamma)}\big\|\nabla U\big\|^2_{\overline{H}^3}.\nonumber
\end{align}
\end{Lemma}
{\bf Proof:}
As in Lemma \ref{lem:pressure estimate}, taking $L^2$ inner product with $p$ in (\ref{e:pressure equation}), we have
\begin{align}\label{e:equality on Euler pressure}
-\big<\triangle p,p\big>
=\big<(\partial_xF+\partial_yG)-(\partial_xR_1+\partial_yR_2),p\big>
\end{align}
Integrating by parts and using the boundary value of $\partial_y p$ in $(\ref{e:pressure equation})_2$ yield that
\beno
-\big<\triangle p,p\big>=\|\nabla p\|_2^2+\int_{\partial\R^2_+}(-\beta\varepsilon^{2-\gamma}\partial_xu-\varepsilon^3\partial_xg_0+R_{2}+\varepsilon^2\partial_tf-\varepsilon^4\partial_{xx}f
-u^a\partial_xv)pdx.
\eeno
By the commutator estimate (\ref{eq:y alpha commutate}) and (\ref{eq:boundness of f}) of (H3),
 the left side of (\ref{e:equality on Euler pressure}) is bigger than
\ben\label{eq:pressure lower order left}
(\frac23-C_0\delta)\big\|\nabla p\big\|^2_2+\int_{\partial R^2_+}pR_2dx-\int_{\partial R^2_+}p(u^a\partial_xv)dx
-C_0\varepsilon^4-C_0\beta^2\varepsilon^{2(2-\gamma)}\big\|\nabla U\big\|^2_2,
\een
where we used that $f(t,x)=\partial_x\overline{f}$ and
\textcolor[rgb]{0.00,0.00,0.00}{\beno
\int_{\partial\R^2_+}\varepsilon^2\partial_tf pdx&=&\int_{\partial\R^2_+}\varepsilon^2\partial_t\partial_x\overline{f} pdx=-\int_{\R^2_+}\varepsilon^2\partial_t\partial_x\overline{f} \partial_y(e^{-y}p)dxdy\\
&\leq& C_0\varepsilon^2\big(\|\partial_t\partial_x\overline{f}\|_{L^2_x}\|\partial_y p\|_{L^2}+\|\partial_t \bar{f}\|_{L^2_x}\|\partial_xp\|_2\big)\leq C_0\varepsilon^4+\frac{1}{10}\big\|\nabla p\big\|^2_2.
\eeno}

Recall that $G(t,x,0)=u^a\partial_xv(t,x,0)$. By (\ref{e:uniform boundness for error}) and integration by parts, the right hand side of (\ref{e:equality on Euler pressure}) can be controlled by
\ben\label{eq:pressure lower order right}
\frac12\big\|\nabla p\big\|^2_2+C_0\big\|(F,G)\big\|^2_2+C_0\varepsilon^4-\int_{\partial R^2_+}p(u^a\partial_xv)dx+\int_{\partial R^2_+}pR_2dx.
\een
Thus, canceling the boundary term between (\ref{eq:pressure lower order left}) and (\ref{eq:pressure lower order right}), taking $\delta$ small, we arrive at
\begin{align}
\big\|\nabla p\big\|^2_2\leq C_0\big\|(F,G)\big\|^2_2+C_0\varepsilon^4+C_0\varepsilon^{2(2-\gamma)}\big\|\nabla U\big\|^2_2.\nonumber
\end{align}

At last, we are aimed to estimate the term $\big\|(F,G)\big\|^2_2$.
Recalling (\ref{q:some new quantity}), by (\ref{e:uniform boundness for approximate solution}), $\tilde{v}|_{y=0}=0$ and
$(v^a-\varepsilon^2f(t,x)e^{-y})|_{y=0}=0,$ we can obtain the following estimate directly
\beno
\big\|(F,G)\big\|^2_2\leq C_0\big(\|U\|^2_2+\|\nabla U\|^2_2+1\big)\Big(\big\|U\big\|^2_{\overline{H}^1}+\varepsilon^4\Big).\nonumber
\eeno
For example, we only deal with the terms $\tilde{v}^a\partial_yu$, $\tilde{v}\partial_yu^a$ and $\tilde{v}\partial_yu$.

Firstly, by (H1) and (\ref{e:uniform boundness for approximate solution}), we have
\beno
\big\|\tilde{v}^a\partial_yu\big\|_2=\Big\|\frac{\tilde{v}_a}{\psi}Zu\Big\|_2\leq C_0\|u\|_{\overline{H}^1}.
\eeno
Moreover, by Hardy inequality and divergence-free property, we deduce that
\beno
\big\|\tilde{v}\partial_yu^a\big\|_2=\Big\|\frac{\tilde{v}}{\psi}Zu^a\Big\|_2\leq C_0\Big(\big\|U\big\|_{\overline{H}^1}+\varepsilon^2\Big).
\eeno
Finally, by Sobolev embedding and divergence free condition, we have
\beno
\big\|\tilde{v}\partial_yu\big\|_2=\big\|\frac{\tilde{v}}{\psi}Zu\big\|_2&\leq& C_0\Big(\big\|\partial_y\tilde{v}\big\|_\infty+\big\|\tilde{v}\big\|_\infty\Big)\|u\|_{\overline{H}^1}
\\
&\leq& C_0\big(\|\nabla U\|_{\overline{H}^4}+\| U\|_{\overline{H}^4}+\varepsilon^2\big)\big(\big\|U\big\|_{\overline{H}^4}+\varepsilon^2\big).
\eeno

Hence, the proof is completed.
\endproof

\subsection{Velocity estimate in Sobolev space}

In this subsection, our goal is to prove the following proposition.

\begin{Proposition}[Velocity estimate in $L^2$]\label{prop:Velocity estimate in Sobolev space}
There exist $\delta>0$ and $\varepsilon_0>0$ such that for any $\varepsilon\in (0,\varepsilon_0)$, we have
\begin{align}
\frac{1}{2}\frac{d}{dt}\big\|U\big\|^2_2+\frac{\varepsilon^2}{2}\big\|\nabla U\big\|^2_2
\leq C_0(\delta)\big(\|U\|^2_{\overline{H}^4}+\|\nabla U\|^2_{\overline{H}^4}+1\big)\Big(\big\|U\big\|^2_{\overline{H}^4}+\varepsilon^4\Big)+C_0\|U\|^2_{\overline{H}^5}.\nonumber
\end{align}
\end{Proposition}

The proof of Proposition \ref{prop:Velocity estimate in Sobolev space} is an immediate corollary of the following Lemma \ref{e:conormal estimate on velocity} and Lemma \ref{lem:low order pressure part}.

\begin{Lemma}\label{e:conormal estimate on velocity}
There exist $\delta_0,\varepsilon_0>0$ such that for any $\delta\in (0,\delta_0)$ and $\varepsilon\in (0, \varepsilon_0)$, we have
\begin{align}
\frac{1}{2}\frac{d}{dt}\big\|U\big\|^2_2+\frac{\varepsilon^2}{2}\big\|\nabla U\big\|^2_2
\leq C_0(\delta)\Big(\big\|U\big\|^2_2+\varepsilon^4\Big)+C_0\|U\|^2_{\overline{H}^1}
+C_0\varepsilon^2\big\|\nabla p\big\|^2_2.\nonumber
\end{align}
\end{Lemma}
{\bf Proof:}
Taking $L^2(R^2_+)$ inner product with $U$ in (\ref{eq:error2}), we arrive at
\begin{align}
\frac{1}{2}\frac{d}{dt}\big\|U\big\|^2_2
\leq& \big|\big<u^{a}\partial_xU+\tilde{v}^a\partial_yU,U\big>\big|
+\big|\big<u\partial_xU^a+\tilde{v}\partial_yU^a,U\big>\big|+\big|\big<u\partial_xU+\tilde{v}\partial_yU,U\big>\big|\nonumber\\
&+\big|\big<\nabla p,U\big>\big|+\big|\big<\widetilde{R},U\big>\big|
=\sum_{i=1}^5N_i.\nonumber
\end{align}

{\bf Step 1. Estimate of $N_1$.}
Note that $\tilde{v}^a|_{y=0}=0$,
it follows from integration by parts, divergence-free condition and (\ref{eq:boundness of f}) that
\beno
N_{1}\leq \Big|\int_{R^2_+}\varepsilon^2f(t,x)e^{-y}|U|^2dxdy\Big|
\leq C_0\varepsilon^2\big\|U\big\|^2_2.
\eeno

{\bf Step 2. Estimate of $I_2$.} Note the fact $\tilde{v}|_{y=0}=0$.
By (\ref{e:uniform boundness for approximate solution}), divergence free condition and Hardy inequality, we get
\beno
N_2&\leq& C_0(\big\|U\big\|^2_2+\varepsilon^4)+\varepsilon^{1-\gamma}\big|\big<\tilde{v}\partial_yu^{p},u\big>\big|\\
&\leq& C_0(\big\|U\big\|^2_2+\varepsilon^4)+\varepsilon^{1-\gamma}\Big|\Big<\frac{1}{y}\int_0^y\partial_{y'}
\tilde{v}dy'z\partial_zu^{p},u\Big>\Big|
\leq C_0\big(\big\|U\big\|^2_{\overline{H}^1}+\varepsilon^4\big),
\eeno
where we used $\|z\partial_zu^{p}\|_{\infty}\leq C_0.$

{\bf Step 3. Estimate of $N_3$.}
By integration by parts, divergence free condition, (\ref{eq:boundness of f}) and $\tilde{v}|_{y=0}=0$, we obtain
\begin{align}
N_{3}\leq \Big|\int_{R^2_+}\varepsilon^2f(t,x)e^{-y}|U|^2dxdy\Big|\leq C_0\varepsilon^2\big\|U\big\|^2_2.\nonumber
\end{align}

{\bf Step 4. Estimate of $N_4$ and $N_5$.}
By integration by parts, divergence free condition, (\ref{eq:boundness of f}) and  $v(t,x,0)=-\varepsilon^2f(t,x)$, we have
\beno
N_4\leq \varepsilon^2\Big|\int_{\partial R^2_+}pfdx\Big|
\leq C_0\varepsilon^2\big\|\nabla p\big\|_2.
\eeno
By (\ref{e:uniform boundness for error}), it is easy to get
$$N_5\leq C_0\big(\big\|U\big\|^2_2+\varepsilon^4\big). $$

{\bf Step 5. Estimate of dissipative term.}
By integration by parts and the boundary condition (\ref{eq:boundary condition of u}), we have
\beno
-\big<\triangle U,U\big>&=&\big\|\nabla U\big\|^2_2
+\beta\varepsilon^{-\gamma}\int_{\partial R^2_+}u^2dx+\varepsilon\int_{\partial R^2_+}ug_0dx-\varepsilon^2\int_{\partial R^2_+}\partial_yv fdx\\
&\geq& \Big(\frac23-C_0\delta\Big)\big\|\nabla U\big\|^2_2-C_0\beta^2\varepsilon^{-2\gamma}\big\| U\big\|^2_2-C_0\varepsilon^2.\nonumber
\eeno
Thus, choosing $\delta_0$ small such that when $\delta\leq \delta_0$,  we arrive at
\beno
-\varepsilon^2\big<\triangle U,U\big>_2\geq \frac{\varepsilon^2}{2}\big\|\nabla U\big\|^2_2-C_0\varepsilon^{2(1-\gamma)}\big\| U\big\|^2_2-C_0\varepsilon^4
\eeno

Collecting these estimates of Step 1-Step 5, we finish the proof.
\endproof

\section{Vorticity estimates in  Gevrey class and Sobolev space}

\subsection{Vorticity estimate in  Gevrey class}

Recall the equation (\ref{eq:modified vorticity equation}) of $\eta$, and we are going to prove the following proposition.

\begin{Proposition}[Vorticity estimate in Gevrey class]\label{prop:Vorticity estimate in Gevrey class}
There exist $\delta>0$ and $\varepsilon_0>0$ such that for any $\varepsilon\in (0, \varepsilon_0)$, there holds
\beno
&&\frac{1}{2}\frac{d}{dt}\big(\big\|\eta\big\|^2_{X^2}-\|\eta\|^2_2\big)+\lambda\Big(\big\|\eta\big\|^2_{Y^2}-\|\eta\|^2_2\Big)
+\frac{\varepsilon^2}{2}\Big(\big\|\nabla \eta\big\|^2_{X^2}-\|\nabla\eta\|^2_2\Big)\\[3pt]
&&\leq C_0\|\eta\|^2_{Y^2\cap X^2}+C_0\|u\|_{X^3}^{\frac12}\| \partial_yu\|_{X^2}^{\frac12}\|\eta\|^2_{X^2\cap Y^2}\\[3pt]
&&\quad +C_0\|u\|_{X^3}\| \eta\|_{X^2}^{\frac32}\| \partial_y\eta\|_{X^2}^{\frac12}
+C_0\delta\varepsilon^{2(2-2\gamma)}\big(\big\|\nabla U\big\|^2_{\overline{H}^3}+\|u\|^2_{Y^3}+\|\partial_yu\|^2_{Y^2}\big) \\[3pt]
&&\quad+C_0(\delta)\|(\widetilde{v},\partial_y\widetilde{v})\|^{\frac12}_{\overline{H}^3}
\|(\partial_y\widetilde{v},\partial_{yy}\widetilde{v})\|^{\frac12}_{\overline{H}^3}
\| \eta\|^2_{X^2\cap Y^2}\\[3pt]
&&\quad+C_0\| \widetilde{v}\|^{\frac12}_{X^2}\|\partial_y\widetilde{v}\|_{X^2}^{\frac12}
\|\partial_y\eta\|_{X^2}\|\eta\|_{X^2}+C_0(\delta)\| \eta\|^{\frac32}_{X^2}\| \partial_y\eta\|_{X^2}^{\frac12}\|(\widetilde{v},\partial_y\widetilde{v})\|_{X^2}\\
&&\quad+C_0(\varepsilon^{-2}\|U\|^2_{X^2\cap Y^2}+\varepsilon^2)
+C_0\delta\varepsilon^{-2\gamma}\big(\|U\|^2_{\overline{H}^4}+\|\nabla U\|^2_{\overline{H}^4}\big)\big(\big\|U\big\|^2_{\overline{H}^4}+\varepsilon^4\big)\\[3pt]
&&\quad+C_0\delta\varepsilon^{-2\gamma}(\|U\|_{X^3}\|\partial_yU\|_{X^2}+1)\big(\varepsilon^4+\|U\|^2_{X^3\cap Y^3}\big)+\frac{\varepsilon^2}{10}\|\partial_y \eta\|^2_{\overline{H}^1}.
\eeno
\end{Proposition}

We need the following lemma, which is similar to Lemma \ref{lem:x product estimate}-\ref{e:singular term estimate}.

\begin{Lemma}\label{e:x product estimate'}
Let
\beno
A'&=&\sum_{m=2}^{\infty}\frac{\rho(t)^{2(m-2)}}{((m-2)!)^{\frac{2}{\gamma}}}\sum_{|\alpha|=m}\big|\langle \partial^\alpha(u\partial_xv),\partial^\alpha v\rangle\big|.
\eeno
Then there hold
\begin{itemize}
\item[(a)]
\beno
 A'\leq C_0\|u\|_{X^3}^{\frac12}\| \partial_yu\|_{X^2}^{\frac12}\|v\|^2_{X^2\cap Y^2}
 +C_0\|u\|_{X^3}\| v\|_{X^2}^{\frac32}\| \partial_yv\|_{X^2}^{\frac12};
\eeno
\item[(b)]
\beno
 A'&\leq&
C_0\|u\|_{X^3}^{\frac12}\|\partial_yu\|_{X^3}^{\frac12}\|v\|^2_{X^2\cap Y^2}.
\eeno
\end{itemize}
\end{Lemma}

\no{\bf Proof.}
The proof of $(b)$ is the same as that in Lemme \ref{lem:x product estimate}, and we omit it here. The proof of $(a)$  is also similar to that in Lemma \ref{lem:x product estimate} except the estimate of $I_{22}$, and we sketch it.

Firstly, using Sobolev inequality and integrating by parts, we have
\beno
\sum_{m=2}^{\infty}\frac{\rho(t)^{2(m-2)}}{((m-2)!)^{\frac{2}{\gamma}}}\sum_{|\alpha|=m}\big|\langle u\partial_x\partial^\alpha v,\partial^\alpha v\rangle\big|
\leq C_0\|u\|_{\overline{H}^2}^{\frac12}\|\partial_yu\|_{\overline{H}^2}^{\frac12}\|v\|_{X^2}^2,
\eeno
which gives
\beno
A'\leq  I'+C_0\|u\|_{\overline{H}^2}^{\frac12}\|\partial_yu\|_{\overline{H}^2}^{\frac12}\|v\|_{X^2}^2,
\eeno
where
\beno
I'=\sum_{m=2}^{\infty}\frac{\rho(t)^{2(m-2)}}{((m-2)!)^{\frac{2}{\gamma}}}\sum_{|\alpha|=m}\sum_{0<\beta\leq\alpha}C_\alpha^\beta\big|\langle \partial^\beta u\partial_x\partial^{\alpha-\beta}v,\partial^\alpha v\rangle\big|,
\eeno
and $I'$ could be estimated as follows
\beno
I'
&\leq&\sum_{m=2}^{\infty}\frac{\rho(t)^{2(m-2)}}{((m-2)!)^{\frac{2}{\gamma}}}\sum_{|\alpha|=m}\sum_{1\leq|\beta|\leq 2,\beta\leq\alpha}C_\alpha^\beta\|\partial^\beta u\partial_x\partial^{\alpha-\beta}v\|_{2}\|\partial^\alpha v\|_{2}\\
&&+\sum_{m=5}^{\infty}\frac{\rho(t)^{2(m-2)}}{((m-2)!)^{\frac{2}{\gamma}}}\sum_{|\alpha|=m}\sum_{j=3}^{m-2}\sum_{|\beta|=j,\beta\leq\alpha}C_\alpha^\beta\|\partial^\beta u\partial_x\partial^{\alpha-\beta}v\|_{2}\|\partial^\alpha v\|_{2}\\
&&+\sum_{m=2}^{\infty}\frac{\rho(t)^{2(m-2)}}{((m-2)!)^{\frac{2}{\gamma}}}\sum_{|\alpha|=m}\sum_{m-1\leq|\beta|\leq m,\beta\leq\alpha}C_\alpha^\beta\|\partial^\beta u\partial_x\partial^{\alpha-\beta}v\|_{2}\|\partial^\alpha v\|_{2}
\doteq\sum_{i=1}^3I'_{i}
\eeno

Next we handle them term by term.\smallskip

{\bf Step 1. Estimate of $I'_1$.}
The same argument as $I_1$ in Lemma \ref{lem:x product estimate} gives
\beno
I'_1\leq C_0\|u\|_{\overline{H}^3}^{\frac12}\|\partial_yu\|_{\overline{H}^3}^{\frac12}\|v\|^2_{Y^2\cap X^2}.
\eeno

{\bf Step 2. Estimate of $I'_2$.} Similar to $I_2$ in Lemma \ref{lem:x product estimate}, we decompose $I'_2$ into two parts
\beno
I'_2&=&\sum_{m=6}^{\infty}\frac{\rho(t)^{2(m-2)}}{((m-2)!)^{\frac{2}{\gamma}}}\sum_{|\alpha|=m}\sum_{j=3}^{[\frac{m}{2}]}\sum_{|\beta|=j,\beta\leq\alpha}
C_\alpha^\beta\|\partial^\beta u\partial_x\partial^{\alpha-\beta}v\|_{2}\|\partial^\alpha v\|_{2}\\
&&+\sum_{m=5}^{\infty}\frac{\rho(t)^{2(m-2)}}{((m-2)!)^{\frac{2}{\gamma}}}\sum_{|\alpha|=m}\sum_{j=[\frac{m}{2}]+1}^{m-2}\sum_{|\beta|=j,\beta\leq\alpha}C_\alpha^\beta\|\partial^\beta u\partial_x\partial^{\alpha-\beta}v\|_{2}\|\partial^\alpha v\|_{2}=I'_{21}+I'_{22}.
\eeno
The same argument as $I_{21}$ in Lemma \ref{lem:x product estimate} gives
$$I'_{21}\leq C_0\| u\|^{\frac12}_{X^2}\| \partial_yu\|_{X^2}^{\frac12}\|v\|^2_{Y^2}.$$
Then $I'_{22}$ can be bounded by
\beno
&&\sum_{m=5}^{\infty}\frac{\rho(t)^{2(m-2)}}{((m-2)!)^{\frac{2}{\gamma}}}\sum_{|\alpha|=m}\sum_{j=[\frac{m}{2}]+1}^{m-2}
\sum_{|\beta|=j,\beta\leq\alpha}C_\alpha^\beta\|\partial^\beta u\|_2\|\partial_x\partial^{\alpha-\beta}v\|^{\frac12}_{2}\|\partial_x\partial_y\partial^{\alpha-\beta}v\|^{\frac12}_{2}\|\partial^\alpha v\|_{2}\\
&&+\sum_{m=5}^{\infty}\frac{\rho(t)^{2(m-2)}}{((m-2)!)^{\frac{2}{\gamma}}}\sum_{|\alpha|=m}\sum_{j=[\frac{m}{2}]+1}^{m-2}
\sum_{|\beta|=j,\beta\leq\alpha}C_\alpha^\beta\|\partial^\beta u\|_2\|\partial_{xx}\partial^{\alpha-\beta}v\|^{\frac12}_{2}\|\partial_{xx}\partial_y\partial^{\alpha-\beta}v\|^{\frac12}_{2}\|\partial^\alpha v\|_{2}.
\eeno
We only estimate the first term, which can be controlled from similar arguments as in (\ref{eq:middle term-Young inequ}) by
\beno
&& \sum_{m=5}^{\infty}\frac{\rho(t)^{(m-2)}}{((m-2)!)^{\frac{1}{\gamma}}}\sum_{j=[\frac{m}{2}]+1}^{m-2}
C_m^j\Big(\sum_{|\alpha|=m-j+1}\|\partial^\alpha v\|_2\|\partial_y \partial^\alpha v\|_2\Big)^{\frac12}\Big(\sum_{|\alpha|=j}
\|\partial^\alpha u\|^2_{2}\Big)^{\frac12}|v|_{m,2}\\
&&\leq C_0\sum_{m=6}^{\infty}\sum_{j=[\frac{m}{2}]+1}^{m-2}
\frac{C_m^j((m-j-1)!)^{\frac{1}{\gamma}}((j-3)!)^{\frac{1}{\gamma}}}{((m-2)!)^{\frac{1}{\gamma}}}\\
&&\quad \times\frac{\rho(t)^{\frac{(m-j-1)}{2}}}{((m-j-1)!)^{\frac{1}{2\gamma}}}\Big(\sum_{|\alpha|=m-j+1}\|\partial_y \partial^\alpha v\|^2_2\Big)^{\frac12}|v|^{\frac12}_{m-j+1,2}|u|_{j,3}|v|_{m,2}\\
&&\quad+C_0\|v\|^{\frac12}_{\overline{H}^3}\|\partial_yv\|^{\frac12}_{\overline{H}^3}\|u\|_{X^3} \|v\|_{X^2}\\
&&\leq C_0\sum_{m=6}^{\infty}\sum_{j=[\frac{m}{2}]+1}^{m-2}
\frac{j^{-\frac{1}{\gamma}}\rho(t)^{\frac{(m-j-1)}{2}}}{((m-j-1)!)^{\frac{1}{2\gamma}}}\Big(\sum_{|\alpha|=m-j+1}\|\partial_y \partial^\alpha v\|^2_2\Big)^{\frac12}|v|^{\frac12}_{m-j+1,2}|u|_{j,3}|v|_{m,2}\\
&&\quad+C_0\|v\|^{\frac12}_{\overline{H}^3}\|\partial_yv\|^{\frac12}_{\overline{H}^3}\|u\|_{X^3} \|v\|_{X^2}\\
&&\leq C_0\| v\|^{\frac32}_{X^2}\| \partial_yv\|_{X^2}^{\frac12}\|u\|_{X^3},
\eeno
where we used
$$\frac{C_m^j((m-j-1)!)^{\frac{1}{\gamma}}((j-3)!)^{\frac{1}{\gamma}}}{((m-2)!)^{\frac{1}{\gamma}}}\leq C_0\Big(C_{m-2}^{j-1}\Big)^{1-\frac{1}{\gamma}}j^{-\frac{1}{\gamma}}$$
for $j\in \{[\frac{m}{2}]+1,\cdots, m-2\}$.

Finally, we obtain
\beno
I'_2\leq   C_0\| v\|^{\frac32}_{X^2}\| \partial_yv\|_{X^2}^{\frac12}\|u\|_{X^3}.
\eeno

{\bf Step 3. Estimate of $I'_3$.} Following the proof of $I_3$ in Lemma \ref{lem:x product estimate} line by line, we get
\beno
I'_3\leq C_0\|v\|_{\overline{H}^3}^{\frac12}\|\partial_yv\|_{\overline{H}^3}^{\frac12}\|u\|_{X^2}\|v\|_{X^2}.
\eeno
\endproof

\no{\bf Proof of Proposition \ref{prop:Vorticity estimate in Gevrey class}.}
Acting $\partial^\alpha$ on both sides of (\ref{eq:modified vorticity equation}), multiplying $\frac{\rho(t)^{2(|\alpha|-2)}}{((|\alpha|-2)!)^{\frac{2}{\gamma}}}\partial^\alpha \eta$ and summing over $|\alpha|\geq 2$, we arrive at
\beno
&&\frac{1}{2}\frac{d}{dt}\big(\big\|\eta\big\|^2_{X^2}-\|\eta\|^2_2\big)+\lambda\Big(\big\|\eta\big\|^2_{Y^2}-\|\eta\|^2_2\Big)-\varepsilon^2
\sum_{m=2}^{\infty}\frac{\rho(t)^{2(m-2)}}{(m-2)!)^{\frac{2}{\gamma}}}\sum_{|\alpha|=m}\big<\partial^\alpha \triangle \eta, \partial^\alpha  \eta\big>\\
&&\leq \Big|\sum_{m=2}^{\infty}\frac{\rho(t)^{2(m-2)}}{((m-2)!)^{\frac{2}{\gamma}}}\sum_{|\alpha|=m}\big<\partial^\alpha(
(u^a\partial_x\eta+\widetilde{v}^a\partial_y\eta), \partial^\alpha \eta\big>\Big|\\
&&\quad+\Big|\sum_{m=2}^{\infty}\frac{\rho(t)^{2(m-2)}}{((m-2)!)^{\frac{2}{\gamma}}}\sum_{|\alpha|=m}\big<\partial^\alpha(
u\partial_x\eta^a+\widetilde{v}\partial_y\eta^a), \partial^\alpha  \eta\big>\Big|\\
&&\quad+\Big|\sum_{m=2}^{\infty}\frac{\rho(t)^{2(m-2)}}{((m-2)!)^{\frac{2}{\gamma}}}\sum_{|\alpha|=m}\big<\partial^\alpha(
u\partial_x\eta
+\widetilde{v}\partial_y\eta), \partial^\alpha \eta\big>\Big|
+\Big|\sum_{m=2}^{\infty}\frac{\rho(t)^{2(m-2)}}{((m-2)!)^{\frac{2}{\gamma}}}\sum_{|\alpha|=m}\big<\beta
\frac{\partial^\alpha\partial_xp}{\varepsilon^\gamma}, \partial^\alpha  \eta\big>\Big|\\
&&\quad+\Big|\sum_{m=2}^{\infty}\frac{\rho(t)^{2(m-2)}}{((m-2)!)^{\frac{2}{\gamma}}}\sum_{|\alpha|=m}\big<\partial^\alpha
(\partial_yR_1-\partial_xR_2+\varepsilon^2fe^{-y}\partial_yu^a+\varepsilon^2\partial_xfe^{-y}\partial_yv^a-\frac{\beta R_1}{\varepsilon^\gamma}+h), \partial^\alpha \eta\big>\Big|\\
&&=\sum_{i=1}^5M_i.
\eeno

{\bf Step 1. \textcolor[rgb]{0.00,0.00,0.00}{Estimate of $M_1+M_3.$}}
Using Lemma \ref{e:x product estimate'} and Lemma \ref{e:y product estimate}, we get
\beno
&&M_1+M_3\\
&&\leq C_0\|\eta\|^2_{Y^2\cap X^2}+C_0\|u\|_{X^3}^{\frac12}\| \partial_yu\|_{X^2}^{\frac12}\|\eta\|^2_{X^2\cap Y^2} +C_0\|u\|_{X^3}\| \eta\|_{X^2}^{\frac32}\| \partial_y\eta\|_{X^2}^{\frac12}
\\
&&\quad+C_0(\delta)\|(\widetilde{v},\partial_y\widetilde{v})\|^{\frac12}_{\overline{H}^3}
\|(\partial_y\widetilde{v},\partial_{yy}\widetilde{v})\|^{\frac12}_{\overline{H}^3}
\| \eta\|^2_{X^2\cap Y^2}\\
&&\quad+C_0\| \widetilde{v}\|^{\frac12}_{X^2}\|\partial_y\widetilde{v}\|_{X^2}^{\frac12}
\|\partial_y\eta\|_{X^2}\|\eta\|_{X^2}+C_0(\delta)\| \eta\|^{\frac32}_{X^2}\| \partial_y\eta\|_{X^2}^{\frac12}\|(\widetilde{v},\partial_y\widetilde{v})\|_{X^2}
\eeno

{\bf Step 2. Estimate of $M_2$.} On the one hand,
recalling the definition of $\eta^a$ (\ref{eq:relationship of vorticity}), and we have
\beno
&&\sum_{m=2}^{\infty}\frac{\rho(t)^{2(m-2)}}{((m-2)!)^{\frac{2}{\gamma}}}\sum_{|\alpha|=m}\langle \partial^\alpha(u\partial_x\eta^a),\partial^\alpha \eta\rangle\\
&&=\sum_{m=2}^{\infty}\frac{\rho(t)^{2(m-2)}}{((m-2)!)^{\frac{2}{\gamma}}}\sum_{|\alpha|=m}\langle \partial^\alpha(u\partial_xw^a),\partial^\alpha \eta\rangle-\beta\sum_{m=2}^{\infty}\frac{\rho(t)^{2(m-2)}}{((m-2)!)^{\frac{2}{\gamma}}}\sum_{|\alpha|=m}\Big< \partial^\alpha\Big(\frac{u\partial_xu^a}{\varepsilon^\gamma}\Big),\partial^\alpha \eta\Big>\\
&&=M_{21}+M_{22}.
\eeno
Applying Lemma \ref{e:x sublinear product estimate} with $k=2$ and (\ref{e:uniform boundness for approximate solution}) to the term $M_{22}$, we have
\beno
|M_{22}|\leq C_0\varepsilon^{-2\gamma}\|u\|^2_{ X^2}+C_0\|\eta\|^2_{X^2\cap Y^2}.
\eeno
Meanwhile, a straightforward computation yields
\beno
|M_{21}|\leq \sum_{m=2}^{\infty}\frac{\rho(t)^{2(m-2)}}{((m-2)!)^{\frac{2}{\gamma}}}\sum_{|\alpha|=m}\left|\Big< \partial^\alpha\Big(\frac{u\partial_{xz}u^{p}}{\varepsilon^{\gamma}}+u\partial_{xy}u^{e}-u\partial_{xx}v^a\Big),\partial^\alpha \eta\Big>\right|,
\eeno
where we used $u^a=u^{e}+\varepsilon^{1-\gamma} u^{p}.$
Thus, the same argument as $M_{22}$ implies that
\beno
|M_{21}|\leq C_0\varepsilon^{-2\gamma}\|u\|^2_{X^2}+C_0\|\eta\|^2_{X^2\cap Y^2}.
\eeno
Summing up $M_{21}$ and $M_{22}$, we arrive at
\beno
\left|\sum_{m=2}^{\infty}\frac{\rho(t)^{2(m-2)}}{((m-2)!)^{\frac{2}{\gamma}}}\sum_{|\alpha|=m}\langle \partial^\alpha(u\partial_x\eta^a),\partial^\alpha \eta\rangle\right|\leq C_0\varepsilon^{-2\gamma}\|u\|^2_{X^2}+C_0\|\eta\|^2_{X^2\cap Y^2}.
\eeno

On the other hand, we have
\ben\label{e:linear trouble term}
&&\sum_{m=2}^{\infty}\frac{\rho(t)^{2(m-2)}}{((m-2)!)^{\frac{2}{\gamma}}}\sum_{|\alpha|=m}\langle \partial^\alpha(\tilde{v}\partial_y\eta^a),\partial^\alpha v\rangle\nonumber\\
&&=\sum_{m=2}^{\infty}\frac{\rho(t)^{2(m-2)}}{((m-2)!)^{\frac{2}{\gamma}}}\sum_{|\alpha|=m}\langle \partial^\alpha(\tilde{v}\partial_yw^a),\partial^\alpha v\rangle-\beta\sum_{m=2}^{\infty}\frac{\rho(t)^{2(m-2)}}{((m-2)!)^{\frac{2}{\gamma}}}\sum_{|\alpha|=m}\Big< \partial^\alpha\Big(\tilde{v}\frac{\partial_yu^a}{\varepsilon^\gamma}\Big),\partial^\alpha v\Big>,\nonumber\\
\een
Since the most singular term in (\ref{e:linear trouble term}) is
\beno
\sum_{m=2}^{\infty}\frac{\rho(t)^{2(m-2)}}{((m-2)!)^{\frac{2}{\gamma}}}\sum_{|\alpha|=m}\varepsilon^{1-\gamma}\langle \partial^\alpha(\tilde{v}\partial_{yy}(u^{p})),\partial^\alpha u\rangle,
\eeno
and the others are similar, hence we only estimates this term.
Note that $\partial_{yy}(u^{p})=\varepsilon^{-1}\partial_y(\partial_zu^{p})$, and by (\ref{e:uniform boundness for approximate solution}),
$\|\widetilde{Z}^k\partial_zu^{p}\|_{\infty}$ is bounded for $k\in N$. Hence,
by Lemma \ref{e:singular term estimate}, we obtain
\beno
&&\sum_{m=2}^{\infty}\frac{\rho(t)^{2(m-2)}}{((m-2)!)^{\frac{2}{\gamma}}}\sum_{|\alpha|=m}\varepsilon^{1-\gamma}\langle \partial^\alpha(\tilde{v}\partial_{yy}(u^{p})),\partial^\alpha \eta\rangle\\
&&\leq C_0\|\eta\|^2_{X^2\cap Y^2}+C_0(\varepsilon^{-2}\|U\|^2_{X^2\cap Y^2} +\varepsilon^2).
\eeno
Consequently,
\beno
M_2\leq C_0\|\eta\|^2_{X^2\cap Y^2}+C_0(\varepsilon^{-2}\|U\|^2_{X^2\cap Y^2}+\varepsilon^2).
\eeno

{\bf Step 3. Estimate of $M_4$.}
It is easy to get
\beno
M_4\leq C_0\delta \varepsilon^{-2\gamma}\|\partial_xp\|^2_{X^2}
+C_0(\delta)\|\eta\|^2_{X^2}.
\eeno
Thus, by Lemma \ref{lem:pressure estimate}, we get
\beno
M_4&\leq& C_0(\delta)\|\eta\|^2_{X^2}
+C_0\delta\varepsilon^{-2\gamma}\big(\|U\|^2_{\overline{H}^4}+\|\nabla U\|^2_{\overline{H}^4}+1\big)\Big(\big\|U\big\|^2_{\overline{H}^4}+\varepsilon^4\Big)\\[3pt]
&&+C_0\delta\varepsilon^{2(2-2\gamma)}\big\|\nabla U\big\|^2_{\overline{H}^3}+C_0\delta\varepsilon^{2(2-2\gamma)}\big(\|u\|^2_{Y^3}+\|\partial_yu\|^2_{Y^2}\big) \nonumber\\[3pt]
&&+C_0\delta\varepsilon^{-2\gamma}(\|U\|_{X^3}\|\partial_yU\|_{X^2}+1)\big(\varepsilon^4+\|U\|^2_{X^3\cap Y^3}\big).
\eeno

{\bf Step 4. Estimate of $M_5$.}
Noting that $\gamma\leq 1$,  by (\ref{e:uniform boundness for approximate solution})-(\ref{eq:boundness of f}), it is easy to obtain
\beno
M_5
\leq C_0\varepsilon^2+C_0\|\eta\|^2_{X^2}+C_0\|U\|^2_{X^2}.
\eeno

{\bf Step 5. Estimate of dissipative term.}
Noting that $\eta$ vanishes on the boundary and by integrating by parts,  we have
\beno
&&-\sum_{m=2}^{\infty}\frac{\rho(t)^{2(m-2)}}{((m-2)!)^{\frac{2}{\gamma}}}\sum_{|\alpha|=m}\big<\partial^\alpha \triangle \eta, \partial^\alpha  \eta\big>\\
&\geq&\Big(\big\|\nabla \eta\big\|^2_{X^2}-\|\nabla\eta\|_2^2\Big)-
\Big|\sum_{m=2}^{\infty}\frac{\rho(t)^{2(m-2)}}{((m-2)!)^{\frac{2}{\gamma}}}\sum_{|\alpha|=m}\big<[\partial^\alpha ,\partial_y] \partial_y\eta, \partial^\alpha \eta\big>\Big|\\
&&-
\Big|\sum_{m=2}^{\infty}\frac{\rho(t)^{2(m-2)}}{((m-2)!)^{\frac{2}{\gamma}}}\sum_{|\alpha|=m}\big<\partial^\alpha\partial_y\eta, [\partial^\alpha,\partial_y]  \eta\big>\Big|
\geq (1-C_0\delta)\Big(\big\|\nabla \eta\big\|^2_{X^2}-\|\nabla\eta\|_2^2\Big).
\eeno
Thus, fixed $\delta$ small, we arrive at
\beno
-\sum_{m=2}^{\infty}\frac{\rho(t)^{2(m-2)}}{((m-2)!)^{\frac{2}{\gamma}}}\sum_{|\alpha|=m}\big<\partial^\alpha \triangle \eta, \partial^\alpha  \eta\big>\geq \frac{1}{2}\Big(\big\|\nabla \eta\big\|^2_{X^2}-\|\nabla\eta\|_2^2\Big).
\eeno

Collecting the estimates in Step 1-Step 5, the proof is completed.
\endproof

\subsection{Vorticity estimate in Sobolev space}

\begin{Proposition}[Vorticity estimate in $L^2$]\label{prop:Vorticity estimate in Sobolev space}
There exist $\delta>0$ and $\varepsilon_0>0$ such that for any $\varepsilon\in (0,\varepsilon_0)$, we have
\begin{align}
\frac{1}{2}\frac{d}{dt}\big\|\eta\big\|^2_2+\frac{\varepsilon^2}{2}\big\|\nabla \eta\big\|^2_2
\leq & C_0 \varepsilon^{-2}\big(\|U\|^2_{\overline{H}^4}+\|\nabla U\|^2_{\overline{H}^4}+1\big)\Big(\big\|U\big\|^2_{\overline{H}^4}+\varepsilon^4\Big) \nonumber\\
&+C_0\delta\beta^2\varepsilon^{4(1-\gamma)}\big\|\nabla U\big\|^2_{\overline{H}^3}+C_0(\delta)\big(1+\|U\|_{\overline{H}^4}+\big\|\nabla U\big\|_{\overline{H}^4}\big)\big\|\eta\big\|^2_{\overline{H}^3}.\nonumber
\end{align}
\end{Proposition}

\no{\bf Proof:}
Taking $L^2(R^2_+)$ inner product with $\eta$ in (\ref{eq:modified vorticity equation}), we arrive at
\begin{align}
&\frac{1}{2}\frac{d}{dt}\|\eta\|^2_2
+\varepsilon^2\|\nabla \eta\|_2^2\nonumber\\
\leq&\big|\big<u^{a}\partial_x\eta+\tilde{v}^a\partial_y\eta,\eta\big>\big|+\big|\big<u\partial_x\eta^a+\tilde{v}\partial_y\eta^a,\eta\big>\big|
+\big|\big<u\partial_x\eta+\tilde{v}\partial_y,\eta\big>\big|
+\beta\big|\big<\varepsilon^{-\gamma}\partial_x p,\eta\big>\big|\nonumber\\
&+\big|\big<\partial_yR_1-\partial_xR_2+\varepsilon^2fe^{-y}\partial_yu^a+\varepsilon^2\partial_xfe^{-y}\partial_yv^a-\beta \varepsilon^{-\gamma} R_1+h,\eta\big>\big|=\sum_{i=1}^5\tilde{N}_i.\nonumber
\end{align}

{\bf Step 1. Estimate of $\tilde{N}_1$.}
The same argument as $N_1$ gives
$$\tilde{N}_1\leq  C_0\varepsilon^2\big\|\eta\big\|^2_2.$$

{\bf Step 2. Estimate of $\tilde{N}_2$.}
Due to (\ref{eq:relationship of vorticity}), we have
\beno
\tilde{N}_2
&\leq& \Big|\Big<u\frac{\beta\partial_xu^a}{\varepsilon^\gamma},\eta\Big>\Big|+\big|\big<u\partial_xw_a,\eta\big>\big|
+\beta\Big|\Big<\tilde{v}\frac{\partial_yu^a}{\varepsilon^\gamma},\eta\Big>\Big|+\big|\big<\tilde{v}\partial_yw_a,\eta\big>\big|.
\eeno

In the above, one of the most singular terms is
$$\varepsilon^{1-\gamma}\big|\big<\tilde{v}\partial_{yy}u^{p},u\big>\big|.$$
We only estimate this term. The same argument as $N_2$ yields
\begin{align}
\varepsilon^{1-\gamma}\big|\big<\tilde{v}\partial_{yy}u^{p},u\big>\big|
=\varepsilon^{-\gamma}\big|\Big<\frac{v+\varepsilon^2f(t,x)e^{-y}}{\psi}Z\partial_zu^{p}),u\Big>\big|
\leq C_0\varepsilon^{-2\gamma}\big\|U\big\|^2_{\overline{H}^1}+C_0\varepsilon^2.\nonumber
\end{align}

{\bf Step 3. Estimate of $\tilde{N}_3$.}
The same argument as $N_3$ gives
\begin{align}
\tilde{N}_3\leq C_0\big(1+\|U\|_{\overline{H}^4}+\big\|\nabla U\big\|_{\overline{H}^4}\big)\big\|\eta\big\|^2_{\overline{H}^3}.\nonumber
\end{align}

{\bf Step 4. Estimate of $\tilde{N}_4$ and $\tilde{N}_5$.}
It is direct to get
\begin{align}
\tilde{N}_4\leq C_0\delta \varepsilon^{-2\gamma}\|\partial_xp\|^2_2+C_0(\delta)\big\|\eta\big\|^2_2.\nonumber
\end{align}
By Lemma \ref{lem:low order pressure part}, we have
\begin{align}
\tilde{N}_4\leq C_0 \varepsilon^{-2\gamma}\big(\|U\|^2_{\overline{H}^4}+\|\nabla U\|^2_{\overline{H}^4}+1\big)\big(\big\|U\big\|^2_{\overline{H}^4}+\varepsilon^4\big)
+C_0\delta\varepsilon^{4(1-\gamma)}\big\|\nabla U\big\|^2_{\overline{H}^3}+C_0(\delta)\big\|\eta\big\|^2_2.\nonumber
\end{align}

By (\ref{e:uniform boundness for error}), we similarly have
$$\tilde{N}_5\leq C_0\big\|\eta\big\|^2_2+C_0\varepsilon^2. $$

Collecting these estimates in Step 1-Step 4, we complete the proof.
\endproof

\section{Construction of the approximate solution}
\setcounter{equation}{0}
In this section, we will construct the approximate solutions of the Navier-Stokes equations (\ref{eq:NS})
by using the matched asymptotic expansion method  for $\gamma=\frac12,1$, and the presentation of reminders $R_1,R_2$ are shown. The same argument can be used to construct approximate solution for any rational number $\gamma\in [0,1]$.

\subsection{The case of $\gamma=\frac12$}

In this subsection, we consider the case of  $\gamma=\frac12$, and compute the specific presentation of the reminders $R_1,R_2$. Moreover, the Assumptions (H) will be verified to be reasonable in this case.

We make the following formal expansions
\beno
&&u^{\varepsilon}(t,x,y)=\sum_{j\geq 0}\varepsilon^{\frac{j}{2}}u_e^{(j)}(t,x,y),\\
&&v^{\varepsilon}(t,x,y)=\sum_{j\geq 0}\varepsilon^{\frac{j}{2}}v_e^{(j)}(t,x,y),\\
&&p^{\varepsilon}(t,x,y)=\sum_{j\geq 0}\varepsilon^{\frac{j}{2}}p_e^{(j)}(t,x,y),
\eeno
away from the boundary and
\beno
&&u^{\varepsilon}(t,x,y)=\sum_{j\geq 0}\varepsilon^{\frac{j}{2}}[u_e^{(j)}(t,x,y)+u_p^{(j)}(t,x,z)],\\
&&v^{\varepsilon}(t,x,y)=\sum_{j\geq 0}\varepsilon^{\frac{j}{2}}[v_e^{(j)}(t,x,y)+v_p^{(j)}(t,x,z)],\\
&&p^{\varepsilon}(t,x,y)=\sum_{j\geq 0}\varepsilon^{\frac{j}{2}}[p_e^{(j)}(t,x,y)+p_p^{(j)}(t,x,z)]
\eeno
near the boundary with $z=\frac{y}{\varepsilon}$ and the following conclusions hold
\beno
&&u_e^{(i)}(t,x,y)=v_e^{(i)}(t,x,y)=p_e^{(i)}(t,x,y)=0,\quad i=1,2\\
&&u_p^{(0)}(t,x,y)=0,\quad v_p^{(i)}(t,x,y)=0,\quad i=0,1,2,\\
&&p_p^{(i)}(t,x,y)=0,\quad i=0,1,2,3,4.
\eeno

More generally, for $(u_e^{(j)},v_e^{(j)},p_e^{(j)})$ with $j\geq 0$ we have
\begin{align} \label{eq:Ej}(\rm Ej)\,\, \left\{
\begin{aligned}
&\partial_t u_e^{(j)}+\sum_{k=0}^j(u_e^{(k)}\partial_x+ v_e^{(k)}\partial_y) u_e^{(j-k)}+\partial_x p_e^{(j)}=\triangle u_e^{(j-4)},\\
&\partial_t v_e^{(j)}+\sum_{k=0}^j(u_e^{(k)}\partial_x+ v_e^{(k)}\partial_y) v_e^{(j-k)}+\partial_y p_e^{(j)}=\triangle v_e^{(j-4)},\\
&\partial_x u_e^{(j)}+\partial_y v_e^{(j)}=0,\\
&v_e^{(j)}|_{y=0}=-v_p^{(j)}|_{y=0},\\
&(u_e^{(j)},v_e^{(j)})|_{t=0}=(0,0).
\end{aligned}
\right. \end{align}
Here $u_e^{k}=v_e^{k}=0$ if $k<0.$
Concluding the above analysis, we have
 $(u_p^{(1)}, v_p^{(3)})$ satisfies the following \textcolor[rgb]{0.00,0.00,0.00}{linear Prandtl-type} equation
\begin{align} \label{eq:u-p 1}(\rm u_p^{(1)})\,\, \left\{
\begin{aligned}
&\partial_tu_p^{(1)}
+u_p^{(1)}\partial_x\overline{u_e^{(0)}}
+\overline{u_e^{(0)}}\partial_xu_p^{(1)}+z\overline{\partial_yv_e^{(1)}}\partial_zu_p^{(1)}
=\partial_z^2 u_p^{(1)},\\
&\partial_x u_p^{(1)}+\partial_z v_p^{(3)}=0,\\
&\partial_zu_p^{(1)}(t,x,0)=\beta u_e^{(0)}(t,x,0),\\
&u_p^{(1)}(t,x,\infty)=v_p^{(3)}(t,x,\infty)=0,\\
&u_p^{(1)}(0,x,z)=0.
\end{aligned}
\right. \end{align}

Similarly, $(u_p^{(2)}, v_p^{(4)})$ satisfies
\begin{align} \label{eq:u-p 2}(\rm u_p^{(2)})\,\, \left\{
\begin{aligned}
&\partial_tu_p^{(2)}+\overline{u_e^{(0)}}\partial_xu_p^{(2)}+u_p^{(2)}\partial_x\overline{u_e^{(0)}}
+u_p^{(1)}\partial_xu_p^{(1)}+(\overline{v_e^{(3)}}+v_p^{(3)})\partial_zu_p^{(1)}\\
&\quad +z\overline{\partial_yv_e^{(0)}}\partial_zu_p^{(2)}=\partial_z^2 u_p^{(2)},\\
&\partial_x u_p^{(2)}+\partial_z v_p^{(4)}=0,\\
&\partial_zu_p^{(2)}(t,x,0)=\beta(u_e^{(1)}(t,x,0)+u_p^{(1)}(t,x,0))-\partial_yu_e^{(0)}(t,x,0),\\
&u_p^{(2)}(t,x,\infty)=v_p^{(4)}(t,x,\infty)=0,\\
&u_p^{(2)}(0,x,z)=0.
\end{aligned}
\right. \end{align}

Moreover,  $(u_p^{(3)}, v_p^{(5)})$ satisfies
\begin{align} \label{eq:u-p 3}(\rm u_p^{(3)})\,\, \left\{
\begin{aligned}
&\partial_tu_p^{(3)}+\overline{u_e^{(0)}}\partial_xu_p^{(3)}+u_p^{(3)}\overline{\partial_xu_e^{(0)}}+u_p^{(1)}\partial_xu_p^{(2)}+u_p^{(2)}\partial_xu_p^{(1)}
+v_p^{(3)}\overline{\partial_yu_e^{(0)}}\\
&\quad +\sum_{k=3}^{4}(\overline{v_e^{(k)}}+v_p^{(k)})\partial_zu_p^{(5-k)}+z\overline{\partial_{xy}u_e^{(0)}} u_p^{(1)}+ z\overline{\partial_{y}u_e^{(0)}} \partial_xu_p^{(1)}\\
&\quad +z\overline{\partial_{y}v_e^{(0)}}\partial_{z}u_p^{(3)}
+\frac{z^2}{2}\overline{\partial_{yy}v_e^{(0)}}\partial_zu_p^{(1)}
=\partial_z^2 u_p^{(3)},\\
&\partial_x u_p^{(3)}+\partial_z v_p^{(5)}=0,\\
&\partial_zu_p^{(3)}(t,x,0)=\beta(u_e^{(2)}(t,x,0)+u_p^{(2)}(t,x,0))-\partial_yu_e^{(1)}(t,x,0),\\
&u_p^{(3)}(t,x,\infty)=v_p^{(5)}(t,x,\infty)=0,\\
&u_p^{(3)}(0,x,z)=0.
\end{aligned}
\right. \end{align}

Finally, the estimate of $p_p^{(5)}$ is needed, which can be solved by  the equation of $v_p^{(3)}$:
\ben\label{eq:v-p 3}
&&\partial_tv_p^{(3)}+\overline{u_e^{(0)}}\partial_xv_p^{(3)}
+u_p^{(3)}\overline{\partial_xv_e^{(0)}}+v_p^{(3)}\overline{\partial_yv_e^{(0)}}+z\overline{\partial_{xy}v_e^{(0)}}u_p^{(1)}
+z\overline{\partial_{y}v_e^{(0)}}\partial_{z}v_p^{(3)}
+\partial_zp_p^{(5)}= \partial_z^2 v_p^{(3)}\nonumber\\
\een

These equations can be solved by the following way
\beno
(u_e^{(0)},v_e^{(0)})\rightarrow (u_p^{(1)},v_p^{(3)})\rightarrow (u_e^{(3)},v_e^{(3)})\rightarrow (u_p^{(2)},v_p^{(4)})\rightarrow(u_e^{(4)},v_e^{(4)})\rightarrow (u_p^{(3)},v_p^{(5)})
\eeno

The approximate solutions are stated as follows
\beno
&&u^{a}(t,x,y)=\sum_{j=0,3}\varepsilon^{\frac{j}{2}}u_e^{(j)}(t,x,y)+\sum_{j=1}^3\varepsilon^{\frac{j}{2}}u_p^{(j)}(t,x,z)= u^e(t,x,y)+\varepsilon^{\frac12} u^p(t,x,z),\\
&&v^{a}(t,x,y)=\sum_{j=0,3}\varepsilon^{\frac{j}{2}}v_e^{(j)}(t,x,y)+\sum_{j=3}^5\varepsilon^{\frac{j}{2}}v_p^{(j)}(t,x,z)= v^e(t,x,y)+\varepsilon^{\frac32} v^p(t,x,z),\\[5pt]
&&p^{a}(t,x,y)=\sum_{j=0,3}\varepsilon^{\frac{j}{2}}p_e^{(j)}(t,x,y)+\varepsilon^{\frac{5}{2}}p_p^{(5)}(t,x,z).
\eeno

Let
\ben\label{d:difinition of f}
f(t,x)=\overline{v_p^{(4)}}+\varepsilon^{\frac{1}{2}}\overline{v_p^{(5)}},
\een
and
\ben\label{d:difinition of g0}
g_0(t,x)=\beta u_e^{(3)}(t,x,0)+\beta u_p^{(3)}(t,x,0)-\varepsilon^{\frac12}\partial_yu_e^{(3)}(t,x,0),
\een
then the approximate solution $(u^a,v^a,p^a)$ satisfies
\begin{align} \label{eq:approximate}\,\, \left\{
\begin{aligned}
&\partial_t u^a+ u^a\partial_x u^a+ (v^a-\varepsilon^2f(t,x)e^{-y})\partial_y u^a+\partial_x p^a-\varepsilon^2\triangle u^a=-R_1,\\
&\partial_t v^a+ u^a\partial_x v^a+ (v^a-\varepsilon^2f(t,x)e^{-y})\partial_y v^a+\partial_y p^a-\varepsilon^2\triangle v^a=-R_2,\\
&\partial_x u^a+\partial_y v^a=0,\\
&v^a(t,x,0)=\varepsilon^2f(t,x),\\
&\partial_y u^a(t,x,y)|_{y=0}=\beta\varepsilon^{-\gamma} u^a(t,x,0)-\varepsilon g_0(t,x),\\
&(u^a, v^a)(0,x,y)=(0,0).
\end{aligned}
\right. \end{align}

Especially, direct computations show that the reminders $R_1,R_2$ have the following representation:
\ben\label{error:1}
-R_1
&=&\varepsilon^{3}(u_e^{(3)}\partial_x+ v_e^{(3)}\partial_y) u_e^{(3)}-\varepsilon^2\triangle u^e-
\varepsilon^{\frac{5}{2}}\partial_{xx}u^p+\varepsilon^{\frac{5}{2}}\partial_xp_p^{(5)}\nonumber\\[5pt]
&&-\varepsilon^2f e^{-y}\partial_y u^e- \varepsilon^2 e^{-y} \overline{v_p^{(5)}}\partial_z u^p- \varepsilon^2 e^{-y} \overline{v_p^{(4)}}\partial_z(u_p^{(2)}+\varepsilon^{\frac12}u_p^{(3)} ) \nonumber\\[5pt]
&&+\varepsilon^2\Big[u_e^{(3)}\partial_x u^p+\partial_xu_e^{(3)} u^p+\sum_{j=4}^6\varepsilon^{\frac{j-4}{2}}\sum_{k=1}^3u_p^{(k)}\partial_xu_p^{(j-k)}+
v_e^{(3)}\partial_z u_p^{(3)}+v_p^{(4)}\partial_y u^e\Big]\nonumber\\[5pt]
&&+\varepsilon^2\Big[\varepsilon^{\frac{1}{2}}v_p^{(5)}\partial_y u^e+\varepsilon v_p^{(3)}\partial_y u_e^{(3)}+
v_p^{(3)}\partial_z u_p^{(3)}+v_p^{(4)}\partial_z u_p^{(2)}+\varepsilon^{\frac{1}{2}}v_p^{(4)}\partial_z u_p^{(3)}+v_p^{(5)}\partial_z u^p\Big]\nonumber\\[5pt]
&&+\varepsilon^{\frac32}\Big[(u_e^{(0)}-\overline{u_e^{(0)}})\partial_xu_p^{(3)}+u_p^{(3)}(\partial_xu_e^{(0)}-\overline{\partial_xu_e^{(0)}})
+(v_e^{(3)}-\overline{v_e^{(3)}})\partial_zu_p^{(2)}\nonumber\\[5pt]
&&\quad\quad  +(v_e^{(4)}-\overline{v_e^{(4)}})\partial_zu_p^{(1)}+v_p^{(3)}(\partial_yu_e^{(0)}-\overline{\partial_yu_e^{(0)}})
+\overline{v_p^{(4)}}(1-e^{-y})\partial_zu_p^{(1)}\Big]\nonumber\\[5pt]
&&+\varepsilon\Big[(u_e^{(0)}-\overline{u_e^{(0)}})\partial_xu_p^{(2)}+u_p^{(2)}(\partial_xu_e^{(0)}-\overline{\partial_xu_e^{(0)}})
+(v_e^{(3)}-\overline{v_e^{(3)}})\partial_zu_p^{(1)}\Big]\nonumber\\[5pt]
&&+\varepsilon^{\frac12}\Big[(\partial_xu_e^{(0)}-\overline{\partial_xu_e^{(0)}}-\overline{y\partial_{xy}u_e^{(0)}})u_p^{(1)}
+(u_e^{(0)}-\overline{u_e^{(0)}}-\overline{y\partial_yu_e^{(0)}})\partial_xu_p^{(1)}\nonumber\\
&&\quad \quad +(v_e^{(0)}-y\overline{\partial_{y}v_e^{(0)}})\partial_{z}u_p^{(3)}\Big]
\nonumber\\[5pt]
&&+(v_e^{(0)}-y\overline{\partial_yv_e^{(0)}})\partial_zu_p^{(2)}
+\varepsilon^{-\frac{1}{2}}(v_e^{(0)}-y\overline{\partial_yv_e^{(0)}}-\frac{y^2}{2}\overline{\partial_{yy}v_e^{(0)}})\partial_zu_p^{(1)}
\een
and
\ben\label{error:2}
-R_2
&=& \varepsilon^{3}(u_e^{(3)}\partial_x+ v_e^{(3)}\partial_y) v_e^{(3)}-\varepsilon^2\triangle v^e-
\varepsilon^{3}\partial_{xx}v^p-
\varepsilon^{2}\partial_{zz}v_p^{(4)}- \varepsilon^{\frac{5}{2}}v_p^{(5)}-\varepsilon^2f(t,x)e^{-y}\partial_y v^e\nonumber\\[5pt]
&&+\varepsilon^{\frac{5}{2}}u_e^{(3)}\partial_xv^p+\varepsilon^{2}u_e^{(0)}\partial_xv_p^{(4)}+\varepsilon^{\frac{5}{2}}u_e^{(0)}\partial_xv_p^{(5)}+\varepsilon^{2}\partial_xv_e^{(3)}u^p
\nonumber\\[5pt]
&&+\varepsilon^{\frac{5}{2}}\partial_yv_e^{(3)}v^p+ (\varepsilon^{2}v_p^{(4)}+\varepsilon^{\frac{5}{2}}v_p^{(5)})\partial_yv_e^{(0)}\nonumber\\[5pt]
&&+\varepsilon^{\frac32}\Big[(u_e^{(0)}-\overline{u_e^{(0)}})\partial_xv_p^{(3)}+\partial_xv_e^{(0)}u_p^{(3)}+u^p\partial_xv^p+v_e^{(3)}\partial_zv^p
+v_e^{(0)}\partial_zv_p^{(5)}\nonumber\\[5pt]
&&\quad \quad +v_p^{(3)}(\partial_yv_e^{(0)}-\overline{\partial_yv_e^{(0)}})\Big]\nonumber\\
&&+\varepsilon \Big[\partial_xv_e^{(0)}u_p^{(2)}+v_e^{(0)}\partial_zv_p^{(4)}+ v^p\partial_zv^p-f(t,x)e^{-y}\partial_z v^p\Big]\nonumber\\
&&+\varepsilon^{\frac{1}{2}}\Big[(\partial_xv_e^{(0)}-y\overline{\partial_{xy}v_e^{(0)}})u_p^{(1)} + (v_e^{(0)}- y\overline{\partial_{y}v_e^{(0)}})\partial_zv_p^{(3)} \Big].
\een

In fact, the approximate solutions and the remainders $R_1, R_2$ satisfy the Assumption $(H)$ if the initial data $(u_0,v_0)\in X_e^{20}$. More precisely, we have the following three lemmas, whose proofs will be presented in the next section.

\begin{Lemma}\label{e:uniform boundness for approximate solution 12}
Let $\tilde{\partial}=\partial_x^{\alpha_1}\tilde{Z}^{\alpha_2}$ with $\tilde{Z}^k=(\delta z)^k\partial_z^k$ for $k\in {\mathbb N}$. There exists $T_a> 0$ such that for any $t\in [0,T_a]$, there holds
\beno
&&\sum_{m=3}^{\infty}\frac{\rho(t)^{2(m-3)}}{((m-3)!)^{\frac{2}{\gamma}}}\sum_{m-3\leq |\alpha|\leq m+6}\big\|\partial_{x,y}^\alpha(u_e^{(j)},v_e^{(j)})\big\|^2_2\leq C_0,\quad j=0,3;\\
&&\sum_{s=0}^2\sum_{m=3}^{\infty}\frac{\rho(t)^{2(m-3)}}{((m-3)!)^{\frac{2}{\gamma}}}\sum_{m-3\leq |\alpha|\leq m+6}\|e^{z^2}\tilde{\partial}_{x,z}^\alpha\partial^s_z(u_p^{(j)},v_p^{(j+2)})\big\|_{2}\leq C_0,\quad j=1,2,3.
\eeno
\end{Lemma}

\begin{Lemma}\label{e:uniform boundness for error 12}
There exists $T_a> 0$ such that for any $t\in [0,T_a]$, there holds
\beno
\big\|(R_1, R_2)\big\|^2_{X^3}\leq C_0\varepsilon^4,\quad \big\|\nabla(R_1, R_2)\big\|^2_{X^2}\leq C_0\varepsilon^2.
\eeno
\end{Lemma}

\begin{Lemma}\label{e:Formulation and Uniform boundness of $f, g_0$}
There exist universal constant $C_0$  and $\overline{f}(t,x)$ such that
$
f(t,x)=\partial_x\overline{f}
$
with
\beno
\|\partial_tf\|_{X^3_x}+\|f\|_{X^5_x}+\|\partial_t\bar{f}\|_{L^2_x}+\|g_0\|_{X^5_x}+\|\partial_tg_0\|_{X^3_x}\leq C_0.
\eeno
\end{Lemma}

\subsection{The case of $\gamma=1$}

In this subsection, we construct approximate solutions, derive the equations of approximate solutions and compute the remainders when $\gamma=1$.

Approximate solutions are constructed similarly as above, but more simple, and the process can be found in \cite{WWZ} beside the boundary conditions of Prandtl equation. We only give a derivation of boundary condition for Prandtl equation.\smallskip

Making the boundary layer(Prandtl layer) expansions
\beno
&&u^{\varepsilon}(t,x,y)=\sum_{j\geq 0}\varepsilon^{j}[u_e^{(j)}(t,x,y)+u_p^{(j)}(t,x,z)],\\
&&v^{\varepsilon}(t,x,y)=\sum_{j\geq 0}\varepsilon^{j}[v_e^{(j)}(t,x,y)+v_p^{(j)}(t,x,z)],\\
&&p^{\varepsilon}(t,x,y)=\sum_{j\geq 0}\varepsilon^{j}[p_e^{(j)}(t,x,y)+p_p^{(j)}(t,x,z)]
\eeno
and matched boundary condition requires that
\beno
u_p^{(i)}(t,x,z)\rightarrow 0, \ \ v_p^{(i)}(t,x,z)\rightarrow 0, \ \ p_p^{(i)}(t,x,z)\rightarrow 0,\ \ \ \text{as} \ z\rightarrow+\infty.
\eeno
While, the boundary condition of $(u^\varepsilon, v^\varepsilon)$ on $y=0$ requires that
\beno
&&\partial_zu_p^{(0)}(t,x,0)=\beta(u_e^{(0)}(t,x,0)+u_p^{(0)}(t,x,0)),\\
&&\partial_zu_p^{(i)}(t,x,0)=\beta(u_e^{(i)}(t,x,0)+u_p^{(i)}(t,x,0))-\partial_yu_e^{(i-1)}(t,x,0),i=1,2,...\\
&&v_e^{(i)}(t,x,0)=-v_p^{(i)}(t,x,0),\ i=0,1,\cdots.
\eeno

The same argument as in \cite{WWZ} gives
\begin{eqnarray}\label{e:prandtl equation}
\left \{
\begin {array}{ll}
\partial_tu_p^{(0)}-\partial_{zz}u_p^{(0)}+u_p^{(0)}\partial_xu_e^{(0)}(t,x,0)+\big(u_p^{(0)}+u_e^{(0)}(t,x,0)\big)\partial_xu_p^{(0)}\\
\qquad+\big(v_p^{(1)}+v_e^{(1)}(t,x,0)+ z\partial_yv_e^{(0)}(t,x,0)\big)\partial_zu_p^{(0)}=0,\\[3pt]
\partial_xu_p^{(0)}+\partial_zv_p^{(1)}=0,\\[3pt]
\lim\limits_{z\rightarrow +\infty}(u_p^{(0)}, v_p^{(1)})(t,x,z)=0,\\[3pt]
\partial_zu_p^{(0)}(t,x,0)=\beta(u_e^{(0)}(t,x,0)+u_p^{(0)}(t,x,0))\\
u_p^{(0)}(0,x,y)=0,\\[3pt]
\end{array}
\right.
\end{eqnarray}
\begin{Remark}
Let $\widetilde{u}_p^{(0)}(t,x,z)=u_p^{(0)}(t,x,z)+u_e^{(0)}(t,x,0)$, and
\beno
\widetilde{v}_p^{(1)}(t,x,z)=v_p^{(1)}(t,x,z)+v_e^{(1)}(t,x,0)+ z\partial_yv_e^{(0)}(t,x,0).
\eeno
Then by Bernoulli law,
\begin{align}
\partial_tu_e^{(0)}(t,x,0)+u_e^{(0)}(t,x,0)\partial_xu_e^{(0)}(t,x,0)+\partial_x p_e^{(0)}(t,x,0)=0,\nonumber
\end{align}
we arrive at
\begin{eqnarray}
\left \{
\begin{array}{lll}
&\partial_t\widetilde{u}_p^{(0)}-\partial_{zz}\widetilde{u}_p^{(0)}+\widetilde{u}_p^{(0)}\partial_x\widetilde{u}_p^{(0)}+\widetilde{v}_p^{(1)}\partial_z\widetilde{u}_p^{(0)}
+\partial_x p_e^{(0)}(t,x,0)=0,\\[3pt]
&\partial_x\widetilde{u}_p^{(0)}+\partial_z\widetilde{v}_p^{(1)}=0,\\[3pt]
&\widetilde{u}_p^{(0)}(0,x,z)=u_e^{(0)}(0,x,0),\\[3pt]
&\lim\limits_{z\rightarrow +\infty}\widetilde{u}_p^{(0)}(t,x,z)=u_e^{(0)}(t,x,0),\\[3pt]
&\partial_z\widetilde{u}_p^{(0)}(t,x,0)=\beta\widetilde{u}_p^{(0)}(t,x,0),\ \ \widetilde{v}_p^{(1)}(t,x,0)=0,
\end{array}
\right.
\end{eqnarray}
which is a nonlinear Prandtl equation with Robin boundary condition.
\end{Remark}

Similarly, we can obtain the equation for $(u_p^{(1)},v_p^{(2)})$.
\begin{Remark}
These equations can be solved in the following way
$$(u_e^{(0)},v_e^{(0)})\rightarrow (u_p^{(0)},v_p^{(1)})\rightarrow (u_e^{(1)},v_e^{(1)})\rightarrow (u_p^{(1)},v_p^{(2)}). $$
\end{Remark}

Now let us define  the approximate solutions ($u^a,v^a,p^a$) as following:
\begin{align*}
&u^a(t,x,y):=\sum\limits_{i=0}^{1}\varepsilon^iu_e^{(i)}(t,x,y)+\sum\limits_{i=0}^{1}\varepsilon^iu_p^{(i)}\Big(t,x,\frac{y}{\varepsilon}\Big),\\&
v^a(t,x,y):=\sum\limits_{i=0}^{1}\varepsilon^iv_e^{(i)}(t,x,y)+\sum\limits_{i=1}^{2}\varepsilon^iv_p^{(i)}\Big(t,x,\frac{y}{\varepsilon}\Big),\\&
p^a(t,x,y):=\sum\limits_{i=0}^{1}\varepsilon^ip_e^{(i)}(t,x,y)+\varepsilon^2p_p^{(2)}\Big(t,x,\frac{y}{\varepsilon}\Big).
\end{align*}

Set
$$f(t,x):=\int_{0}^\infty\partial_xu_p^{(1)}(t,x,z)dz,\quad g_0(t,x)=-\partial_yu_e^{(1)}(t,x,0).$$
A straightforward computation gives
that the approximate solution $(u^a,v^a,p^a)$ satisfies
\begin{align}
\left\{
\begin{array}{lll}
\partial_t u^a+ u^a\partial_x u^a+ (v^a-\varepsilon^2f(t,x)e^{-y})\partial_y u^a+\partial_x p^a-\varepsilon^2\triangle u^a=-R_1,\\[5pt]
\partial_t v^a+ u^a\partial_x v^a+ (v^a-\varepsilon^2f(t,x)e^{-y})\partial_y v^a+\partial_y p^a-\varepsilon^2\triangle v^a=-R_2,\\[5pt]
\partial_x u^a+\partial_y v^a=0,\\[5pt]
(u^a, v^a)(0,x,y)=(0,0),\\[5pt]
v^a(t,x,0)=\varepsilon^2f(t,x),\\[5pt]
\partial_y u^a(t,x,y)|_{y=0}=\beta\varepsilon^{-1} u^a(t,x,0)-\varepsilon g_0(t,x).
\end{array}
\right.\nonumber
\end{align}
where the reminders $R_1,R_2$ has the same form as in \cite{WWZ}.

\begin{Remark}
i)In this case, the norms in the Assumptions (H) are just the analytical norms, and the above approximate solutions can also be verified to satisfy the Assumptions (H) with the initial analytic data similar to the case $\gamma=\frac12$. \\
ii)Note that at this time the equation of $u_p^{(0)}$ is a nonlinear Prandtl equation with Robin boundary condition, which can be solved in the analytic setting(see also \cite{DJ}), and we omit the proof of local well-posedness result.\\
iii)Since the nonlinear Prandtl equations with Robin boundary condition occur in this case(the local well-posedness in analytic setting was established in \cite{DJ}), the vanishing viscosity limit can only be verified in the analytic setting which is similar as in previous papers of dealing with no-slip conditions \cite{SC2,WWZ}. In this paper, we handle this case and the cases $\gamma<1$ in Gevrey class together, and give a unified proof.
\end{Remark}

\section{Appendix}

In this section, our goal is to prove the local well-posedness of approximate solutions and verify the reasonability of the Assumptions (H). We consider the case $\gamma=\frac12$. To be specific, we will prove the well-posedness of the Euler system (\ref{e:Euler equation}) and the \textcolor[rgb]{0.00,0.00,0.00}{linear Prandtl-type equation} (\ref{eq:u-p 1}), and
Lemma \ref{e:uniform boundness for approximate solution 12}-\ref{e:Formulation and Uniform boundness of $f, g_0$}.
Note that the proofs of well-posedness of the linearized equations (\ref{eq:Ej}) for $j=3$ and (\ref{eq:u-p 2}), (\ref{eq:u-p 3}) is similar, and we omit them.

\subsection{Well-posedness of the Euler system in Gevrey class}

First, let us introduce some semi-norms.
Set $\partial_{x,y}^\alpha:=\partial_x^{\alpha_1}\partial_y^{\alpha_2}, k\in {\mathbb N}$ and  define
\beno
&&\|f\|^2_{X_e^k}=\sum_{m=k}^\infty\frac{\rho_e(t)^{2(m-k)}}{((m-k)!)^{\frac{2}{\gamma}}}\sum_{|\alpha| =m}\big\|\partial_{x,y}^\alpha f\big\|^2_2+\|f\|_2^2,\\
&&\|f\|^2_{Y_e^k}=\sum_{m=k+1}^\infty\frac{\rho_e(t)^{2(m-k)}(m-k)}{((m-k)!)^{\frac{2}{\gamma}}}\sum_{|\alpha|= m}\big\|\partial_{x,y}^\alpha f\big\|^2_2+\|f\|_2^2,
\eeno
where $\rho_e(t)=2-\lambda_e t\geq 1$ and  $\lambda_e>0$, to be decided later. Moreover, $X^{k}_e(\R^2_+)$ states
\beno
X^k_e(\R^2_+)=\big\{u; \|u\|^2_{X_e^k}< \infty\big\}.
\eeno

\begin{Lemma}\label{e:x product estimate''}
Let $u\cdot n|_{y=0}=0$, and
\beno
A''&=&\sum_{m=19}^{\infty}\frac{\rho_e(t)^{2(m-19)}}{((m-19)!)^{\frac{2}{\gamma}}}\sum_{|\alpha|=m}\big|\langle \partial^\alpha(u\cdot\nabla v),\partial^\alpha v\rangle\big|.
\eeno
Then there holds
\beno
 A''\leq C_0\|(u, v)\|_{X_e^{19}}\big(\|v\|^2_{ Y_e^{19}}+\|(u,v)\|^2_{X_e^{19}}\big).
\eeno
\end{Lemma}

The proof is similar to Part $(a)$ of Lemma \ref{lem:x product estimate}, and we omit it here.

Since well-posedness of the Euler system in Sobolev space or Gevrey class is well-known, here we sketch it in our frame for completeness. The a priori energy estimates are stated as follows.
\begin{Proposition}\label{e:unifrom boundness for Euler approximate}
Let the initial data $(u_0,v_0)$ satisfy $\partial_xu_0+\partial_yv_0=0$ and $v_0(t,x,0)=0$. Moreover,
\beno
\|(u_0,v_0)\|_{X_e^{20}}^2\leq M<\infty.
\eeno
Then there exists $T_e>0$ such that the Euler system (\ref{e:Euler equation}) has a unique solution $U^e=(u^e,v^e)$ on $[0,T_e]$, which satisfies
\begin{align}
\sup_{0\leq t\leq T_e}\big(\|U^e\|^2_{X_e^{20}}+\|\partial_tU^e\|^2_{X_e^{19}}\big)\leq C_0.\nonumber
\end{align}
\end{Proposition}

\no{\bf Proof.} At first, we consider the vorticity equation of the system (\ref{e:Euler equation}):
\begin{align}\label{e:vorticity equation1}
\left\{
\begin{array}{ll}
\partial_t w^e+ U^e\cdot \nabla w^e=0,\\
U^e\cdot n|_{y=0}=0.
\end{array}
\right.
\end{align}
 Acting $\partial^\alpha_{x,y}$ on both sides of (\ref{e:vorticity equation1}), then taking $L^2$ inner product with $\frac{\rho_e(t)^{2(|\alpha|-19)}}{((|\alpha|-19)!)^{\frac{2}{\gamma}}}\partial_{x,y}^\alpha w^e$, integrating by parts and summing over \textcolor[rgb]{0.00,0.00,0.00}{$|\alpha|\geq 19$}, we arrive at
\beno
\frac12\frac{d}{dt}\|w^e\|^2_{X_e^{19}}+\lambda_e\|w^e\|^2_{Y_e^{19}}\leq C_0\sum_{m=19}^\infty\frac{\rho_e(t)^{2(m-19)}}{((m-19)!)^{\frac{2}{\gamma}}}\sum_{|\alpha|= m}\big|\big<\partial^\alpha_{x,y}(U^e\cdot\nabla w^e),\partial_{x,y}^\alpha w^e\big>\big|.
\eeno
\textcolor[rgb]{0.00,0.00,0.00}{By Lemma \ref{e:x product estimate''}}, the right side term can be controlled by
\beno
 C_0\|(U^e, w^e)\|_{X_e^{19}}\big(\|(U^e, w^e)\|^2_{X_e^{19}}+\| w^e\|^2_{Y_e^{19}}\big).
\eeno
Thus, we arrive at
\ben\label{eq:w e X19 estimate}
\frac12\frac{d}{dt}\|w^e\|^2_{X_e^{19}}+\lambda_e\|w^e\|^2_{Y_e^{19}}\leq C_0\|(U^e, w^e)\|_{X_e^{19}}\big(\|(U^e, w^e)\|^2_{X_e^{19}}+\| w^e\|^2_{Y_e^{19}}\big).
\een

Similarly, the energy estimate in Sobolev space $H^{18}$ can be obtained as follows
\ben\label{eq:w e H18 estimate}
\frac12\frac{d}{dt}\|w^e\|^2_{H^{18}}\leq C_0\|U^e\|_{H^{18}}\|w^e\|^2_{H^{18}}.
\een
To close the estimates, we have to recover the estimates of the velocity. Due to
\beno
\textcolor[rgb]{0.00,0.00,0.00}{\triangle v^e=-\partial_xw^e},\quad v^e|_{y=0}=0,
\eeno
it is easy to obtain
\beno
\|\nabla v^e\|^2_{X_e^{18}}\leq C_0\|w^e\|^2_{X_e^{18}}.
\eeno
Note that
\beno
\partial_yu^e=w^e+\partial_xv^e,\quad \partial_xu^e=-\partial_yv^e,
\eeno
and we also have
\beno
\|\nabla u^e\|^2_{X_e^{18}}\leq C_0\|w^e\|^2_{X_e^{18}}.
\eeno
Therefore,
\ben\label{eq:U e X19}
\|U^e\|^2_{X_e^{19}}\leq C_0\|w^e\|^2_{X_e^{18}}.
\een
Similarly, a straightforward computation gives
\ben\label{eq:U e H18}
\|U^e\|^2_{H^{20}}\leq C_0\Big(\|w^e\|^2_{H^{19}}+\|U^e\|^2_2\Big),\quad \|w^e\|^2_{X_e^{18}}\leq C_0\|w^e\|^2_{X_e^{19}}.
\een
On the other hand, it is direct to get from the velocity equation (\ref{e:Euler equation}) that
\ben\label{eq:U e L 2}
\frac12\frac{d}{dt}\|U^e\|^2_{L^2}= 0.
\een
Finally, from (\ref{eq:w e X19 estimate})-(\ref{eq:U e L 2}) we deduce that
\begin{align}\label{e:final energy estimate}
&\frac12\frac{d}{dt}\Big(\|w^e\|^2_{X_e^{19}}+\|U^e\|^2_{L^2}\Big)+\lambda_e\|w^e\|^2_{Y_e^{19}}\nonumber\\
&\leq C_0\Big(\|w^e\|_{X_e^{19}}+\|U^e\|_{L^2}\Big)\big\|w^e\|^2_{Y_e^{19}}
+C_0\Big(\|w^e\|_{X_e^{19}}+\|U^e\|_{L^2}\Big)^3.
\end{align}
Here we take $\lambda_e=4C_0M^{\frac12}$ and $T_e>0$ so that
\beno
\rho_e(t)\geq 1 \quad {\rm for } \quad t\in [0,T_e],\quad C_0T_eM\leq \frac12.
\eeno
Thus, from (\ref{e:final energy estimate}), a continuity  argument yields that
\beno
\sup_{0\leq t\leq T_e}\big(\|w^e\|^2_{X_e^{19}}+\|U^e\|^2_{L^2}\big)\leq C_0.
\eeno
Due to (\ref{eq:U e X19}) and (\ref{eq:U e H18}), there also holds
\beno
\sup_{0\leq t\leq T_e}\|U^e\|^2_{X_e^{20}}\leq C_0.
\eeno
Moreover, from the equations (\ref{e:vorticity equation1}), it is not difficult to get
\beno
\sup_{0\leq t\leq T_e}\|\partial_tw^e\|^2_{X_e^{18}}\leq C_0.
\eeno
Consequently, similar arguments imply that
$$\sup_{0\leq t\leq T_e}\|\partial_tU^e\|^2_{X_e^{19}}\leq C_0. $$

\subsection{Well-posedness of \textcolor[rgb]{0.00,0.00,0.00}{linear Prandtl-type equations} in Gevrey class}
To prove the well-posedness of \textcolor[rgb]{0.00,0.00,0.00}{linear Prandtl-type equations} in Gevrey class, we introduce the following semi-norms and functional spaces:
\beno
\|U^p\|^2_{X^k_p}&:=&\sum_{m=k}^\infty\frac{\rho_p(t)^{2(m-k)}}{((m-k)!)^{\frac{2}{\gamma}}}\sum_{|\alpha|= m}\int_{R^2_+}\Big| e^{\phi_p(t,z)}\widetilde{\partial}^\alpha U^p\Big|^2dxdz+\|U^p\|^2_{L^2_p},\\
\|U^p\|^2_{Y^k_p}&:=&\sum_{m=k+1}^\infty\frac{\rho_p(t)^{2(m-k)}(m-k)}{((m-k)!)^{\frac{2}{\gamma}}}\sum_{|\alpha|= m}\int_{R^2_+}\Big| e^{\phi_p(t,z)}\widetilde{\partial}^\alpha U^p\Big|^2dxdz\\
&&+\sum_{m=k}^\infty\frac{\rho_p(t)^{2(m-k)}}{((m-k)!)^{\frac{2}{\gamma}}}\sum_{|\alpha|= m}\int_{R^2_+}\Big| ze^{\phi_p(t,z)}\widetilde{\partial}^\alpha U^p\Big|^2dxdz+\|U^p\|^2_{\overline{L}^2_p},\\
\|U^p\|^2_{\overline{L}^2_p}&:=&\int_{R^2_+}\Big|z e^{\phi_p(t,z)} U^p\Big|^2dxdz,\quad
\|U^p\|^2_{L^2_p}:=\int_{R^2_+}\Big| e^{\phi_p(t,z)} U^p\Big|^2dxdz
\eeno
where
\beno
\rho_p(t)=2-\lambda_pt\geq 1,\quad \phi_p(t,z)=\rho_p(t)z^2,\quad \widetilde{\partial}^\alpha=\partial_x^{\alpha_1}\widetilde{Z}^{\alpha_2},\quad
\widetilde{Z}^k=(\delta z)^k\partial_z^k,
\eeno
and $\lambda_p>0$ to be decided later.
Moreover, we denote $X^k_p$ by
\beno
X^k_p=\big\{u:\|u\|^2_{X^k_p}< \infty\big\}.
\eeno

\begin{Lemma}\label{e:x product estimate'''}
Let $u$ be only the function of $x$, and
\beno
A'''&=&\sum_{m=17}^{\infty}\frac{\rho_p(t)^{2(m-17)}}{((m-17)!)^{\frac{2}{\gamma}}}\sum_{|\alpha|=m}\big|\langle \widetilde{\partial}^\alpha(u\partial_xv),e^{2\phi_p}\widetilde{\partial}^\alpha v\rangle\big|,\\
B'''&=&\sum_{m=17}^{\infty}\frac{\rho_p(t)^{2(m-17)}}{((m-17)!)^{\frac{2}{\gamma}}}\sum_{|\alpha|=m}\big|\langle \widetilde{\partial}^\alpha(uv),e^{2\phi_p}\widetilde{\partial}^\alpha v\rangle\big|.
\eeno
Then there hold
\beno
 A'''\leq
C_0\left(\sum_{m=17}^{\infty}\frac{\rho_p(t)^{2(m-17)}}{((m-17)!)^{\frac{2}{\gamma}}}\sum_{k=1}^{m+2}\| u\|_{H^{k}_x}^2\right)^{\frac12}\|v\|^2_{X_p^{17}\cap Y_p^{17}},
\eeno
and
\beno
 B'''\leq
C_0\left(\sum_{m=17}^{\infty}\frac{\rho_p(t)^{2(m-17)}}{((m-17)!)^{\frac{2}{\gamma}}}\sum_{k=1}^{m+2}\| u\|_{H^{k}_x}^2\right)^{\frac12}\|v\|^2_{X_p^{17}\cap Y_p^{17}}.
\eeno
\end{Lemma}

The estimate of $A'''$ is similar as Part $(b)$ of Lemma \ref{lem:x product estimate}, and  $B'''$ is similar to Lemma \ref{e:x sublinear product estimate}. We omit the details here.

\begin{Proposition}\label{e:well posedness of prandtl}
Let $(u^e,v^e)$ be given in Proposition \ref{e:unifrom boundness for Euler approximate}. There exists $T_p>0$ such that the system (\ref{eq:u-p 1}) has a unique solution $U^p=(u^p,v^p)$ on $[0,T_p]$, which satisfies
\beno
\sup_{0\leq t\leq T_p}\big(\|U^p\|^2_{X^{17}_p}+\|\partial_{zz}U^p\|^2_{X^{14}_p}+\|\partial_{t}U^p\|^2_{X^{14}_p}\big)\leq C_0.
\eeno
\end{Proposition}

\no{\bf Proof.}
As in Proposition \ref{e:unifrom boundness for Euler approximate}, we only derive a priori estimates. Recalling (\ref{eq:u-p 1}), $u^p=u_p^{(1)}$ satisfies
\begin{align}\label{e:lineaized prandtl equation}
\left\{
\begin{aligned}
&\partial_tu^p
+u^p\partial_x\overline{u^e}
+\overline{u^e}\partial_xu^p+z\overline{\partial_yv_e}\partial_zu^p
=\partial_z^2 u^p,\\[3pt]
&\partial_zu^p(t,x,0)=\beta u^e(t,x,0), \quad u^p(t,x,\infty)=0,\\
&u^p(0,x,z)=0.
\end{aligned}
\right.
\end{align}
Acting $\widetilde{\partial}^\alpha$ on both sides of (\ref{e:lineaized prandtl equation}), then taking $L^2$ inner product with $\frac{\rho_p(t)^{2(|\alpha|-17)}}{((|\alpha|-17)!)^{\frac{2}{\gamma}}}e^{2\phi_p(t,z)}\widetilde{\partial}^\alpha u^p$ and \textcolor[rgb]{0.00,0.00,0.00}{summing over $|\alpha|\geq 17$}, we arrive at
\ben\label{eq:u p energy estimate}
&&\frac12\frac{d}{dt}\big(\|u^p\|^2_{X^{17}_p}-\|u^p\|^2_{L^2_p}\big)+\lambda_p\Big(\big\|u^p\big\|^2_{Y^{17}_p}-\|U^p\|^2_{\overline{L}^2_p}\Big)\\
&&\quad-\sum_{m=17}^\infty\frac{\rho_p(t)^{2(m-17)}}{((m-17)!)^{\frac{2}{\gamma}}}\sum_{|\alpha|= m}\big<\widetilde{\partial}^\alpha\partial_{zz}u^p,e^{2\phi_p} \widetilde{\partial}^\alpha u^p\big>\nonumber\\
&&\leq C_0\sum_{m=17}^\infty\frac{\rho_p(t)^{2(m-17)}}{((m-17)!)^{\frac{2}{\gamma}}}\sum_{|\alpha|= m}\big|\big<\widetilde{\partial}^\alpha(u_p\partial_x\overline{u^e}
+\overline{u^e}\partial_xu^p+z\overline{\partial_yv_e}\partial_zu^p),e^{2\phi_p} \widetilde{\partial}^\alpha u^p\big>\big|
\een
Similar to the dissipative term as in Lemma \ref{e:Gevery derivate estimate on velocity}, we obtain
\ben\label{eq:u p dissipative estimate}
&&\sum_{m=17}^\infty\frac{\rho_p(t)^{2(m-17)}}{((m-17)!)^{\frac{2}{\gamma}}}\sum_{|\alpha|= m}\big<\widetilde{\partial}^\alpha\partial_{zz}u^p,e^{2\phi_p} \widetilde{\partial}^\alpha u^p\big>\nonumber\\
&&\geq \Big(\frac34-C_0\delta\Big)\Big(\|\partial_zu^p\|^2_{X^{17}_p}-\|\partial_zu^p\|^2_{L^2_p}\Big)-C_0\|U^p\|^2_{\overline{L}^2_p}\nonumber\\
&&\quad-\sum_{m=17}^\infty\frac{\rho_p(t)^{2(m-17)}}{((m-17)!)^{\frac{2}{\gamma}}}\int_{\partial R^2_+}\partial_x^m \partial_zu^p\partial_x^mu^pdx\nonumber\\
&&\geq \Big(\frac23-C_0\delta\Big)\Big(\|\partial_zu^p\|^2_{X^{17}_p}-\|\partial_zu^p\|^2_{L^2_p}\Big)
-C_0\|U^p\|^2_{\overline{L}^2_p}\nonumber\\
&&\quad-C_0\sum_{m=17}^\infty\frac{\rho_p(t)^{2(m-17)}}{((m-17)!)^{\frac{2}{\gamma}}}\|\partial^m_x \overline{u^e}\|_2^2.
\een
By Proposition \ref{e:unifrom boundness for Euler approximate}, Sobolev embedding inequality, and Lemma \ref{e:x product estimate'''}, we have
\ben\label{eq:estimate of ue up}
\Big|\sum_{m=17}^\infty\frac{\rho_p(t)^{2(m-17)}}{((m-17)!)^{\frac{2}{\gamma}}}\sum_{|\alpha|= m}\big<\widetilde{\partial}^\alpha(u_p\partial_x\overline{u^e}+\overline{u^e}\partial_xu^p
),e^{2\phi_p} \widetilde{\partial}^\alpha u^p\big>_2\Big|\leq C_0\|u^p\|^2_{X^{17}_p\cap Y^{17}_p}.
\een
Note that for any $m\in {\rm N}_+$,
\ben
[\widetilde{Z}^m,z]f=m\delta z \widetilde{Z}^{m-1}f.
\een
Using Proposition \ref{e:unifrom boundness for Euler approximate} and  Sobolev embedding inequality again, it follows  from a similar proof of Lemma \ref{e:x product estimate'''} that
\ben\label{eq:estimate of u p}
&&\Big|\sum_{m=17}^\infty\frac{\rho_p(t)^{2(m-17)}}{((m-17)!)^{\frac{2}{\gamma}}}\sum_{|\alpha|= m}\big<\widetilde{\partial}^\alpha(
z\overline{\partial_yv^e}\partial_zu^p),e^{2\phi_p} \widetilde{\partial}^\alpha u^p\big>\Big|\nonumber\\
&&\leq\Big|\sum_{m=17}^\infty\frac{\rho_p(t)^{2(m-17)}}{((m-17)!)^{\frac{2}{\gamma}}}\sum_{|\alpha|= m}\sum\limits_{\beta\leq \alpha}C_\alpha^\beta\big<\widetilde{\partial}^\beta(
\overline{\partial_yv^e})\widetilde{\partial}^{\alpha-\beta}(\partial_zu^p),e^{2\phi_p} z\widetilde{\partial}^\alpha u^p\big>\Big|\nonumber\\
&&\quad+\Big|\sum_{m=17}^\infty\frac{\rho_p(t)^{2(m-17)}}{((m-17)!)^{\frac{2}{\gamma}}}\sum_{|\alpha|= m}\sum\limits_{\beta\leq \alpha}C_\alpha^\beta\big<\widetilde{\partial}^\beta(
\overline{\partial_yv^e})[\widetilde{\partial}^{\alpha-\beta},z]\partial_zu^p),e^{2\phi_p} \widetilde{\partial}^\alpha u^p\big>\Big|\nonumber\\
&&\leq C_0\|u^p\|^2_{Y^{17}_p\cap X^{17}_p}+\frac{1}{10}\|\partial_zu^p\|^2_{X^{17}_p}.
\een

Taking $\delta$ small, by (\ref{eq:u p energy estimate})-(\ref{eq:estimate of u p}), we arrive at
\beno
&&\frac12\frac{d}{dt}\Big(\|u^p\|^2_{X^{17}_p}-\|u^p\|^2_{L^2_p}\Big)+(\lambda_p-C_0)\big\|u^p\big\|^2_{Y^{17}_p}+
\frac12\Big(\|\partial_zu^p\|^2_{X^{17}_p}-\|\partial_zu^p\|^2_{L^2_p}\Big)\\
&\leq& C_0\big(\|u^p\|^2_{X^{17}_p}+1\big)+\frac{1}{10}\|\partial_zu^p\|^2_{ H_p^{16}}.
\eeno
Moreover, it is easy to get
\beno
&&\frac12\frac{d}{dt}\|u^p\|^2_{L^2_p}+(\lambda_p-C_0)\big\|u^p\big\|^2_{\overline{L}^0_p}+\frac12\|\partial_zu^p\|^2_{L^2_p}\leq C_0\big(\|u^p\|^2_{L^2_p}+1\big).
\eeno
Consequently, we arrive at
\beno
\frac12\frac{d}{dt}\|u^p\|^2_{X^{17}_p}+(\lambda_p-C_0)\big\|u^p\big\|^2_{Y^{17}_p}
+\frac13\|\partial_zu^p\|^2_{X^{17}_p}
\leq C_0\big(\|u^p\|^2_{X^{17}_p}+1\big).
\eeno
Then, by taking $\lambda_p\geq C_0$, there exists $0<T_p\leq T_e$ such that
$$\sup_{0\leq t\leq T_p}\|u^p\|^2_{ X^{17}_p}\leq C_0.$$
Furthermore, due to
\beno
\partial_xu^p+\partial_zv^p=0,\quad \lim_{z\rightarrow \infty}v^p(t,x,z)=0,
\eeno
we have
\beno
\sup_{0\leq t\leq T_p}\|v^p\|^2_{ X^{16}_p}\leq C_0.
\eeno
Applying $\partial_t$ to the equation (\ref{e:lineaized prandtl equation}), similar energy estimates yield that
$$\sup_{0\leq t\leq T_p}\|\partial_t(u^p,v^p)\|^2_{ X^{14}_p}\leq C_0,$$
and using the equation (\ref{e:lineaized prandtl equation}) again, we get
$$\sup_{0\leq t\leq T_p}\|\partial_{zz}(u^p,v^p)\|^2_{ X^{14}_p}\leq C_0.$$
\endproof

\subsection{Verification of the Assumptions (H) for the case of $\gamma=\frac12$}

Making similar arguments as in Proposition \ref{e:unifrom boundness for Euler approximate} and Proposition \ref{e:well posedness of prandtl}, we can show that the linearized Euler equations (\ref{eq:Ej}) and \textcolor[rgb]{0.00,0.00,0.00}{linear Prandtl-type equations} (\ref{eq:u-p 2}), (\ref{eq:u-p 3}) with the associated initial and boundary condition can also be solved, in $X_e^{18}$, $X_p^{16}$ and $X_p^{15}$, respectively. In conclusion, we obtain the following estimates for approximate solutions:
\ben\label{eq:estimate of approximate solu}
&&\|(u_e^{(0)},v_e^{(0)})\|_{X_e^{19}}+\|\partial_t(u_e^{(0)},v_e^{(0)})\|_{X_e^{18}}\leq C_0,\nonumber\\
&&\|(u_e^{(3)},v_e^{(3)})\|_{X_e^{18}}+\|\partial_t(u_e^{(3)},v_e^{(3)})\|_{X_e^{17}}\leq C_0,\nonumber\\
&&\|(u_p^{(1)},v_p^{(3)})\|_{X_p^{16}}+\|\partial_t(u_p^{(1)},v_p^{(3)})\|_{X_p^{14}}+\|\partial_{zz}(u_p^{(1)},v_p^{(3)})\|_{X_p^{14}}\leq C_0,\nonumber\\
&&\|(u_p^{(2)},v_p^{(4)})\|_{X_p^{15}}+\|\partial_t(u_p^{(2)},v_p^{(4)})\|_{X_p^{13}}+\|\partial_{zz}(u_p^{(2)},v_p^{(4)})\|_{X_p^{13}}\leq C_0,\nonumber\\
&&\|(u_p^{(3)},v_p^{(5)})\|_{X_p^{14}}+\|\partial_t(u_p^{(3)},v_p^{(5)})\|_{X_p^{12}}+\|\partial_{zz}(u_p^{(3)},v_p^{(5)})\|_{X_p^{12}}\leq C_0.
\een
Then Lemma \ref{e:uniform boundness for approximate solution 12} follows easily.

Using (\ref{error:1}) and (\ref{error:2}), it is also easy to prove Lemma \ref{e:uniform boundness for error 12}.

Note that (\ref{d:difinition of f}) and
\beno
\overline{v_p^{(4)}}=-\int_0^\infty\partial_xu_p^{(2)}(t,x,z)dz,\quad \overline{v_p^{(5)}}=-\int_0^\infty\partial_xu_p^{(3)}(t,x,z)dz,
\eeno
and hence,
\beno
f(t,x)=-\partial_x\Big(\int_0^\infty u_p^{(2)}(t,x,z)dz+\varepsilon^{\frac12}\int_0^\infty u_p^{(3)}(t,x,z)dz\Big)=\partial_x\overline{f}(t,x).
\eeno
By (\ref{eq:estimate of approximate solu}), a direct computation yields that
\beno
\|f\|_{X^5_x}+\|\partial_tf\|_{X^3_x}+\|\partial_t\bar{f}\|_{L^2_x}(t)\leq C_0.
\eeno
Finally, due to (\ref{d:difinition of g0}) and (\ref{eq:estimate of approximate solu}), we also get
\beno
\|g_0\|_{X^5_x}+\|\partial_tg_0\|_{X^3_x}\leq C_0.
\eeno
Accordingly, Lemma \ref{e:Formulation and Uniform boundness of $f, g_0$} is proved. Thus, the Assumptions $(H)$ in Section 3 for $\gamma=\frac12$ is satisfied. \endproof

\bigskip

\noindent {\bf Acknowledgments.}
W. Wang was supported by NSFC under grant 11671067,
 the Fundamental Research Funds for the Central Universities and China Scholarship Council. Z. Zhang was partially supported by NSF of China under Grant 11425103.


\end{document}